\documentclass{llncs}[13pt]
\usepackage{mathrsfs}
\usepackage[utf8]{inputenc}
\usepackage[english]{babel}
\usepackage{multicol}
\usepackage{pdfpages}
\usepackage{wrapfig}
\usepackage{color}
\usepackage{bbm}
\usepackage{tabularx}
\usepackage{caption}
\usepackage{bm}
\usepackage{mathrsfs}
\usepackage{amssymb}
\usepackage{enumitem}
\usepackage{pdfpages}
\usepackage{float}
\usepackage{graphicx}
\usepackage{csquotes}
\usepackage{hyperref}
\usepackage{amsmath}
\usepackage[margin=1.30in]{geometry}
\usepackage{ulem}

\newtheorem{theo}{Theorem}

\newtheorem{lem}{Lemma}

\usepackage{ulem}
\usepackage{mathrsfs}
\usepackage[utf8]{inputenc}
\usepackage[english]{babel}
\usepackage{multicol}
\usepackage{pdfpages}
\usepackage{wrapfig}
\usepackage{color}
\usepackage{tabularx}
\usepackage{caption}
\usepackage{bm}
\usepackage{mathrsfs}
\usepackage{amssymb}
\usepackage{enumitem}
\usepackage{pdfpages}
\usepackage{float}
\usepackage{graphicx}
\usepackage{csquotes}
\usepackage{hyperref}
\usepackage{amsmath}
\usepackage{ulem}
\usepackage{lipsum}

\numberwithin{equation}{section}
\numberwithin{theorem}{section}
\numberwithin{cor}{section}
\numberwithin{lemma}{section}
\numberwithin{proposition}{section}
\numberwithin{cor}{section}
\numberwithin{eg}{section}
\numberwithin{examp}{section}

\begin{document}

\title{Higher Order Refinements by Bootstrap in Lasso and other Penalized Regression Methods}
\titlerunning{SOC in Post Model-selection Inference}  % abbreviated title (for running head)
\author{Debraj Das \inst{1},
Arindam Chatterjee\inst{2} \and S. N. Lahiri\inst{3}
}
\authorrunning{Das et al.} % abbreviated author list (for running head)


\institute{Indian Institute of Technology, Kanpur\\
\email{rajdas@iitk.ac.in}
\and
Indian Statistical Institute, Delhi\\
\email{cha@isid.ac.in}
\and
Washington University in St. Louis\\
\email{S.LAHIRI@WUSTL.EDU}
}

\maketitle              % typeset the title of the contribution

\begin{abstract}
Selection of important covariates and to drop the unimportant ones from a high-dimensional regression model is a long standing problem and hence have received lots of attention in the last two decades. After selecting the correct model, it is also important to properly estimate the existing parameters corresponding to important covariates. In this spirit, Fan and Li (2001) proposed Oracle property as a desired feature of a variable selection method. Oracle property has two parts; one is the variable selection consistency (VSC) and the other one is the asymptotic normality. Keeping VSC fixed and making the other part stronger, Fan and Lv (2008) introduced the strong oracle property. In this paper, we consider different penalized regression techniques  which are VSC and classify those based on oracle and strong oracle property. We show that both the residual and the perturbation bootstrap methods are second order correct for any penalized estimator irrespective of its class. Most interesting of all is the Lasso, introduced by Tibshirani (1996). Although Lasso is VSC, it is not asymptotically normal and hence fails to satisfy the oracle property.
\end{abstract}
\tableofcontents

\section{Introduction}
\large
Consider the multiple linear regression model
\begin{equation}\label{eqn:model}
y_{i} = \bm{x}'_{i}\bm{\beta}_n + \epsilon_{i}, \; \;\;\;\;     i = 1,\dots,n,
\end{equation}
where $y_1,\ldots,y_n$ are responses, $\epsilon_1,\ldots,\epsilon_n$ are independent and identically distributed (iid) random variables, $\bm{x}_1,\ldots,\bm{x}_n$ are known non-random design vectors, and $\bm{\beta}=(\beta_{1,n},\ldots, \beta_{p,n})$ is the $p$-dimensional vector of regression parameters. When the dimension $p$ is large, it is common to approach regression model (\ref{eqn:model}) with the assumption that the vector $\bm{\beta}_n$ is sparse, that is that the set $\mathcal{A}_n= \{j:\beta_{j,n}\neq 0\}$ has cardinality $p_0 = |\mathcal{A}_n|$ much smaller than $p$, meaning that only a few of the covariates are ``active''.
%The underlying regression method is expected to choose the appropriate model, that is to detect the set $\mathcal{A}_n$. The simplest way is to conduct least square estimation after performing  best subset selection using AIC or BIC. However best subset selection has two major drawbacks. First of all it is computationally inefficient since increase in the number of regression parameters, number of models to compare increases exponentially. It is also not stable meaning that small change in the data may result in completely different choice of the model.
A widely used approach to handle the underlying sparsity is to add a penalty term to the least square criterion function. The resulting estimator is called penalized regression estimator. A general definition of a penalized estimator $\hat{\bm{\beta}}_n$ is 
\begin{equation}\label{eqn:def}
\bm{\hat{\beta}}_n = \operatorname*{arg\,min}_{\bm{t}}\Bigg[\sum_{i=1}^{n}(y_i - \bm{x}'_i \bm{t})^2
+n\sum_{j=1}^{p}P_{\lambda_n, j}(|t_j|)\Bigg]
\end{equation}
where $\lambda_n>0$ is the penalty parameter and 
$P_{\lambda_n, j}(\cdot): [0,\infty)\rightarrow [0,\infty)$ is the penalty function corresponding to $j$th coefficient of $\bm{\beta}_n$. The penalty term $P_{\lambda_n, j}(\cdot)$ should impose some sparsity to the estimated model. However bringing sparsity to the model is not the only requirement from the perspective of choosing the right model. If $\{j:\hat{\beta}_{j,n}\neq 0\}=\mathcal{A}_n^{\mathsf{c}}$, that is if the penalized estimator $\bm{\hat{\beta}}_n$ chooses $\mathcal{A}_n^{\mathsf{c}}$ as the set of indices corresponding to important covariates, then also the resulting inference will be misleading. Therefore it is crucial from the perspective of valid inference that the penalized estimator detects the true set of covariates $\mathcal{A}_n$. In that case $\bm{\hat{\beta}}_n$ is said to be variable selection consistent or VSC. There are large number of penalized estimators in the literature that imposes sparsity to the regression model. Most well known is the  lasso, introduced by Tibshirani (1996). Lasso estimator is simply the minimizer of $l_1$-penalized least-square criterion function, that is $P_{\lambda_n, j}(|t|)=\dfrac{\lambda_n}{n}|t|$, for all $j \in \{1,\dots,p\}$.
%It is defined as the minimizer of $l_1$-penalized least-square criterion function, 
%\begin{equation}\label{eqn:lasso}
%\hat{\beta}_n = \operatorname*{arg\,min}_{t}\Bigg\{\sum_{i=1}^{n}(y_i - x^
%{\prime}_i t)^2
%+\lambda_n\sum_{j=1}^{p}|t_{j}|\Bigg\},
%\end{equation}
%where $\lambda_n>0$ is the penalty parameter.
Lasso is well suited to the sparse setting because of its property that it sets some regression coefficients exactly equal to 0 and hence it automatically leads to variable selection. Another attractive feature of the lasso is its computational feasibility [Fu (1998), Osborne et al. (2000)]. However Lasso is not in general variable selection consistent. The VSC property or its variants in lasso have been studied by many authors in different regularity conditions [cf. Zhao and Yu (2006), Meinshausen and B\"ulman (2006),  Wainwright (2009), Lahiri (2018)]. Recently, Lahiri (2018) found that the irrepresentable condition, along with an upper bound condition on $\lambda_n$ and some conditions on the design vectors, are necessary and sufficient for lasso to be VSC. Another interesting finding of Lahiri (2018) is that lasso can not be VSC and $\sqrt{n}$-consistent simultaneously. More precisely, Lahiri(2018) showed that when $p$ is fixed, for lasso to be VSC one needs $n^{-1/2}\lambda_n$ growing to $\infty$, which does not fall in the regime of Knight and Fu (2000). Knight and Fu (2000) established $\sqrt{n}$-consistency of lasso in fixed dimension assuming $n^{-1/2}\lambda_n \rightarrow \lambda_0 \in [0,\infty)$. Therefore when one is considering lasso, she has to consider VSC and $\sqrt{n}$-consistency separately under different conditions on $\lambda_n$. 

Under the set up of Knight and Fu (2000), that is when $n^{-1/2}\lambda_n \rightarrow \lambda_0 \in [0,\infty)$, lasso is $\sqrt{n}$-consistent, but has complicated asymptotic distribution [cf. Knight and Fu (2000), Wagener and Dette (2011), Camponovo (2015)]. An alternative approach beside using the asymptotic distribution is to look into the bootstrap distribution and utilize it to infer about $\bm{\beta}$. Chatterjee and Lahiri (2010) showed that the usual residual bootstrap fails when $\mathcal{A}_n^{\mathsf{c}}$ is not empty. Subsequently Chatterjee and Lahiri (2011) proposed a modification to the residual bootstrap and established its validity in approximating the distribution of the lasso estimator. Camponovo (2015) established the validity of paired bootstrap in lasso when the design is random. Recently, Das and Lahiri (2019) developed a perturbation bootstrap method in approximating the distribution of the lasso estimator and this approximation works irrespective of the nature of the design and even when the errors are independent, but may not be identically distributed. Although all of these results are helpful for the purpose of inference, the underlying lasso estimators lacks the basic necessity of VSC due to the assumption $n^{-1/2}\lambda_n \nrightarrow \infty$. Additionally, all of these results are in fixed $p$ setting and none is a uniform result. The problem of doing uniform asymptotic inference is still open under $n^{-1/2}\lambda_n \rightarrow \lambda_0 \in [0,\infty)$ even when $p$ is fixed.

In this paper, we consider lasso to be a VSC procedure and then investigate the difficulties that are arising in making valid inferences. We show that the Normal approximation is of no use in making inference; cf. Theorem \ref{thm:oracleIII}. We carefully construct residual and perturbation bootstrap procedures in making valid inferences, keeping lasso VSC. We establish Berry-Essen type result for both the bootstrap procedures for lasso uniformly over the collection of Borel measurable convex sets, even when $p$ grows with $n$; cf. Theorem \ref{thm:boot}. Therefore we can say that bootstrap is somewhat immune towards the effect of considering lasso a VSC procedure, where as $\sqrt{n}$-consistency or normal approximation are not. The reason behind this phenomenon is the substantial bias, incurred due to the assumption $n^{-1/2}\lambda_n\rightarrow \infty$ as $n\rightarrow \infty$. Both the bootstrap procedures can correctly mimic the bias and then correct it by using the fact that the collection of Borel measurable convex sets is closed under translation. Moreover, if we correctly define the studentized versions, then both the residual and perturbation bootstrap are second order correct even in increasing dimension. Second order correctness means that the error of distributional approximation is $o_p(n^{-1/2})$ uniformly over a class of sets, generally the class of all Borel measurable convex sets. See  Theorem \ref{thm:bootIII} for second order results of bootstrap in case of lasso.

%However, all of these bootstrap results are in fixed $p$ setting and none is a uniform result. only result that is available in lasso in increasing dimension is due to Goldsmith (2015). Under some conditions, he showed that if $\bm{\check{\beta}}_n$ is the lasso estimator, then  $\|\check{\bm{\beta}}_n-\bm{\beta}\|_1=O(\lambda_n/b_n^2)$ where $b_n=\inf \{t\geq 0: R^{-\alpha}(t)+r^{-\alpha}(t)=1/n\}$, $R^{-\alpha}(t)=\mathbf{P}(\epsilon_1>t)$, $r^{-\alpha}(t)=\mathbf{P}(-\epsilon_1>-t)$ and $\|\cdot\|_1$ is the $l_1$ norm; see Theorem 4.1.2 in Goldsmith (2015) for details. Although this is a nice result, it is again of little use to infer about $\bm{\beta}$ based on $\check{\bm{\beta}}_n$. 

Building on the ideas of lasso, other penalized methods are developed in the subsequent years. These developments mostly aim to rectify the issues of lasso, viz., to achieve VSC without irrepresentable condition and to make the penalized estimator asymptotically normal. Fan and Li (2001) discussed the guidelines to construct a penalty function to capture underlying sparsity and define the oracle property. A penalized regression method is said to satisfy the oracle property if it is VSC and if the estimator of the coefficients of important covariates perform asymptotically normal with same covariance matrix as the ordinary least square estimator (OLS). In other words the penalized regression method should work like an oracle who knows everything beforehand. In mathematical terms a penalized estimator $\hat{\bm{\beta}}_n=(\beta_{1,n},\dots, \beta_{p,n})^\prime$ is said to satisfy the Oracle property if $\hat{\bm{\beta}}_n$ has following two features:
\begin{enumerate}[label=(\Alph*)]
\item Variable Selection Consistency (VSC):

$\mathbf{P}(\hat{\mathcal{A}}_n=\mathcal{A}_n)\rightarrow 1$ where $\hat{\mathcal{A}}_n = \{1\leq j\leq p: \hat{\beta}_{j,n}\neq 0\}$. 
\item Asymptotic Normality with same precision as OLS:

$\sqrt{n}\bm{C}_{11,n}^{1/2}(\bm{\hat{\beta}}_{n}^{(1)}-{\bm{\beta}}_{n}^{(1)}) \xrightarrow{d}\mathbf{N}\big(\bm{0}, \sigma^2\bm{I}_{p_0}\big)$ when ${\mathcal{A}}_n=\{1,\ldots,p_0\}$.
\end{enumerate}
Here $\bm{C}_{11,n}$ is the upper left $p_0\times p_0$ submatrix of $\bm{C}_n=n^{-1}\sum_{i=1}^{n}\bm{x}_i\bm{x}_i^\prime$ and ``$\xrightarrow{d}$'' denotes convergence in distribution. As mentioned before, lasso does not have the second feature and hence it does not satisfy the oracle Property. Even the asymptotic normal approximation fails drastically for lasso; see theorem \ref{thm:oracleIII} for details.  Fan and Li (2001) developed the non-convex SCAD penalty following their own guidelines and showed that it satisfies the oracle property. A close relative to SCAD is MCP, developed by Zhang (2010), which is another non-convex penalized method with the oracle property. Lasso, being convex and at the boundary of oracle procedures, apparently it seems that one needs to have the penalty function non-convex with singularity at origin to achieve both the criteria of the oracle property. However the crux behind satisfying oracle Property lies in the construction of the SCAD penalty. The main aim of SCAD was to put more weights to the smaller coefficients than the larger ones in the penalty term, keeping the feature of singularity at origin. By exploring these features, Zou (2006) developed an oracle but convex penalized regression method, called adaptive lasso. Adaptive lasso is simply an weighted $l_1$-penalized regression with component specific penalty terms. The penalty function depends on some $\sqrt{n}$-consistent initial estimator. 

There is a deeper theoretical distinction between SCAD \& MCP and the adaptive lasso other than the nature of the optimization problems. To point out this distinction, we need to look into the strong oracle property, introduced by Fan and Lv (2008). A penalized estimator $\hat{\bm{\beta}}_n=(\beta_{1,n},\dots, \beta_{p,n})^\prime$ is said to satisfy the strong oracle property if $\hat{\bm{\beta}}_n$ has following features:
\begin{enumerate}[label=(\Alph*)]
\item Variable Selection Consistency (VSC):

$\mathbf{P}(\hat{\mathcal{A}}_n=\mathcal{A}_n)\rightarrow 1$ where $\hat{\mathcal{A}}_n = \{1\leq j\leq p: \hat{\beta}_{j,n}\neq 0\}$. 
\item Same as OLS on ${\mathcal{A}}_n$:

$\mathbf{P}\big{(}\hat{\bm{\beta}}_{{\mathcal{A}}_n}=\bar{\bm{\beta}}_{{\mathcal{A}}_n}\big{)}\rightarrow 1$, $\bar{\bm{\beta}}_{{\mathcal{A}}_n}$ being OLS of $\bm{\beta}_{{\mathcal{A}}_n}$ assuming $\bm{\beta}_{{\mathcal{A}}^c_n}=\bm{0}$.
\end{enumerate}
Here $\bm{\alpha}_{I}$ denotes the sub-vector of the vector $\bm{\alpha}$ with entries in the index set $I$. As a generalization of adaptive lasso, Zou and Li (2008) developed a class of weighted $l_1$-penalized estimators. Depending on the choice of the weights, the underlying one step estimator either satisfies the oracle property or satisfies the strong oracle property.
%SCAD and MCP both satisfy strong oracle property where as adaptive lasso and one step estimators do not.  Now let us consider the general definition of a penalized regression estimator to put all the aforementioned estimators in the same frame. A penalized estimator $\hat{\bm{\beta}}_n$ is defined as
%\begin{equation*}
%\bm{\hat{\beta}}_n = \operatorname*{arg\,min}_{\bm{t}}\Bigg[\sum_{i=1}^{n}(y_i - \bm{x}'_i \bm{t})^2
%+n\sum_{j=1}^{p}P_{\lambda_n, j}(|t_j|)\Bigg]
%\end{equation*}
%where $\lambda_n>0$ is the penalty parameter and 
%$P_{\lambda_n, j}(\cdot): [0,\infty)\rightarrow [0,\infty)$ is the penalty function corresponding to $j$th coefficient of $\bm{\beta}_n$.
Suppose $P^\prime_{\lambda_n, j}(\cdot)$ denotes the derivative of $P_{\lambda_n, j}(\cdot)$ in (\ref{eqn:def}). Then we can put all the aforementioned penalized estimators in the framework of (\ref{eqn:def}) in the following way:
\begin{enumerate}[label=(\arabic*)]
\item Lasso:
$P_{\lambda_n, j}^\prime(t)= \dfrac{\lambda_n}{n} , \; t>0, \; \text{for all}\; j$
\item SCAD:
$P^{\prime}_{\lambda_n, j}(t)= \lambda_n \mathbbm{1}(t\leq \lambda_n) +  \dfrac{(a\lambda_n - t)_+}{a-1}\mathbbm{1}(t > \lambda_n),\; t>0,\; a>2,\; \text{for all}\; j$
\item MCP: 
$P^{\prime}_{\lambda_n, j}(t)= (\lambda_n - a^{-1}t)_+,\; t>0,\; a>1,\; \text{for all}\; j$
\item Adaptive Lasso:
$P_{\lambda_n,j}^\prime(t)= \dfrac{\lambda_n}{ |\tilde{\beta_j}|^{\gamma}}, \; t>0, \; \gamma>0$
\item One Step Estimators:
$P_{\lambda_n,j}^\prime(t)= \tilde{P}^{\prime}_{\lambda_n}(|\tilde{\beta}_{j}|) ,\; t>0 $
\end{enumerate}
where $\mathbbm{1}(\cdot)$ is the indicator function and $\tilde{P}^{\prime}_{\lambda_n}(\cdot)$ is the derivative of some penalty function. One can consider $\tilde{P}_{\lambda_n}(\cdot)$ to be some non-convex penalty, like SCAD or MCP. $\tilde{\bm{\beta}}_n=(\tilde{\beta}_1,\dots, \tilde{\beta}_p)$ is some preliminary estimator, like OLS when $p\leq n$ and lasso or ridge estimator when $p>n$. Since oracle and strong oracle properties assume OLS as the benchmark, one can also use OLS for the model, selected by some variable selection procedure like lasso. The resulting estimator is called the Post-model selection OLS and was introduced by Belloni and  Chernozhukov (2013). Therefore, based on how close the penalized regression estimator is to the OLS, We can classify all the aforementioned penalized regression estimators in the following three classes:
\begin{enumerate}[label=(\Roman*)]
\item When the strong oracle property holds.
\item When the Oracle property holds, but not the strong oracle property.
\item When only VSC holds.
\end{enumerate}
SCAD, MCP and Post-model selection OLS fall in the class I, where as adaptive lasso falls in class II. One step estimator falls in Class I or II depending on choice of $\tilde{P}_{\lambda_n}(\cdot)$. If $\tilde{P}_{\lambda_n}(\cdot)$ is either SCAD or MCP, then one step estimator falls in class I. On the other hand if $\tilde{P}^{\prime}_{\lambda_n}(\cdot)=\lambda_n\tilde{P}^{\prime}(\cdot)$ with $\tilde{P}(\theta)=\theta^q$, $0<q<1$, or $\tilde{P}(\theta)=\log\theta$ then it falls in class II. As mentioned above, lasso falls in class III. For the first two classes of estimators, one can use the oracle normal approximation to make statistical inference. In addition to how good the normal approximation is, we also explore how good the bootstrap is. We show that both residual and perturbation bootstrap are second order correct for first two classes; see theorems \ref{thm:bootI} and \ref{thm:bootII} for details.

%\item \textcolor{blue}{$\gamma>0$}
%\item \textcolor{blue}{$\bm{\hat{\beta}}_n$} is the \textcolor{blue}{LSE} restricted on the set \textcolor{blue}{$\{(t_1,\dots,t_p):$}\\
%\textcolor{blue}{$\sum_{j=1}^{p}P_{\lambda_n}(|t_j|)\leq \check{\lambda}_n\}$} for some $\check{\lambda}_n$.

%Another important aspect of lasso is its computational feasibility in high dimensional regression problems [cf. Efron et al. (2004), Friedman et al. (2007), Fu (1998), Osborne et al. (2000)]. 

 To describe the findings related to the rate of convergence more elaborately, without loss of generality assume that $\mathcal{A}_n=\{j:\beta_{j,n}\neq 0\}=\{1,\dots,p_0\}$. $p$ \& $p_0$ can grow with $n$. Define $\bm{T}_n=\sqrt{n}\bm{D}_n(\hat{\bm{\beta}}_n - \bm{\beta})$ where $\bm{D}_n$ is a $q\times p$ matrix with $tr(\bm{D}_n\bm{D}_n^\prime)=O(1)$ ($q$ being fixed). Also define $\bm{\Sigma}_n=\bm{D}_n^{(1)}\bm{C}_{11,n}^{-1}\bm{D}_n^{(1)\prime}$ where $\bm{D}_n^{(1)}$ consists of first $p_0$  columns of $\bm{D}_n$. Under some regularity conditions, $\bm{\Sigma}_n$ is the asymptotic variance covariance matrix of $\bm{T}_n$. We will consider the following quantity to measure the error of oracle normal approximation:
\begin{align*}
\Delta_n = \sup_{\bm{B}\in \mathcal{C}_q}|\mathbf{P}(\bm{T}_n\in \bm{B})-\Phi(\bm{B}; \bm{\Sigma}_n)|
\end{align*}
where $\mathcal{C}_q$ is the collection of Borel measurable convex sets of $\mathcal{R}^q$ and $\Phi(\cdot; \bm{A})$ denotes the zero mean normal measure of the set $\bm{A}\in \mathcal{R}^q$. We show that $\Delta_n=O(n^{-1/2})$ for the estimators in class I, meaning that under strong oracle property the rate of oracle approximation is same as in case of OLS. However for the estimators in class II, $\Delta_n=O(n^{-1/2+c})$ for some constant $c\geq 0$ and the rate is attained with some $c>0$ under some additional conditions.% This is similar to the results obtained by Chatterjee and Lahiri (2013) (see Theorem 3.1 and 3.2 there) in case of adaptive Lasso.
Therefore when the oracle property holds but the strong oracle property fails, the error rate is worse than that of OLS, implying that the oracle based inference is not as accurate as in case of Class I. Moreover, for lasso, $\Delta_n\rightarrow 1$ as $n\rightarrow \infty$ implying that oracle normal approximation can not be used for the purpose of inference. The reason behind the decline of the rate of convergence as we move from class I to III is the increase in asymptotic order of the bias term. 

In this paper, we consider the residual and perturbation bootstrap approximations of the distribution of $\bm{T}_n$ and show that the error rate is $O_p(n^{-1/2})$ for each of the classes. Thus both the bootstrap methods correct for the effect of the bias in case of class II and III, whereas keeping the  error rate same as Oracle approximation for class I. % although $||\bm{T}_n||$ is not bounded (????????????**$\lambda_n/\sqrt{n}>K.p_0.\log n$ required)
 Moreover if we consider suitable studentized pivot based on $\bm{T}_n$ then the error rate of the bootstrap approximations improve uniformly to $o_p(n^{-1/2})$, that is both residual and perturbation bootstrap approximations are second order correct. These results show that the bootstrap methods enable, for example, the construction of confidence intervals for the nonzero regression coefficients for each of the classes of estimators even when dimension $p$ is increasing with $n$. Additionally the bootstrap inference is much more accurate than that based on the oracle normal approximation when atleast the oracle property holds.% with smaller coverage error than those based on the asymptotic distribution of Lasso.

We conclude this section with a brief literature review on bootstrap methods in penalized regression. In a series of papers Chatterjee and Lahiri showed usual residual bootstrap fails for approximating the distribution of Lasso and consequently developed a modified residual bootstrap [cf. Chatterjee and Lahiri (2010, 2011)] when the errors are iid and design vectors are non-random. In random design case, Camponovo (2015) developed a modified paired bootstrap in Lasso which also works in heteroscedastic regression. Recently Das and Lahiri (2019) developed a perturbation bootstrap method for Lasso which works irrespective of the nature of the design and when errors are heteroscedastic. Chatterjee and Lahiri (2013) established second order correctness of residual bootstrap in adaptive lasso. Perturbation bootstrap was introduced in adaptive lasso and SCAD by Minnier et al. (2011). Recently Das et al. (2018) showed that the perturbation bootstrap proposed by Minnier et al. (2011) fails to be second order correct and developed a modification which enables perturbation bootstrap to achieve second order correctness in adaptive lasso.

The rest of the paper is organized as follows. The bootstrap methods for the Alasso is introduced and discussed in Section \ref{sec:mpb}. Assumptions and explanations of those are presented in Section \ref{sec:assum}. Results on rate of normal approximation is presented in Section \ref{sec:oracle}. Main results, i.e. results concerning the estimation properties of bootstrap are given in Section \ref{sec:mainresults}. %Section \ref{sec:simulation} presents simulation results exploring the finite-sample performance of the bootstrap in comparison with other methods for constructing confidence intervals based on Lasso estimators.
Proofs are presented in Section \ref{sec:proofs}.% Section \ref{sec:1.777} states concluding remarks.

\section{Description of the Bootstrap Methods}
\label{sec:mpb}
\subsection{Residual Bootstrap}
First, let us briefly describe the residual bootstrap method in penalized regression. The residuals
$\{\hat{\epsilon}_1,\dots,\hat{\epsilon}_n\}$ by $\hat{\epsilon}_i = y_i - x_i^\prime \hat{\beta}_n$. Suppose $\bar{\epsilon}_n$ is the mean of the residuals. Then select a random sample $\{\epsilon_1^*, \dots, \epsilon_n^*\}$ from $\{(\hat{\epsilon}_1-\bar{\epsilon}_n),\dots,(\hat{\epsilon}_n-\bar{\epsilon}_n)\}$ and define
\begin{align*}
y_i^* = x_i^\prime \hat{\beta}_n + \epsilon_i^*\;\;\;\;\;, i=1,\dots,n
\end{align*}
Then the residual bootstrap version of the penalized estimator is defined as
\begin{align}\label{eq:relasso}
\hat{\beta}_n^*=\operatorname*{arg\,min}_{t^*}\Big\{\sum_{i=1}^{n}(y_i^* - x^\prime_i t^*)^2
+n\sum_{j=1}^{p}P_{\lambda_n, j}(|t_j^*|)\Big\}.
\end{align}

Note that we have considered actual residuals, not the modified residuals unlike the construction of Chatterjee and Lahiri (2011).

\subsection{Perturbation Bootstrap}
Let $G_1^*,\ldots, G_n^*$ be $n$ independent copies of a non-degenerate random variable $G^* \in [0,\infty)$ having expectation $\mu_{G^*}$. These quantities will serve as perturbation quantities in the construction of the perturbation bootstrapped version of the penalized estimator. %\bl{\sout{The modified perturbation bootstrap in Alasso involves careful construction of the penalized objective function.} 
We define perturbation bootstrap version as the minimizer of a carefully constructed penalized objective function which involves the predicted values $\hat{y}_i = \bm{x}'_i\bm{\hat{\beta}}_n$, $i=1,\dots,n$ as well as the observed values $y_i,\dots,y_n$.
These sets of values appear in the objective function in two perturbed least-squares criteria. %\sout{The modified penalized objective function needs to incorporate the sum of two perturbed least square objective functions, one involving $y_i$, $i\in \{1,\cdots, n\}$ and other with $\hat{y}_i$, $i\in \{1,\cdots, n\}$, see  (\ref{eqn:mpb}).}}

We formally define the perturbation bootstrap version $\bm{\hat{\beta}}_n^{**}$ of the penalized estimator as\begin{align}\label{eqn:boot}
\hat{\bm{\beta}}_n^{**} = \operatorname*{arg\,min}_{t^*}&\Bigg\{\sum_{i=1}^{n}(y_i - x^\prime _i t^*)^2(G^*_i-\mu_{G^*}) \nonumber\\
&+\sum_{i=1}^{n}(\tilde{y}_i-x^\prime_i t^*)^2(2\mu_{G^*}-G_i^*)+\mu_{G^*}n\sum_{j=1}^{p}P_{\lambda_n, j}(|t_j|^*)\Bigg\}.
\end{align}

\par
We point out that the modified perturbation bootstrap estimator can be computed using existing algorithms. 
%Before stating the higher-order results, we mention an important property of our modification from the perspective of computation. Apparently, it seems that the minimization problems (2.1) and (2.2) are complicated. But some manipulations can be incorporated to reduce these complicated forms to much simpler form. The next result states this precisely. 
Define $\bm{L}_1(\bm{t}) = \sum_{i=1}^{n}(y_i - \bm{x}'_i \bm{t})^2(G^*_i-\mu_{G^*})+\sum_{i=1}^{n}(\hat{y}_i-\bm{x}'_i\bm{t})^2(2\mu_{G^*}-G_i^*)+\mu_{G^*}n\sum_{j=1}^{p}P_{\lambda_n, j}(|t_j|)$, $j =1,\cdots,p$. Now set $z_i=\hat{y}_i+\hat{\epsilon}_i\mu_{G^*}^{-1}(G_i^*-\mu_{G^*})$, where $\hat \epsilon_i = y_i - \hat{y}_i$ for $i=1,\dots,n$ and let $\bm{L}_2(\bm{t})=\sum_{i=1}^{n}\big(z_i-\bm{x}'_i\bm{t}\big)^2+n\sum_{j=1}^{p}P_{\lambda_n, j}(|t_j|)$. Then we have $\operatorname*{arg\,min}_{\bm{t}}\bm{L}_1(\bm{t})=\operatorname*{arg\,min}_{\bm{t}}\bm{L}_2(\bm{t})$. This allows us to compute $\hat{\bm{\beta}}^{**}_n$ by minimizing standard objective functions on some pseudo-values. Note that the perturbation bootstrapped estimator can be obtained simply by properly perturbing the residuals in the decomposition $y_i=\hat{y}_i + \hat{\epsilon}_i$, $i = 1,\dots,n$.

\section{Assumptions}\label{sec:assum}

We first introduce some notations required for stating our assumptions and useful for the proofs later. We denote the true parameter vector as $\bm{\beta}_n = (\beta_{1,n},\dots,\beta_{p,n} )'$, where the subscript $n$ emphasizes that the dimension $p:=p_n$ may grow with the sample size $n$.
Set $\mathcal{A}_n=\{j: \beta_{j,n}\neq 0\}$ and  $p_0:=p_{0,n}=|\mathcal{A}_n|$. For simplicity, we shall suppress the subscript $n$ in the notations $p_n$ and $p_{0n}$. Without loss of generality, we shall assume that $\mathcal{A}_n=\{1,\dots,p_0\}$. Let $\bm{C}_n=n^{-1}\sum_{i=1}^{n}\bm{x}_i\bm{x}'_i$ and partition it according to $\mathcal{A}_n = \{1,\dots,p_0\}$ as 
\begin{equation*}
\bm{C}_{n} = \begin{bmatrix}
\bm{C}_{11,n} \;\;\;\bm{C}_{12,n}\\
\bm{C}_{21,n}\;\;\; \bm{C}_{22,n}
\end{bmatrix},
\end{equation*}

where $\bm{C}_{11,n}$ is of dimension $p_0\times p_0$. Define $\tilde{\bm{x}}_i=\bm{C}_n^{-1}\bm{x}_i$ (when $p\leq n$) and $sgn(x) =-1, 0 ,1$ according as $x<0$, $x=0$, $x>0$, respectively. Suppose $\bm{D}_n$ is a known $q\times p$ matrix with $\text{tr}(\bm{D}_n\bm{D}'_n)=O(1)$ and $q$ is not dependent on $n$. Let $\bm{D}_n^{(1)}$ contains the first $p_0$ columns of $\bm{D}_n$.

%Define $\bm{W}_n=n^{-1/2}\sum_{i=1}^{n}\bm{x}_i\epsilon_i$ and $\bm{W}_n=\big(\bm{W}_n^{(1)'},\bm{W}_n^{(2)'}\big)'=(W_{1n},\ldots,W_{pn})'$, where $\bm{W}_n^{(1)}$ consists of the first $p_0$ components of $\bm{W}_n$. Also define a $p_0\times 1$ vector $\tilde{\bm{s}}_n^{(1)}$ which has $j$th component $\sgn(\beta_{j,n})|\tilde{\beta}_{j,n}|$.  Similarly define the bootstrap versions of $\bm{W}_n$ and $\tilde{\bm{s}}_n^{(1)}$ as $\breve{\bm{W}}_n^*$ and $\tilde{\bm{s}}_n^{*(1)}$  by replacing $\bm{\beta}_n$ and $\tilde{\bm{\beta}}_n$ and $\epsilon_i$ respectively by $\bm{\hat{\beta}}_n$, $\tilde{\bm{\beta}}_n^*$ and $\hat{\epsilon}_i$, where $\hat{\epsilon}_i = Y_i - \x_i'\hat{\bm{\beta}}_n$ are Alasso residuals.
Define
\begin{align*}
\bm{S}_{n} = \begin{bmatrix}
\bm{D}_n^{(1)}\bm{C}_{11,n}^{-1}\bm{D}_n^{(1)'}.\sigma^2 \;\;\;\bm{D}_n^{(1)}\bm{C}_{11,n}^{-1}\bar{\bm{x}}^{(1)}_n.\mu_{3}\\
\bar{\bm{x}}^{(1)'}_n\bm{C}_{11,n}^{-1}\bm{D}_n^{(1)'}.\mu_3\;\;\;\;\;\;\;\;\; (\mu_4-\sigma^4)
\end{bmatrix},
\end{align*}
where $\bar{\bm{x}}_n=n^{-1}\sum_{i=1}^{n}\bm{x}_i=(\bar{\bm{x}}^{(1)\prime}_n,\bar{\bm{x}}^{(2)\prime}_n)^\prime$, $\sigma^2=\mathbf{Var}(\epsilon_1)=\mathbf{E}(\epsilon_1^2)$, and where $\mu_3$ and $\mu_4$ are, respectively, the third and fourth central moments of $\epsilon_1$. Define in addition the $q\times p_0$ matrix $\check{\bm{D}}_n^{(1)}=\bm{D}_n^{(1)}\bm{C}_{11,n}^{-1/2}$ and the $p_0 \times 1$ vector $\check{\bm{x}}_i^{(1)}=\bm{C}_{11,n}^{-1/2}\bm{x}_i^{(1)}$. Let $K$ be a positive constant and $r$ be a positive number $\geq 1$ unless otherwise specified. $||\cdot||$ and $||\cdot||_{\infty}$ respectively denote the Euclidean norm and the Sup norm. By $\mathbf{P_*}$ and $\mathbf{E_*}$ we denote, respectively, probability and expectation with respect to the distribution of $G^{*}$ conditional upon the observed data. Write $\max_{j\in\{1,\dots,p_0\}}|\tilde{P}^\prime(|\beta_{j,n}|)| = \tilde{P}_1$, $\max_{j\in\{1,\dots,p_0\}}|\tilde{P}^{\prime\prime}(|\beta_{j,n}|)| = \tilde{P}_2$ and $\max_{j\in\{1,\dots,p_0\}}|\tilde{P}^{\prime\prime\prime}(|\beta_{j,n}|)| = \tilde{P}_3$ where the function $\tilde{P}(\cdot)$ is defined just before Theorem \ref{thm:oracleII} in section \ref{sec:oracle}.

\par
We now introduce our assumptions. 
\begin{enumerate}[label=(A.\arabic*)]
%\item $\mu_{G^*}=\mathbf{E_*}(G_1^*)< \infty$. On a set of probability more than $1-o(n^{-1/2})$ (specified otherwise if needed), for each $j\in \{(p_0+1),\ldots,p\}$, 
%\begin{align*}
%\dfrac{-\lambda_n}{2\sqrt{n}}|\tilde{\beta}_{j,n}|^{-\gamma}\leq \Big(\bm{C}_{21,n}\bm{C}_{11,n}^{-1}\Big)_{j.}\Big[\bm{W}_n^{(1)}-\dfrac{\lambda_n}{2\sqrt{n}}\tilde{s}_n^{(1)}\Big]-W_{j,n}\leq \dfrac{\lambda_n}{2\sqrt{n}}|\tilde{\beta}_{j,n}|^{-\gamma}
%\end{align*}
%And on a set of bootstrap probability more than $1-o_p(n^{-1/2})$ (specified otherwise if needed), for each $j\in \{(p_0+1),\ldots,p\}$, 
%\begin{align*}
%\dfrac{-\lambda_n}{2\sqrt{n}}|\tilde{\beta}_{j,n}^*|^{-\gamma}\leq \Big(\bm{C}_{21,n}\bm{C}_{11,n}^{-1}\Big)_{j.}\Big[\mu_{G^*}^{-1}\breve{\bm{W}}_n^{*(1)}-\dfrac{\lambda_n}{2\sqrt{n}}\tilde{s}_n^{*(1)}\Big]-\mu_{G^*}^{-1}\breve{W}_{j,n}^*\leq \dfrac{\lambda_n}{2\sqrt{n}}|\tilde{\beta}_{j,n}^*|^{-\gamma}
%\end{align*}
\item $\max\{|(\bm{C}_{21,n})_{j\cdot} \bm{C}_{11,n}^{-1} sgn(\bm{\beta}_n^{(1)})|:1\leq j \leq p_0\}<1-\eta$ for some $\eta>0$.
\item Let $\eta_{11,n}$ denote the smallest eigenvalue of the matrix $\bm{C}_{11,n}$.
\begin{enumerate}[label=(\roman*)]
\item $\eta_{11,n}>Kn^{-a}$ for some $a\in [0,1)$.
 \item[(i)$'$] $\eta_{11,n}>Kn^{-a}$ for some $a\in [0,1/4]$.
  \item $\max\{n^{-1}\sum_{i=1}^{n}|x_{i,j}|^{2r}:1\leq j \leq p\}+\{n^{-1}\sum_{i=1}^{n}\big|(\bm{C}_{11,n}^{-1})_{j.}\bm{x}_{i}^{(1)}\big|^{2r}:1\leq j \leq p_0\}+\max\{n^{-1}\sum_{i=1}^{n}|(\bm{C}_{21,n})_{j\cdot} \bm{C}_{11,n}^{-1} x_i^{(1)}|^{2r}:j\in \mathcal{A}_n\} = O(1)$.
  \item $\max\{c_{11,n}^{j,j}: 1\leq j\leq p_0\}=O(1)$, where $c_{11,n}^{j,j}$ is the $(j,j)$th element of $C_{11,n}^{-1}$.
  \item $\max\{n^{-1}\sum_{i=1}^{n}|\tilde{x}_{i,j}|^{2r}:1\leq j \leq p\}= O(1)$ (when $p\leq n$).
 % \item $n^{-1}\sum_{i=1}^{n}||\bm{D}_n^{(1)}\bm{C}_{11,n}^{-1}\bm{x}_i^{(1)}||^{2r}=O(1)$.
  %\item $\max\{n^{-1}\sum_{i=1}^{n}|\tilde{x}_{i,j}|^{2r}:1\leq j \leq p\} = O(1)$, where $\tilde{x}_{i,j}$ is the $j$th element of $\tilde{\bm{x}}_i$. (when $p\leq n$)
	%\item[(iii)$'$]$\max\{c_{11,n}^{j,j}: 1\leq j\leq p_0\}=O(1)$, where $c_{11,n}^{j,j}$ is the $(j,j)$th element of $C_{11,n}^{-1}$. (when $p > n$)	
\end{enumerate}
\item There exists a $\delta \in (0,1)$ such that for all $n>\delta^{-1}$,
\begin{enumerate}[label=(\roman*)]
  \item sup$\{\bm{x}'\check{\bm{D}}_n^{(1)}\check{\bm{D}}_n^{(1)'}\bm{x}:\bm{x}\in \mathcal{R}^{q}, ||\bm{x}||=1\}<\delta^{-1}$. %(Only for PB)
  \item $n^{-1}\sum_{i=1}^{n}||\check{\bm{D}}_n^{(1)}\check{\bm{x}}_i^{(1)}\check{\bm{x}}_i^{(1)'}\check{\bm{D}}_n^{(1)'}||^r = O(1)$. %sup$\{\bm{x}'\bm{D}_n^{(1)}\bm{C}_{11,n}^{-1}\bm{G}_n^{(1)}\bm{C}_{11,n}^{-1}\bm{D}_n^{(1)'}\bm{x}:\bm{x}\in \mathcal{R}^{q}, ||\bm{x}||=1\}<\delta^{-1}$.
  %\item inf$\{\bm{x}'\bm{D}_n^{(1)}\bm{C}_{11,n}^{-1}\bm{D}_n^{(1)'}\bm{x}:\bm{x}\in \mathcal{R}^{q}, ||\bm{x}||=1\}>\delta$.\\
  \item inf$\{\bm{x}'\bm{S}_n\bm{x}:\bm{x}\in \mathcal{R}^{q+1}, ||\bm{x}||=1\}>\delta$.
\end{enumerate}
\item $\max\{|\beta_{j,n}|:j\in \mathcal{A}_n\}=O(1)$.  %and min$\{|\beta_{j,n}|:j\in \mathcal{A}_n\}\geq Kn^{-b}$ for some $b\geq 0$ such that $4b< 1$ and $a+2b\leq 1$, where $a$ is defined as in (A.2)(i).
\item[(A.4)$'$]$\min\{|\beta_{j,n}|:j\in \mathcal{A}_n\}\geq K\lambda_n n^{-1}\min\{p_0^{3/2},n^a\sqrt{p_0}\}$ or $\min\{|\beta_{j,n}|:j\in \mathcal{A}_n\}\geq K\lambda_n n^{-1}\sqrt{p_0}||\bm{C}_{11,n}^{-1}||$ for some $K>0$.
\item \begin{enumerate}[label=(\roman*)]
  \item $\mathbf{E}|\epsilon_1|^{2r}< \infty$. $\mathbf{E}\epsilon_1=0$.
  %\item $\epsilon_1$ satisfies Cramer's condition: $\limsup_{|t|\rightarrow \infty}\mathbf{E}e^{it\epsilon_1}<1$.\\
  \item $(\epsilon_1, \epsilon_1^{2})$ satisfies Cramer's condition:\\ 
\hspace*{5mm}$\limsup_{||(t_1,t_2)||\rightarrow \infty}\mathbf{E}(exp(i(t_1\epsilon_1+t_2\epsilon_1^{2})))<1$.
\end{enumerate}
\item \begin{enumerate}[label=(\roman*)]
  \item $\mathbf{E_*}(G_{1}^{*})^{2r} < \infty$.  
$\mathbf{Var}(G_1^*) =\sigma^2_{G^*}= \mu_{G^*}^2$, $\mathbf{E_*}(G_1^* - \mu_{G^*})^3 = \mu_{G^*}^3$. (Only for PB)
  \item $G_{i}^{*}$ and $\epsilon_i$ are independent for all $1\leq i\leq n$. (Only for PB)
 % \item $(G^*_1 - \mu_{G^*})$ satisfies Cramer's condition:\\ 
%\hspace*{25mm}$\limsup_{|t|\rightarrow \infty}\mathbf{E}(exp(it(G_1^*-\mu_{G^*}))<1$

\item $((G^*_1-\mu_{G^*}), (G_1^* - \mu_{G^*})^{2})$ satisfies Cramer's condition:\\ 
\hspace*{9mm}$\limsup_{||(t_1,t_2)||\rightarrow \infty}\mathbf{E_*}(exp(i(t_1(G_1^*-\mu_{G^*})+t_2(G_1^* - \mu_{G^*})^{2})))<1$. (Only for PB)
\end{enumerate}
\item There exists $\delta_1 \in (0,1)$ such that for all $n>\delta_1^{-1}$,

\begin{enumerate} [label=(\roman*)]
%\item $\sqrt{n}\lambda_n\leq \delta_1^{-1}n^{-\delta_1} \text{min}\Big{\{}\dfrac{n^{-b\gamma}}{p_0},\dfrac{n^{-b\gamma-a/2}}{\sqrt{p_0}}\Big{\}}$.
\item $\sqrt{n}\lambda_n\leq \delta_1^{-1}n^{-\delta_1} \text{min}\Big{\{}\dfrac{\min\{\tilde{P}_1^{-1},\tilde{P}_2^{-1},\sqrt{n}\tilde{P}_3^{-1}\}}{p_0},\dfrac{n^{-a/2}\tilde{P}_1}{\sqrt{p_0}}\Big{\}}$.
\item[(i)$'$]
$p_0. \lambda_n. n^{-3/4}. \sqrt{\log n}.\max\{p_0, n^a\} = o(1)$ as $n\rightarrow \infty$.
\item $\lambda_n n^{(1+\gamma)/2} \geq \delta_1 n^{\delta_1} p_0$ where $\gamma$ is defined in section \ref{oracle}.
\item[(ii)$'$] $\dfrac{\lambda_n}{\sqrt{n}} \geq \delta_1 p_0 \sqrt{\log n}$ %\Big{(}$p_0$ can be dropped if we assume $\max\{n^{-1}\sum_{i=1}^{n}|(\bm{C}_{21,n})_{j\cdot} \bm{C}_{11,n}^{-1} x_i^{(1)}|^{2r}:j\in \mathcal{A}_n\}=O(1)$\Big{)}
\item $p_{0}=o\big(n^{1/2}(\log n)^{-3/2}\big)$.

\end{enumerate}
%\item[(A.7)] There exists $C\in(0,\infty)$ and $\delta_2\in(0,\gamma^{-1}\delta_1)$, $\delta_1$ being defined in the assumption (A.6), such that
%\begin{align*}
%&\mathbf{P}\Big(\max\{\big|\sqrt{n}(\tilde{\beta}_{j,n}-\beta_{j,n})\big|:1\leq j\leq p\}>C.n^{\delta_2}\Big)=o(n^{-1/2})\\
%& \mathbf{P_*}\Big(\max\{\big|\sqrt{n}(\tilde{\beta}_{j,n}^*-\hat{\beta}_{j,n})\big|:1\leq j\leq p\}>C.n^{\delta_2}\Big)=o_p(n^{-1/2})
%\end{align*}

\end{enumerate}
\par

Now we explain the assumptions briefly. Assumption is the weak irrepresentable condition required to achieve general sign consistency of Lasso [cf. Zhao and Yu (2006)]. Clearly this condition holds if $\max\{||(\bm{C}_{21,n})_{j\cdot} \bm{C}_{11,n}^{-1}|| :1\leq j \leq p_0\}<1$. Assumption (A.2) describes the regularity conditions needed on the growth of the design vectors. Assumption (A.1)(i) is a restriction on the smallest eigenvalue of $\bm{C}_{11,n}$. Assumption (A.2)(i) is a weaker condition than assuming that $\bm{C}_{11,n}$ converges to a positive definite matrix. (A.1)(ii) and (iii) are needed to bound the weighted sums of types $\big[\sum_{i=1}^{n}\bm{x}_i\epsilon_i\big]$, $\big[\bm{C}_{11,n}^{-1}\sum_{i=1}^{n}\bm{x}_i^{(1)}\epsilon_i\big]$. For general value of $r$, (A.1)(ii) and (iii) are much weaker than conditioning on $l_r$-norms of the design vectors. Here the value of $r$ is specified by the underlying Edgeworth expansion. 

\par
Assumptions (A.3)(i) bounds the eigenvalues of the matrix $\bm{D}_n^{(1)}\bm{C}_{11,n}^{-1}\bm{D}_n^{(1)'}$ away from infinity. It is necessary to obtain bounds needed in the studentized setup. Assumption (A.3)(ii) is a condition similar to the conditions in (A.3)(ii) and (iii); but involving the $q\times p$ matrix $\bm{D}_n$. This condition is needed for finding two term Edgeworth expansions in original Lasso estimator and also for showing necessary closeness of the covariance matrix estimators $\breve{\bm{\Sigma}}_n, \tilde{\bm{\Sigma}}_n$ [defined in Section \ref{sec:mainresults}] to their population counterparts (for details see Lemma \ref{lem:Sigma}) in case of perturbation bootstrap.  %which bounds the maximum eigenvalue of $\bm{D}_n^{(1)}\bm{C}_{11,n}^{-1}\bm{G}_n^{(1)}\bm{C}_{11,n}^{-1}\bm{D}_n^{(1)'}$, is somewhat restrictive; it may be dropped, but then the maximum growth rate of $p_0$ that can be allowed is $O(n^{1/4})$, so we choose to keep it. 
Assumption (A.2)(iii) bounds the minimum eigen value of the matrix $S_n$ away from $0$. This condition along with the Cramer conditions given in (A.4) and (A.5) enable certain Edgeworth expansions.
\par
Assumption (A.4) specifies the minimum magnitude of the non-zero regression coefficient required for the validity of bootstrap. More precisely, the condition (A.4) along with the condistion (A.1) are essential in obtaining a suitable form of the bootstrapped estimators from (2.1) and (2.2), required to achieve second order correctness.
%. The condition on the minimum is needed to ensure that the non-zero regression coefficients cannot converge to zero faster than the error rate, that is not faster than $O(n^{-1/2})$. We mention that one can assume $b<1/2$ instead of assuming $b<1/4$, but with the price of putting another restriction on the penalty parameter $\lambda_n$. We do not consider such a setting here. {We also want to point  out that it is not possible to relax this minimal signal condition by the bias correction, considered in Section \ref{sec:mainresults}.  With further relaxation, the bias of the Alasso estimator will be larger than the estimation error which is of order $O_p(n^{-1/2})$ and hence second-order correctness cannot be achieved by perturbation bootstrap in more relaxed minimal signal condition.}
\par
Assumption (A.5)(i) is a moment condition on the error term needed for valid Edgeworth expansion. Assumption (A.4)(ii) is Cramer's condition on the errors, which is very common in the literature of Edgeworth expansions; it is satisfied when the distribution of $(\epsilon_1, \epsilon_1^{2})$ has a non-degenerate component which is absolutely continuous with respect to the Lebesgue measure [cf. Hall (1992)]. Assumption (A.4)(ii) is only needed to get a valid Edgeworth expansion for the original Alasso estimator in the studentized setup. Assumptions (A.6)(i) and (iii) are the analogous conditions that are needed on the perturbing random quantities to get a valid Edgeworth expansion in the perturbation bootstrap setting. Assumption (A.5)(ii) is natural, since the $\epsilon_i$ are present already in the data generating process, whereas $G_i^*$ are introduced by the user. One can look for Generalized Beta and Generalized Gamma families for suitable choices of the distribution of $G^*$. The pdf of Generalized Beta family of distributions is
\begin{equation*}
  GB(y;f,g,h,\omega, \rho)=\left\{
  \begin{array}{@{}ll@{}}
    \dfrac{|f|y^{f\omega-1}\Big(1-(1-c)(y/g)^f\Big)^{\rho-1}}{g^{f\omega}B(\omega,\rho)\Big(1+c(y/g)^f\Big)^{\omega+\rho}}\;\;\;& \text{for}\;\; 0<y^f<\dfrac{g^f}{1-h}\\
     0\;\; &  \text{otherwise}
  \end{array}\right.
\end{equation*} 
where $0\leq h\leq 1$ and other parameters are all positive. We interpret $1/0$ as $\infty$. The function $B(\omega,\rho)$ is the beta function. Choices of the distribution of $G^*$ can be obtained by finding solution of $(f,g,h,\omega,\rho)$ from the following two equations
\begin{align*}
\dfrac{B(\omega+2/f,\rho)}{B(\omega,\rho)}& {}_2F_1\big[\omega+2/f,2/f;h;\omega+\rho+2/f\big]\\
&=2\bigg[\dfrac{B(\omega+1/f,\rho)}{B(\omega,\rho)}{}_2F_1\big[\omega+1/f,1/f;h;\omega+\rho+1/f\big]\bigg]^2\\
\text{and}\;\; \dfrac{B(\omega+3/f,\rho)}{B(\omega,\rho)}{}_2&F_1\big[\omega+3/f,3/f;h;\omega+\rho+3/f\big]\\
&=5\bigg[\dfrac{B(\omega+1/f,\rho)}{B(\omega,\rho)}{}_2F_1\big[\omega+1/f,1/f;h;\omega+\rho+1/f\big]\bigg]^3
\end{align*}
where ${}_2F_1$ denotes hypergeometric series. The pdf of Generalized Gamma family of distributions is given by
\begin{equation*}
  GG(y;\omega, \rho, \nu)=\left\{
  \begin{array}{@{}ll@{}}
    \dfrac{(\nu/\omega^\rho)y^{\rho-1}e^{(y/\omega)^{\nu}}}{\Gamma(\rho/\nu)}\;\;\;& \text{for}\;\; y>0\\
     0\;\; &  \text{otherwise}
  \end{array}\right.
\end{equation*} 
where all the parameters are positive and $\Gamma(\cdot)$ denotes the gamma function. For this family, the suitable choices of the distribution of $G^*$ can be obtained by considering any positive value of the parameter $\omega$ and solving the following two equations for $(\rho,\nu)$,
\begin{align*}
&\Big[\Gamma((\rho+2)/\nu)\Big]* \Gamma(\rho/\nu) = 2\Big[\Gamma((\rho+1)/\nu)\Big]^2\\
\text{and}\;\; &\Big[\Gamma((\rho+3)/\nu)\Big]* \Big[\Gamma(\rho/\nu)\Big]^2 = 5\Big[\Gamma((\rho+1)/\nu)\Big]^3.
\end{align*}
One immediate choice of the distribution of $G^*$ from Generalized Beta family is the Beta$(\alpha, \beta)$ distribution with $3\alpha=\beta=3/2$. We have utilized this distribution as the distribution of the perturbing quantities $G_i^*$'s in our simulations, presented in Section \ref{sec:simulation}. Outside these two generalized family of distributions, one possible choice is the distribution of $(M_1+M_2)$ where $M_1$ and $M_2$ are independent and $M_1$ is a Gamma random variable with shape and scale parameters $0.008652$ and $2$ respectively and $M_2$ is a Beta random variable with both the parameters $0.036490$. Another possible choice is the distribution of $(M_3+M_4)$ where $M_3$ and $M_4$ are independent and $M_3$ is an Exponential random variable with mean $\big{(}79-15\sqrt{33}\big{)}/16$ and $M_4$ is an Inverse Gamma random variable with both shape and scale parameters $\big{(}4 + \sqrt{11/3}\big{)}$.

Assumptions (A.7)(i)$'$ can be compared with $\lambda_n = o(n)$ which is required to achieve consistency of the lasso estimator in fixed dimensional setting [cf. Knight and Fu (2000)]. However, when dimension $p$ is increasing with $n$, one needs to have $\lambda_n.p_0.\max\{n^a,p_0\} = o(n)$ to achieve $\|\hat{\bm{\beta}_n}-\bm{\beta}\|_{\infty}=o_p(1)$. To accommodate studentiztion, one needs to assume (A.7)(i)$'$ and  (A.7)(iii). (A.7)(i)$'$ is little stronger than $\lambda_n.n^a.p_0 = o(n)$. On the other hand (A.7)(ii)$'$ along with (A.1) are required to obtain an exact form of the lasso estimator. Note that when $p_0$ is fixed and $a=0$, then one can choose $\lambda_n$ to be $n^l$ where $l\in (1/2, 3/4)$. When $a=0$ and $p_0$ can grow like $n^d$, then we can choose $\lambda_n$ to be of order $n^{7/12}$ %(\rd{$n^{1/2+\delta_5}$})
provided $d\in [0,1/12)$ %(\rd{$n^{1/8-\delta_5}$ for some $\delta_5\in (0,1/8)$}).
Again the assumptions (A.7)(i), (A.7)(ii) and (A.7)(iii) are required to establish second order correctness for the estimator belonging to class II. If we consider the adaptive lasso penalty, then assuming the true non-zero regression coefficients to be bounded away from zero we can consider $p_0 = O(n^{\gamma/5})$ and $\lambda_n = C.n^{-1/2-\gamma/4}$ for some constant $C> 0$, provided $\gamma \in (0,2)$.

\par

\section{Rate of Oracle Normal Approximation}\label{sec:oracle}

In this section, we are going to state the results on the error rate of the oracle normal approximation to the distribution of the penalized regression estimator.  Recall that $\Delta_n = \sup_{\bm{B}\in \mathcal{C}_q}|\mathbf{P}(\bm{T}_n\in \bm{B})-\Phi(\bm{B}; \sigma^2\bm{\Sigma}_n)|$ where $\bm{T}_n=\sqrt{n}\bm{D}_n(\hat{\bm{\beta}}_n - \bm{\beta})$ and $\bm{\Sigma}_n=\bm{D}_n^{(1)}\bm{C}_{11,n}^{-1}\bm{D}_n^{(1)\prime}$ with $\bm{D}_n$, $\bm{D}_n^{(1)}$, $\bm{C}_{11,n}$ being defined earlier. $\mathcal{C}_q$ is the collection of Borel measurable convex sets of $\mathcal{R}^q$ and $\Phi(\cdot; \bm{A})$ denotes the zero mean normal measure of the set $\bm{A}\in \mathcal{R}^q$. Without loss of generality assume that $\mathcal{A}_n=\{j:\beta_{j,n}\neq 0\}=\{1,\dots,p_0\}$. Define the sets $\bm{A}_n=\{\hat{\mathcal{A}}_n=\mathcal{A}_n\}$ and $\bm{B}_n=\{\hat{\bm{\beta}}_n^{(1)}=\bar{\bm{\beta}}_n^{(1)}\}$ where $\hat{\bm{\beta}}_n^{(1)}$ consists of first $p_0$ components of $\hat{\bm{\beta}}_n$ and $\bar{\bm{\beta}}_n^{(1)}$ is the OLS of $\bm{\beta}_n^{(1)}$ assuming $\bm{\beta}_n^{(2)}=\bm{0}$. Note that strong oracle property implies $\mathbf{P}(\bm{A}_n\cap \bm{B}_n)=o(1)$ as $n\rightarrow \infty$. Then we have the following result for the estimators in class I, that is when a penalized estimator satisfies strong oracle property.

\begin{theo}\label{thm:oracleI}
Suppose $\mathbf{P}(\bm{A}_n\cap \bm{B}_n)=1-o(n^{-1/2})$ and the conditions (A.3)(ii), (A.3)(iii) and (A.5)(i) hold with $r=3/2$. Then we have
\begin{align*}
\Delta_n=O(n^{-1/2}).
\end{align*}
\end{theo}

Theorem \ref{thm:oracleI} shows that when the strong oracle property holds, the inference on non-zero regression parameters based on oracle normal approximation has the same level of accuracy in increasing dimensions as in the simpler case of the OLS  when the dimension is fixed and no penalization is used. Note that the condition $\mathbf{P}(\bm{A}_n\cap \bm{B}_n)=1-o(n^{-1/2})$ is somewhat a high level condition and hence needs to be checked for the underlying penalized method which belong to class I. For example in case of MCP, $\mathbf{P}(\bm{A}_n\cap \bm{B}_n)=1-o(n^{-1/2})$ is true provided there exists $\delta \in (0,1)$ such that for all $n>\delta^{-1}$ the following conditions are satisfied:
\begin{enumerate}
\item For all $\bm{x}\in \mathcal{R}^{p_0}$, $\bm{y}\in \mathcal{R}^{p-p_0}$, 
\begin{align*}
\bm{x}^\prime \bm{C}_{12,n}\bm{y}\leq \delta^2 (\bm{x}^\prime \bm{C}_{11,n}\bm{x})(\bm{y}^\prime \bm{C}_{22,n}\bm{y})
\end{align*}
\item (A.2)(i)$'$ is true.
\item $\min\{|\beta_{j,n}|:j=1,\dots, p_0\} \geq \delta^{-1} \gamma \lambda_n$ and $\max\{|\beta_{j,n}|:j=1,\dots, p_0\}=O(1)$.
\item $\lambda_n\geq \dfrac{\delta^{-1}\log n}{\sqrt{n}}$ and $\lambda_n\leq \dfrac{\delta(\log n)^{-1}\sqrt{n}}{p_0}$.
\item $p_0\leq \delta (\log n)^{-3} n^{1-2a}$
\end{enumerate}
In case of one step estimator with $\tilde{P}_{\lambda_n}(\cdot)$ being the SCAD penalty, $\mathbf{P}(\bm{A}_n\cap \bm{B}_n)=1-o(n^{-1/2})$ is true provided there exists $\delta \in (0,1)$ such that for all $n>\delta^{-1}$ the following conditions are satisfied:
\begin{enumerate}
\item $\min\{|\beta_{j,n}|:j=1,\dots, p_0\} \geq \delta^{-1}  \lambda_n \log n$ and $\max\{|\beta_{j,n}|:j=1,\dots, p_0\}=O(1)$.
\item $\sqrt{n}\lambda_n>\delta^{-1}p_0\log n$.
\end{enumerate}

Next consider the estimators belonging to class II, that is when oracle property holds but not the strong oracle property. %\rd{For this class of estimators, hereafter we will only consider that the penalty for $j$th component $\bm{\beta}$ ,$P_{\lambda_n,j}(|\beta_{j,n}|)$, has the form $P_{\lambda_n,j}(|\beta_{j,n}|)= \tilde{P}^{\prime}_{\lambda_n}(|\tilde{\beta}_{j,n}|)|\beta_{j,n}|$ where $\tilde{P}^{\prime}_{\lambda_n}(\cdot)=\lambda_n\tilde{P}^{\prime}(\cdot)$ and $\tilde{\bm{\beta}}_n=(\beta_{1,n},\dots,\beta_{p_0,n})$ is some initial estimator, eg. OLS when $p\leq n$ and Lasso for $p>n$. $\tilde{P}(\cdot)$ is thrice differentiable on $(0,\infty)$ and $\max_j\tilde{P}^{\prime\prime\prime}(t_j)$ is continuous at $\bm{t}=(t_1,\dots,t_{p_0})^\prime=\Big(|\beta_{1,n}|,\dots,|\beta_{p_0,n}|\Big)^\prime$. For any sequence $\theta_n$ converging to $0+$, $\tilde{P}^{\prime}(\theta_n)=O\big(\theta_n^{-\gamma}\big)$ as $n\rightarrow \infty$, for some $\gamma>0$.} 
For this class of estimators, hereafter we will only consider that the penalty for $j$th component $\bm{\beta}$ ,$P_{\lambda_n,j}(|\beta_{j,n}|)$, has the form $P_{\lambda_n,j}(|\beta_{j,n}|)= \tilde{P}^{\prime}_{\lambda_n}(|\tilde{\beta}_{j,n}|)|\beta_{j,n}|$ where $\tilde{P}^{\prime}_{\lambda_n}(\cdot)=\lambda_n\tilde{P}^{\prime}(\cdot)$ and $\tilde{\bm{\beta}}_n=(\beta_{1,n},\dots,\beta_{p_0,n})$ is some initial estimator, eg. OLS when $p\leq n$ and Lasso for $p>n$. $\tilde{P}(\cdot)$ is thrice differentiable on $(0,\infty)$ and $\max_j\tilde{P}^{\prime\prime\prime}(t_j)$ is continuous at $\bm{t}=(t_1,\dots,t_{p_0})^\prime=\Big(|\beta_{1,n}|,\dots,|\beta_{p_0,n}|\Big)^\prime$. For any sequence $\theta_n$ converging to $0+$, $\tilde{P}^{\prime}(\theta_n)=O\big(\theta_n^{-\gamma}\big)$ as $n\rightarrow \infty$, for some $\gamma>0$. Before stating the result we also need to define the form of the bias as it has a significant contribution in the rate of convergence to the oracle normal limit. Define the bias term $\bm{b}_n=-\bm{D}_n^{(1)}\bm{C}_{11,n}^{-1}\bm{s}_n^{(1)}$ where $\bm{s}_n^{(1)}$ is a $p_0\times 1$ vector with $j$th component $s_{j,n}=\sqrt{n}P^\prime_{\lambda_n,j}(|\beta_{j,n}|)sgn(\beta_{j,n})$. Then we have the following result:

\begin{theo}\label{thm:oracleII}
Suppose the conditions (A.2)-(A.4), (A.5)(i), (A.7)(i) and (A.7)(ii) hold with $r=3/2$. Also assume $p\leq n$ and the initial estimator $\tilde{\bm{\beta}}_n=(\beta_{1,n},\dots,\beta_{p_0,n})$ is the OLS.  Then we have
\begin{align*}
\Delta_n=O\Big{(}n^{-1/2}+||\bm{b}_n||+%\rd{\lambda_n\cdot \min\{n^a,p_0\}\cdot \tilde{P}_2+\lambda_n^2p_0^2\tilde{P}_2^2}
\lambda_n\cdot \min\{n^a,p_0\}\cdot \tilde{P}_2+\lambda_n^2p_0^2\tilde{P}_2^2\Big{)}
\end{align*}
where $a$, $\tilde{P}_2$ are defined in section \ref{sec:assum}. %\bl{[red part is $o(n^{-1/2})$ if we want $\|\hat{\bm{b}}_n-\bm{b}_n\|=o_p(n^{-1/2})$ (required in case of bootstrap). One way out is to define the assumption (A.7)(i) separately for Theorem 2 and for the bootstrap theorems.]}
\end{theo}

Theorem \ref{thm:oracleII} gives the description of quantities which determine the rate of convergence. When only the oracle property holds, the rate of convergence to normal limit may be sub-$n^{-1/2}$ depending on the penalty terms and the constants $a,b$. Moreover under some additional conditions the above rate is attained with $d>0$. Therefore the inferences on non-zero regression coefficients based on normal approximation may have less accuracy than that when strong oracle property holds. Theorem \ref{thm:oracleII} includes Theorem 3.1 of Chatterjee and Lahiri (2013) which shows that the error rate for oracle approximation is sub-$n^{-1/2}$ in case of adaptive lasso.

Lastly consider the situation when only VSC holds, but not the oracle property. We will focus only on Lasso under the strong irrepresentable condition (A.1). Here the bias term is $\bm{b}_n^{\dagger}=-\bm{D}_n^{(1)}\bm{C}_{11,n}^{-1}\bm{c}_n^{(1)}$ where $\bm{c}_n^{(1)}$ is a $p_0\times 1$ vector with $j$th component $c_{j,n}=sgn(\beta_{j,n})$. Although the oracle property does not hold in this case, we want to point out how bad the oracl normal approximation is in case of Lasso.   

\begin{theo}\label{thm:oracleIII}
Suppose the conditions (A.1), (A.2)(ii), (A.3)(ii), (A.3) (iii), (A.5)(i), (A.7)(i)$'$ and (A.7)(ii)$'$ hold with $r=3/2$. Also consider the matrix $\bm{D}_n$ such that  $\max_j\Big|\Big(\bm{\Sigma}_n^{-1/2}\Big)_{j\cdot}\Big[\bm{D}_n^{(1)}\bm{C}_{11,n}^{-1}\bm{c}_n^{(1)}\Big]\Big|>\kappa$ for some $\kappa >0$, where $\Big(\bm{\Sigma}_n^{-1/2}\Big)_{j\cdot}$ is the $j$th row of $\bm{\Sigma}_n^{-1/2}$. Then we have
\begin{align*}
\Delta_n\rightarrow 1\; \text{as}\; n\rightarrow \infty.
\end{align*}
\end{theo}

Theorem \ref{thm:oracleIII} clearly shows that one can not even use oracle based normal approximation for the purpose of inference when only VSC holds, unlike the case for class I and class II. Therefore one needs to look into some alternative approximation techniques such as bootstrap to make valid inference. In the next section we will show that both residual and perturbation bootstrap approximations of the distribution of $\bm{T}_n=\sqrt{n}\bm{D}_n(\hat{\bm{\beta}}_n-\bm{\beta})$ have error rate $O_p(n^{-1/2})$ for all the three classes of penalized estimators. Moreover suitable studentization improves the rate to $o_p(n^{-1/2})$, that is both the bootstrap methods are second order correct.

\section{Bootstrap and its Higher Order Properties}
\label{sec:mainresults}
This section is divided into two sub-sections. The first one is on the bootstrap approximation of the distribution of $\bm{T}_n$. The second sub-section introduces studentizations of $\bm{T}_n$ and subsequently describes the higher order asymptotic properties of the residual and perturbation bootstrap methods.

\subsection{Bootstrap Approximation of $T_n$}
Recall that $\bm{T}_n=\sqrt{n}\bm{D}_n(\hat{\bm{\beta}}_n-\bm{\beta}_n)$ where $\bm{D}_n$ is a $q\times p$ matrix with $q$ being fixed. Define the residual bootstrap version of $\bm{T}_n$ as $\bm{T}_n^*=\sqrt{n}\bm{D}_n(\hat{\bm{\beta}}_n^*-\hat{\bm{\beta}}_n)$ and the perturbation bootstrap version of $\bm{T}_n$ as $\bm{T}_n^{**}=\sqrt{n}\bm{D}_n(\hat{\bm{\beta}}_n^{**}-\hat{\bm{\beta}}_n)$. Define $\Delta_n^*=\sup_{\bm{B}\in \mathcal{C}_q}\big{|}\mathbf{P}_*(\bm{T}_n^*\in \bm{B})-\mathbf{P}(\bm{T}_n\in B)\big{|}$ and $\Delta_n^{**}=\sup_{\bm{B}\in \mathcal{C}_q}\big{|}\mathbf{P}_*(\bm{T}_n^{**}\in \bm{B})-\mathbf{P}(\bm{T}_n\in B)\big{|}$. Recall that $\mathcal{A}_n=\{j:\beta_{j,n}\neq 0\}$. Suppose $\bm{\beta}_n^{(1)}$ consists of the components of $\bm{\beta}$ belonging to $\mathcal{A}_n$. Define the residual bootstrap analogues of $\bm{A}_n$ and $\bm{B}_n$ respectively as $\bm{A}_n^*=\{\mathcal{A}_n^*=\mathcal{\hat{A}}_n\}$ and $\bm{B}_n^*=\{\hat{\bm{\beta}}_n^{*(1)}=\bar{\bm{\beta}}_n^{*(1)}\}$, where $\mathcal{A}_n^*=\{j:\hat{\beta}_{j,n}^*\neq 0\}$, $\hat{\bm{\beta}}_n^{*(1)}$ consists of components of $\hat{\bm{\beta}}_n^*$ belonging to $\mathcal{A}_n^*$ and $\bar{\bm{\beta}}_n^{*(1)}$ is the residual bootstrapped OLS of $\bm{\beta}_n^{(1)}$ assuming $\bm{\beta}_n^{(2)}=\bm{0}$. Suppose $\bm{A}_n^{**}$ and $\bm{B}_n^{**}$ are corresponding perturbation bootstrap versions. In case of Class II, that is when only oracle property holds, the penalty term is generally defined based on an initial estimator $\tilde{\bm{\beta}}_n$. For example in case of adaptive lasso, one can use OLS as the initial estimator when $p\leq n$ and Lasso when $p>n$. Corresponding residual bootstrap version $\tilde{\bm{\beta}}_n^*$ of the initial estimator is defined similarly as in original case but after replacing $\{(y_i,\bm{x}_i^\prime):i=1,\dots,p\}$ with $\{(y_i^*,\bm{x}_i^\prime):i=1,\dots,p\}$ and $\bm{\beta}_n$ with $\hat{\bm{\beta}}_n$ where $\{y_1^*,\dots,y_n^*\}$ are as defined in section \ref{sec:mpb}. The perturbation bootstrap version of the initial estimator is defined as $\operatorname*{arg\,min}_{\bm{t}^*}\Big[\sum_{i=1}^{n}(y_i - \bm{x}'_i \bm{t}^*)^2(G^*_i-\mu_{G^*})
+\sum_{i=1}^{n}(\hat{y}_i-\bm{x}'_i\bm{t}^*)^2(2\mu_{G^*}-G_i^*)+\mu_{G^*}\tilde{\lambda}_n\sum_{j=1}^{p}|t_{j}^*|\Big]$ where $\{\hat{y}_1,\dots,\hat{y}_n\}$ are predicted values defined in section \ref{sec:mpb}. $\tilde{\lambda}_n=0$ is for OLS and $\tilde{\lambda}_n>0$ is for Lasso.  Following results show that both the bootstrap approximations attain the optimal error rate $O_p(n^{-1/2})$ for each of the classes of estimators.
\begin{theo}\label{thm:boot}
\begin{enumerate}
\item Consider class I, that is when strong oracle property holds. Suppose $\mathbf{P}(\bm{A}_n\cap \bm{B}_n)=1-o(n^{-1/2})$. Assume $\mathbf{P}_*(\bm{A}_n^*\cap \bm{B}_n^*)=1-o_p(n^{-1/2})$ for residual bootstrap and $\mathbf{P}_*(\bm{A}_n^{**}\cap \bm{B}_n^{**})=1-o_p(n^{-1/2})$ for perturbation bootstrap. Additionally let the conditions (A.3)(ii), (A.3)(iii), (A.5)(i) and (A.6)(i) [only for perturbation bootstrap] hold with $r=3$. Then we have
\begin{align*}
\Delta_n^*, \Delta_n^{**}=O_p(n^{-1/2})
\end{align*}
\item Consider class II, that is when only oracle property holds. Suppose the conditions (A.2)-(A.4), (A.5)(i), (A.7)(i) and (A.7)(ii) hold with $r=3$.   Also assume $p\leq n$ and the initial estimator $\tilde{\bm{\beta}}_n=(\beta_{1,n},\dots,\beta_{p_0,n})$ is the OLS. Then we have
\begin{align*}
\Delta_n^*=O_p(n^{-1/2})
\end{align*}
If in addition the condition (A.6)(i) holds with $r=3$, then 
\begin{align*}
\Delta_n^{**}=O_p(n^{-1/2})
\end{align*}
\item Now consider Lasso. Suppose the conditions (A.1), (A.2)(ii), (A.3)(ii), (A.3) (iii), (A.5)(i), (A.7)(i)$'$ and (A.7)(ii)$'$ hold with $r=3$. Then we have
\begin{align*}
\Delta_n^*=O_p(n^{-1/2})
\end{align*}
If in addition the condition (A.6)(i) holds with $r=3$, then 
\begin{align*}
\Delta_n^{**}=O_p(n^{-1/2})
\end{align*}
\end{enumerate}
\end{theo}

Theorem \ref{thm:boot} shows that before standardization or studentization, the bootstrap approximation has the optimal error rate of $O_p(n^{-1/2})$. This is in contrast with the error rate of the oracle normal approximation established in section \ref{oracle}. The rate of convergence to oracle limit worsens from class I to class III due to the increase in the contribution of the bias term in the estimators from class I to class III. This indicates that the bootstrap approximation is not affected by how large the bias is, even in increasing dimension.

Note that to achieve the error rate $O_p(n^{-1/2})$ for the bootstrap approximation of the distribution of $\bm{T}_n$, we need to assume $\mathbf{P}_*(\bm{A}_n^*\cap \bm{B}_n^*)=1-o_p(n^{-1/2})$ for residual bootstrap and $\mathbf{P}_*(\bm{A}_n^{**}\cap \bm{B}_n^{**})=1-o_p(n^{-1/2})$ for perturbation bootstrap along with  $\mathbf{P}(\bm{A}_n\cap \bm{B}_n)=1-o(n^{-1/2})$. In general the conditions that imply $\mathbf{P}(\bm{A}_n\cap \bm{B}_n)=1-o(n^{-1/2})$ will also imply $\mathbf{P}_*(\bm{A}_n^*\cap \bm{B}_n^*)=1-o_p(n^{-1/2})$ and $\mathbf{P}_*(\bm{A}_n^{**}\cap \bm{B}_n^{**})=1-o_p(n^{-1/2})$ provided some higher order moments of $\epsilon_1$ and $G_1^*$ exist.

\subsection{Studentization of $T_n$ and higher order results}
Recall that $\mathcal{A}_n=\{j:\beta_{j,n}\neq 0\}$ denotes the set of significant regression coefficients. $\hat{\mathcal{A}}_n=\{j:\hat{\beta}_{j,n}\neq 0\}$ is an estimator of $\mathcal{A}_n$. Without loss of generality assume $\mathcal{A}_n=\{1,\dots,p_0\}$.  Recall that $\bm{A}_n=\{\hat{\mathcal{A}}_n=\mathcal{A}_n\}$ and $\bm{B}_n=\{\hat{\bm{\beta}}_n^{(1)}=\bar{\bm{\beta}}_n^{(1)}\}$ where $\hat{\bm{\beta}}_n^{(1)}$ consists of first $p_0$ components of $\hat{\bm{\beta}}_n$ and $\bar{\bm{\beta}}_n^{(1)}$ is the OLS of $\bm{\beta}_n^{(1)}$ assuming $\bm{\beta}_n^{(2)}=\bm{0}$. Note that on the set $\bm{A}_n\cap \bm{B}_n$, $\bm{T}_n=\sqrt{n}\bm{D}_n(\bm{\hat{\beta}}_n-\bm{\beta}_n)=n^{-1/2}\sum_{i=1}^{n}\bm{D}_n^{(1)}\bm{C}_{11,n}^{-1}\bm{x}_i^{(1)}\epsilon_i$ and hence the asymptotic variance of $\bm{T}_n/\sigma$ is $\bm{\Sigma}_n=n^{-1}\sum_{i=1}^{n}\bm{\xi}_i^{(0)}\bm{\xi}_i^{(0)\prime}$ where $\bm{\xi_i}^{(0)}=\bm{D}_n^{(1)}\bm{C}_{11,n}^{-1}\bm{x}_i^{(1)}$.

Recall the set $\hat{\mathcal{A}}_n=\{j: \hat \beta_{j,n}\neq 0\}$ and define  $\hat p_{0,n}=|\hat{\mathcal{A}}_n|$, supposing, without loss of generality, that $\hat{\mathcal{A}}_n=\{1,\ldots,\hat{p}_{0,n}\}$. We then partition the matrix $\bm{C}_n=n^{-1}\sum_{i=1}^{n}\bm{x}_i\bm{x}'_i$ as 
\begin{equation*}
\bm{C}_{n} = \begin{bmatrix}
\hat{\bm{C}}_{11,n} \;\;\;  \hat{\bm{C}}_{12,n}\\
\hat{\bm{C}}_{21,n}\;\;\; \hat{\bm{C}}_{22,n}
\end{bmatrix},
\end{equation*}
where $\hat{\bm{C}}_{11,n}$ is of dimension $\hat{p}_{0,n}\times \hat{p}_{0,n}$. Similarly, we define $\hat{\bm{D}}_n^{(1)}$ as the matrix containing the first $\hat{p}_{0,n}$ columns of $\bm{D}_n$ and we define $\hat{\bm{x}}_i^{(1)}$ as the vector containing the first $\hat{p}_{0,n}$ entries of $\bm{x}_i$. Then define the studentized versions of $\bm{T}_n$ as 
\begin{align*}
\bm{R}_n = \bm{T}_n/\hat{\sigma}_n\;\;\; \text{and}\;\;\; \check{\bm{R}}_n=\check{\sigma}_n^{-1}\hat{\Sigma}_n^{-1/2}\bm{T}_n
\end{align*}
where $\hat{\sigma}_n^2=n^{-1}\sum_{i=1}^{n}(\hat{\epsilon}_i-\bar{\epsilon}_n)^2$, $\check{\sigma}_n^2=n^{-1}\sum_{i=1}^{n}\hat{\epsilon}_i^2$ and $\hat{\bm{\Sigma}}_n=n^{-1}\sum_{i=1}^{n}\hat{\bm{\xi}}_i^{(0)}\hat{\bm{\xi}}_i^{(0)\prime}$ with $\hat{\bm{\xi}}_i^{(0)}=\hat{\bm{D}}_n^{(1)}\hat{\bm{C}}_{11,n}^{-1}\hat{\bm{x}}_i^{(1)}$ and $\epsilon_i=y_i-\bm{x}_i^{\prime}\hat{\bm{\beta}}_n$.
We are going to use $\bm{R}_n$ and $\check{\bm{R}}_n$ respectively in case of residual and perturbation bootstrap methods. Define the residual bootstrap version of $\bm{R}_n$ as $\bm{R}_n^*=\bm{T}_n^*/\sigma_n^*$ and the perturbation bootstrap version of $\check{\bm{R}}_n$ as $\check{\bm{R}}_n^*=\sigma_n^{**-1}\check{\sigma}_n\tilde{\Sigma}_n^{-1/2}\bm{T}_n^{**}$ where $\sigma_n^{*2}=n^{-1}\sum_{i=1}^{n}\epsilon_i^{*2}$, $\sigma_n^{**2}=\mu_{G^*}^{-2}n^{-1}\sum_{i=1}^{n}\epsilon_i^{**2}(G_i^*-\mu_{G^*})^2$, $\tilde{\bm{\Sigma}}_n=n^{-1}\sum_{i=1}^{n}\hat{\bm{\xi}}_i^{(0)}\hat{\bm{\xi}}_i^{(0)\prime}\hat{\epsilon}_i^2$  with $\{\epsilon_1^*, \dots, \epsilon_n^*\}$ are as defined in the section \ref{sec:mpb} and $\epsilon_i^{**}=y_i-\bm{x}_i^{\prime}\hat{\bm{\beta}}_n^{**}$. Then we have the following result for the estimators which belong to class I:

\begin{theo}\label{thm:bootI}
Suppose $\mathbf{P}(\bm{A}_n\cap \bm{B}_n)=1-o(n^{-1/2})$, $\mathbf{P}_*(\bm{A}_n^*\cap \bm{B}_n^*)=1-o_p(n^{-1/2})$ for residual bootstrap and $\mathbf{P}_*(\bm{A}_n^{**}\cap \bm{B}_n^{**})=1-o_p(n^{-1/2})$ for perturbation bootstrap. Additionally the conditions (A.3), (A.5) (A.7)(iii) hold with $r=4$. Then we have
\begin{align*}
\sup_{\bm{B}\in \mathcal{C}_q}\big{|}\mathbf{P}_*(\bm{R}_n^*\in \bm{B})-\mathbf{P}(\bm{R}_n\in \bm{B})\big{|}=o_p(n^{-1/2}).
\end{align*}
If in addition the condition (A.6) holds with $r=4$, then the following is also true:
\begin{align*}
\sup_{\bm{B}\in \mathcal{C}_q}\big{|}\mathbf{P}_*(\check{\bm{R}}_n^{*}\in \bm{B})-\mathbf{P}(\check{\bm{R}}_n\in \bm{B})\big{|}=o_p(n^{-1/2}).
\end{align*}
\end{theo}

Now consider class II, that is when only the oracle property holds. Here we need to assume some moderate deviation bounds on the initial estimator and its bootstrap versions. 
%For example in case of adaptive lasso, one can use OLS as the initial estimator ($\tilde{\bm{\beta}}_n$) when $p\leq n$ and Lasso when $p>n$. Corresponding residual bootstrap version $\tilde{\bm{\beta}}_n^*$ of the initial estimator is defined similarly as in original case but after replacing $\{(y_i,\bm{x}_i^\prime):i=1,\dots,p\}$ with $\{(y_i^*,\bm{x}_i^\prime):i=1,\dots,p\}$ and $\bm{\beta}_n$ with $\hat{\bm{\beta}}_n$ where $\{y_1^*,\dots,y_n^*\}$ are as defined in section \ref{sec:mpb}. The perturbation bootstrap version of the initial estimator is defined as $\operatorname*{arg\,min}_{\bm{t}^*}\Big[\sum_{i=1}^{n}(y_i - \bm{x}'_i \bm{t}^*)^2(G^*_i-\mu_{G^*})+\sum_{i=1}^{n}(\hat{y}_i-\bm{x}'_i\bm{t}^*)^2(2\mu_{G^*}-G_i^*)+\mu_{G^*}\tilde{\lambda}_n\sum_{j=1}^{p}|t_{j}^*|\Big]$ where $\{\hat{y}_1,\dots,\hat{y}_n\}$ are predicted values defined in section \ref{sec:mpb}. $\tilde{\lambda}_n=0$ is for OLS and $\tilde{\lambda}_n>0$ is for Lasso.
To obtain second order correctness in case of estimators in class II one needs to have 
\begin{align}\label{eqn:initialmoderate}
&\mathbf{P}\Big(\max\{\big|\sqrt{n}(\tilde{\beta}_{j,n}-\beta_{j,n})\big|:1\leq j\leq p\}>C.n^{\delta_2}\Big)=o(n^{-1/2})\\\nonumber
& \mathbf{P_*}\Big(\max\{\big|\sqrt{n}(\tilde{\beta}_{j,n}^*-\hat{\beta}_{j,n})\big|:1\leq j\leq p\}>C.n^{\delta_2}\Big)=o_p(n^{-1/2})
\end{align}
for some $C\in(0,\infty)$.
$\delta_2\in(0,\gamma^{-1}\delta_1)$ if $\gamma \geq 1$ and $\delta_2\in(0,\delta_1)$ if $\gamma < 1$, where $\delta_1$ being defined in the assumption (A.7) and $\gamma$ is defined in the form of the penalty term considered for class II. See Chatterjee and Lahiri (2013) and Das et al. (2018) for details regarding the requirements on the initial estimator and its bootstrap version for establishing higher order results in adaptive lasso. The following theorem shows that both residual and perturbation bootstrap are second order correct for class II:

\begin{theo}\label{thm:bootII}
Suppose the conditions (A.2)-(A.5) and (A.7) hold with $r=4$. Additionally assume that the moderate deviation bounds like (\ref{eqn:initialmoderate}) holds for the initial estimator and its bootstrap version. Then we have
\begin{align*}
\sup_{\bm{B}\in \mathcal{C}_q}\big{|}\mathbf{P}_*(\bm{R}_n^*\in \bm{B})-\mathbf{P}(\bm{R}_n\in \bm{B})\big{|}=o_p(n^{-1/2}).
\end{align*}
If in addition the condition (A.6) holds with $r=4$, then the following is also true:
\begin{align*}
\sup_{\bm{B}\in \mathcal{C}_q}\big{|}\mathbf{P}_*(\check{\bm{R}}_n^{*}\in \bm{B})-\mathbf{P}(\check{\bm{R}}_n\in \bm{B})\big{|}=o_p(n^{-1/2}).
\end{align*}
\end{theo}

Next consider Lasso, that is when only VSC holds, not the oracle property. We have seen in the section \ref{sec:oracle} that the oracle normal approximation is of no use in this case. However one can use either residual or perturbation bootstrap approximation of $\bm{T}_n=\sqrt{n}\bm{D}_n(\hat{\bm{\beta}}_n-\bm{\beta}_n)$ to make inferences. Now the question is whether it is possible to establish second order correctness for Lasso. Apparently it seems that the error rate can not be improved to $o_p(n^{-1/2})$ from $O_p(n^{-1/2})$ due to substantially large bias. Therefore only option here is to correct $\bm{T}_n$ for bias and then consider suitable studentization. Note that for Lasso the form of $\bm{T}_n$ is 
\begin{align*}
\bm{T}_n=n^{-1/2}\sum_{i=1}^{n}\bm{D}_n^{(1)}\bm{C}_{11,n}^{-1}\bm{x}_i^{(1)}\epsilon_i -\dfrac{\lambda_n}{2\sqrt{n}}\bm{D}_n^{(1)}\bm{C}_{11,n}^{-1}\bm{s}_n^{\dagger(1)}
\end{align*}
where $\bm{s}_n^{\dagger(1)}$ is a $p_0\times 1$ vector with $j$th component $sgn(\beta_{j,n})$. The bias is $\bm{b}_n^\dagger=-\dfrac{\lambda_n}{2\sqrt{n}}\bm{D}_n^{(1)}\bm{C}_{11,n}^{-1}\bm{s}_n^{\dagger(1)}$. An estimator of $\bm{b}_n^\dagger$ is $\hat{\bm{b}}^\dagger_n=-\dfrac{\lambda_n}{2\sqrt{n}}\hat{\bm{D}}_n^{(1)}\hat{\bm{C}}_{11,n}^{-1}\hat{\bm{s}}_n^{\dagger(1)}$ where $\hat{\bm{D}}_n^{(1)}$ and $\hat{\bm{C}}_{11,n}$ are defined earlier in this section and $\hat{\bm{s}}_n^{\dagger(1)}$ is a $p_0\times 1$ vector with $j$th component $sgn(\hat{\beta}_{j,n})$. Therefore the bias corrected studentized versions of $\bm{T}_n$ are

\begin{align*}
\breve{\bm{R}}_n = \hat{\sigma}_n^{-1}(\bm{T}_n-\hat{\bm{b}}_n^\dagger)\;\;\; \text{and}\;\;\; \tilde{\bm{R}}_n=\check{\sigma}_n^{-1}\hat{\Sigma}_n^{-1/2}(\bm{T}_n-\hat{\bm{b}}^\dagger_n)
\end{align*}

As opposed to $\bm{R}_n$ and $\check{\bm{R}}_n$, we are going to use $\breve{\bm{R}}_n$ and $\tilde{\bm{R}}_n$ respectively for residual and perturbation bootstrap. The residual bootstrap version of $\breve{\bm{R}}_n$ and the perturbation bootstrap version of $\tilde{\bm{R}}_n$ are respectively
\begin{align*}
\breve{\bm{R}}_n^* = \sigma_n^{*-1}(\bm{T}_n^*-\hat{\bm{b}}_n^\dagger)\;\;\; \text{and}\;\;\; \tilde{\bm{R}}_n^*=\sigma_n^{**-1}\check{\sigma}_n\tilde{\Sigma}_n^{-1/2}(\bm{T}_n^{**}-\hat{\bm{b}}_n^\dagger).
\end{align*}
Then we have the following higher order results for lasso

\begin{theo}\label{thm:bootIII}
Suppose the conditions (A.1), (A.2)$'$,(A.3), (A.4)$'$, (A.5) and (A.7)$'$ hold with $r=4$. Then we have
\begin{align*}
\sup_{\bm{B}\in \mathcal{C}_q}\big{|}\mathbf{P}_*(\breve{\bm{R}}_n^*\in \bm{B})-\mathbf{P}(\breve{\bm{R}}_n\in \bm{B})\big{|}=o_p(n^{-1/2}).
\end{align*}
If in addition the condition (A.6) holds with $r=4$, then the following is also true:
\begin{align*}
\sup_{\bm{B}\in \mathcal{C}_q}\big{|}\mathbf{P}_*(\tilde{\bm{R}}_n^{*}\in \bm{B})-\mathbf{P}(\tilde{\bm{R}}_n\in \bm{B})\big{|}=o_p(n^{-1/2}).
\end{align*}
\end{theo}

\begin{remark}
Theorem \ref{thm:bootI}-\ref{thm:bootIII} show that both the residual and perturbation bootstrap approximations are second order correct irrespective of the underlying class of the penalized estimators. These are remarkable results in view of the theorems \ref{thm:oracleI}-\ref{thm:oracleIII}. This indicate that bootstrap methods are somehow immune towards the effect of bias, unlike the case of oracle normal approximation. Most interesting case is the Lasso where one can not even use the oracle approximation as the error converges to $1$, see theorem \ref{thm:oracleIII}.
\end{remark}

\begin{remark}
Note that in case of estimators of class I and II, no bias correction is necessary for the bootstrap to achieve second order correctness. However for lasso, bias correction is indispensable. The reason is that the order of the bias term $\bm{b}_n^{\dagger}$ is higher than the order of $\bm{T}_n-\bm{b}_n^{\dagger}$. Studentization without any bias correction is of no help in improving the rate of convergence of the bootstrap from $O_p(n^{-1/2})$ to $o_p(n^{-1/2})$, unlike the situation in case of the estimators belonging to class II. This is the same reason why normal approximation fails drastically.  
\end{remark}

\begin{remark}\label{rem3}
The rate of bootstrap approximation in Theorem \ref{thm:bootI} can be improved to $O_p(n^{-1})$ with the same studentized pivots, provided $\mathbf{P}(\bm{A}_n\cap \bm{B}_n)=1-o(n^{-1})$ and $\mathbf{P}_*(\bm{A}_n^*\cap \bm{B}_n^*)=1-o_p(n^{-1})$ or $\mathbf{P}_*(\bm{A}_n^{**}\cap \bm{B}_n^{**})=1-o_p(n^{-1})$ are satisfied, that is, if the strong oracle property in the original and bootstrap regime are true with higher probability. This are generally true if suitable value of $r$ can be assumed in Theorem \ref{thm:bootI}. In case of lasso, the rate can be improved to $O_p(n^{-1})$ with the same bias corrected studentized pivots if all the conditions, mentioned in the theorem, are true with $r=6$. However, for the estimators in class II, the rate can not in general be improved to $O_p(n^{-1})$ using the same studentized pivots. The reason being the non-trivial contribution of the bias term. Here one needs to correct for the bias both in the original and bootstrap regime, similar to the case in Lasso, to achieve $O_p(n^{-1})$. In case of residual bootstrap, the bias corrected original and bootstrap pivots can be defined as
\begin{align*}
\dot{\bm{R}}_n= \hat{\sigma}^{-1}_n\big[\bm{T}_n+ \breve{b}_n\big]\;\;\; \text{and}\;\;\; \dot{\bm{R}}_n^*=\sigma^{*-1}_n\big[\bm{T}_n^*+ \breve{b}_n^*\big].
\end{align*}
Similarly the perturbation bootstrap version of the bias corrected pivots can be defined as
\begin{align*}
\ddot{\bm{R}}_n= \check{\sigma}_n^{-1}\hat{\Sigma}_n^{-1/2}\big[\bm{T}_n+ \breve{b}_n\big]\;\;\; \text{and}\;\;\; \ddot{\bm{R}}_n^*=\sigma_n^{**-1}\check{\sigma}_n\tilde{\Sigma}_n^{-1/2}\big[\bm{T}_n^{**}+ \breve{b}_n^{**}\big].
\end{align*}
Here, $\breve{\bm{b}}_n=\hat{\bm{D}}_n^{(1)}\hat{\bm{C}}_{11,n}^{-1}\breve{\bm{s}}_n^{(1)}$ with $\breve{\bm{s}}_n^{(1)}$ being a $\hat{p}_0\times 1$ vector with $j$th component $\breve{s}_{j,n}=\sqrt{n}P^\prime_{\lambda_n,j}(|\tilde{\beta}_{j,n}|)sgn(\hat{\beta}_{j,n})$, $\tilde{\beta}_{j,n}$ being the $j$th component of the initial estimator $\tilde{\bm{\beta}}_n$ and $\hat{\bm{D}}_n^{(1)}$ \& $\hat{\bm{C}}_{11,n}$
 are defined in subsection 5.2. $\breve{b}_n^*=\bm{D}_n^{*(1)}\bm{C}_{11,n}^{*-1}\breve{\bm{s}}_n^{*(1)}$ where $\bm{D}_n^{*(1)}$ \& $\bm{C}_{11,n}^{*}$ are same as $\hat{\bm{D}}_n^{(1)}$ \& $\hat{\bm{C}}_{11,n}$ but after replacing $\bm{\beta}_n$, $\hat{\bm{\beta}}_n$, $\tilde{\bm{\beta}}_n$ respectively by $\hat{\bm{\beta}}_n$, $\hat{\bm{\beta}}_n^*$ and $\tilde{\bm{\beta}}_n^*$. $\breve{b}_n^{**}=\bm{D}_n^{**(1)}\bm{C}_{11,n}^{**-1}\breve{\bm{s}}_n^{**(1)}$ is the perturbation bootstrap version of $\breve{b}_n$, defined in similar fashion as $\breve{b}_n^{*}$. All other notations are defined in subsection 5.2. Results for adaptive lasso after bias correction can be found in Chatterjee and Lahiri (2013) and Das et al. (2019).
\end{remark}

\subsection{Symmetric Bootstrap Confidence Intervals}
The second order results of the previous subsection directly imply that the one-sided bootstrap confidence intervals perform much better in terms of coverage error than one-sided intervals based on normal approximation. However the second order results can not separate the two-sided bootstrap intervals from their normal counterparts in terms of coverage error and hence two sided bootstrap confidence intervals need to be studied separately. For detailed discussion on this issue, see Chapter 2 of Hall (1992). Among all the two-sided bootstrap intervals most promising one is the symmetric bootstrap confidence interval, due to its coverage error of order $O_p(n^{-2})$ in most classical setups. See for example Hall (1988) for theoretical aspects of symmetric confidence intervals based on Efron's bootstrap.  

Suppose $H_n$ is a pivotal quantity for the parameter $\theta_n$. $H_n^*$ is the bootstrap version of the original pivot $H_n$. Define the original quantile $h_{n,\alpha}=\inf\{x:\mathbf{P}\big(|H_n|\leq x\big)\geq 1-\alpha\}$ and the bootstrap quantile $\hat{h}_{n,\alpha}=\inf\{x:\mathbf{P}_*\big(|H_n^*|\leq x\big)\geq 1-\alpha\}$, for some $\alpha \in (0,1)$. Now note that without loss of generality we can simply define $h_{n,\alpha}$ and $\hat{h}_{n,\alpha}$ respectively as solutions of $\mathbf{P}\big(|H_n|\leq x\big) = 1-\alpha$ and $\mathbf{P}_*\big(|H_n^*|\leq x\big) = 1-\alpha$, due to the Cramer's conditions [cf. assumption (A.5)(ii) for $\epsilon_i$'s, (A.6)(iii) on $G_1^*$'s for perturbation bootstrap and the restricted Cramer's condition Lemma \ref{lem:residualCramer} for $\epsilon_i^*$'s in case of the residual bootstrap]. Cramer's condition implies that the heaviest atom of the distribution of the underlying pivot has mass $O(e^{-\delta_6n})$ for some $\delta_6>0$ [cf. Theorem 2.3 in Hall (1992)], which is negligible compared to the desired coverage error $O(n^{-2})$. This is the same reason why there is no difference between openness and closedness of any interval considered in this paper. By symmetric bootstrap confidence interval of $\theta_n$ based on $H_n$ and $H_n^*$, here we mean the interval $I_{n,(1-\alpha)}$ with the property that the event $\{\theta \in I_{n,(1-\alpha)}\}$ is same as the event $\{|H_n|\leq \hat{h}_{n,\alpha}\}$. Clearly, the event $\{|H_n|\leq \hat{h}_{n,\alpha}\}$ is an estimator of the ideal event $\{|H_n|\leq h_{n,\alpha}\}$ which corresponds to the exact symmetric confidence interval. In most of the situations $H_n=a_n(\hat{\theta}_n-\theta_n)$ for some estimator $\hat{\theta}_n$ and scaling $a_n$, resulting $I_{n.(1-\alpha)}$ to be symmetric around $\hat{\theta}_n$. 

In this subsection we are interested on the coverage accuracy of symmetric bootstrap confidence intervals of a linear combination of the components of $\bm{\beta}_n$. Hence $\theta_n=\bm{D}_n\bm{\beta}_n$ with $\bm{D}_n$ being a $p-$dimensional row vector throughout this subsection. While exploring symmetric confidence interval based on residual bootstrap, we are going to consider the pair $\big(H_n,H_n^*\big)$ to be $\big(\bm{R}_n, \bm{R}_n^*\big)$, $\big(\dot{\bm{R}}_n,\dot{\bm{R}}_n^*\big)$ or $\big(\breve{\bm{R}}_n, \breve{\bm{R}}_n^*)$ according as the estimator $\hat{\bm{\beta}}_n$ falls in class I, II or III. Similarly for perturbation bootstrap, $\big(H_n,H_n^*\big)$ is $\big(\check{\bm{R}}_n, \check{\bm{R}}_n^*\big)$, $\big(\ddot{\bm{R}}_n,\ddot{\bm{R}}_n^*\big)$ or $\big(\tilde{\bm{R}}_n, \tilde{\bm{R}}_n^*)$ according as the estimator $\hat{\bm{\beta}}_n$ falls in class I, II or III. %[\rd{bias correction will also do the trick for class II, because only the sign will matter, see page 18 of the supplementary of Chatterjee and Lahiri (2013)}]
We denote the residual bootstrap pair by $\big(H_n^r,H_n^{r*}\big)$ and the perturbation bootstrap pair by $\big(H_n^p,H_n^{p*}\big)$. We are going to see that the symmetric bootstrap confidence interval based on the pair $\big(H_n^r,H_n^{r*}\big)$ result in a two-sided interval with coverage error $O(n^{-2})$, where as the coverage error remains $O(n^{-1})$ for the symmetric confidence interval based on $\big(H_n^p,H_n^{p*}\big)$. We introduce some correction term in the form of the symmetric confidence interval based on $\big(H_n^p,H_n^{p*}\big)$ which will result in an error of $O(n^{-2})$. Let us discuss on how to get hold of the correction term in case of perturbation bootstrap and also on why we are achieving error of order $O(n^{-2})$, before formally stating the result. Note that the asymptotic variances of $H_n^p$ and $H_n^{p*}$ are both equal to 1, where as that of $H_n^r$ and $H_n^{r*}$ may not be 1. However for notational simplicity, we assume that the asymptotic variances of the pivots $H_n^r$ and $H_n^{r*}$ are also equal to 1. All the following arguments in this subsection will go through for arbitrary variances.

Define, $\mathbf{P}\big(|H_n^r|\leq h_{n,\alpha}^r\big)=1-\alpha$ and $\mathbf{P}_*\big(|H_n^{r*}|\leq \hat{h}_{n,\alpha}^r\big)=1-\alpha$. Therefore $(1-\alpha)\%$ symmetric bootstrap confidence interval of $\theta_n$ based on $\big(H_n^r,H_n^{r*}\big)$ is $I^r_{n,1-\alpha}$ where the event $\{\theta_n \in I_{n,(1-\alpha)}^r\}$ is same as the event $\{|H_n^r|\leq \hat{h}_{n,\alpha}^r\}$. Hence it is enough to have $$\mathbf{P}\big(|H_n^r|\leq \hat{h}_{n,\alpha}^r\big)=O(n^{-2}).$$ Now by Edgeworth expansion theory it is easy to show that under suitable conditions,
\begin{align*}
\mathbf{P}\big(H_n^r\leq x\big)=\Phi(x) + n^{-1/2}q_1^r(x)\phi(x)+n^{-1}q_2^r(x)\phi(x)+n^{-3/2}q_3^r(x)\phi(x)+O(n^{-2})
\end{align*}
uniformly in $x$. Here $q_i^r$ is even or odd polynomial if $i$ is odd or even. Also the coefficients of $q_i^r$ depends on $H_n^r$ through its moments of order $i+2$ or less, $i=1,2,3$. Therefore we have,
\begin{align*}
\mathbf{P}\big(|H_n^r|\leq x\big)=2\Phi(x) -1 +2n^{-1}q_2^r(x)\phi(x)+O(n^{-2})
\end{align*}
uniformly in $x$. If $\Phi(z_\alpha)=1-\alpha/2$, then inverting the above expression we have 
$h_{n,\alpha}^r = z_{\alpha} - n^{-1}q_2^r(z_{\alpha}) + O(n^{-2})$. Similarly, $\hat{h}_{n,\alpha}^r = z_{\alpha} - n^{-1}\hat{q}_2^r(z_{\alpha}) + O_p(n^{-2})$ for some odd polynomial $\hat{q}_2^r$. Hence $\hat{h}_{n,\alpha}^r = h_{n,\alpha}^r-n^{-3/2}V_{n,\alpha}^r + O_p(n^{-2})$ where $V_{n,\alpha}^r = \sqrt{n}\big\{\hat{q}_2^r(z_{\alpha})-q_2^r(z_{\alpha})\big\}$. Therefore upon application of delta method for Edgeworth expansions [cf. Section 2.7 of Hall (1992)], we have
\begin{align}\label{eqn:sym1}
\mathbf{P}\big(\theta_n\in I_{n,(1-\alpha)}^r\big)=\mathbf{P}\big(|H_n^r|\leq \hat{h}_{n,\alpha}^r\big) =& \mathbf{P}\Big(H_n^r+n^{-3/2}V_{n,\alpha}^r \leq h_{n,\alpha}^r\Big)\nonumber\\
&- \mathbf{P}\Big(H_n^r-n^{-3/2}V_{n,\alpha}^r \leq - h_{n,\alpha}^r \Big) + O(n^{-2}),
\end{align}
provided there exist Edgeworth expansions of both $H_n^r+n^{-3/2}V_{n,\alpha}^r$ and $H_n^r-n^{-3/2}V_{n,\alpha}^r$ upto order $n^{-2}$. Now following the arguments similar to Hall (1988), it can be established that
\begin{align*}
&\mathbf{P}\Big(H_n^r+n^{-3/2}V_{n,\alpha}^r\leq x\Big)=\mathbf{P}\Big(H_n^r\leq x\Big) + n^{-3/2}s_1^r(x) + O(n^{-2})\\
& \mathbf{P}\Big(H_n^r-n^{-3/2}V_{n,\alpha}^r\leq x\Big)=\mathbf{P}\Big(H_n^r\leq x\Big) - n^{-3/2}s_1^r(x) + O(n^{-2}),
\end{align*}
both uniformly in $x$, since the quantity $V_{n,\alpha}^r$ is properly centered and scaled. Here $s_1^r$ is an odd polynomial. Using these expansions in (\ref{eqn:sym1}), we have $\mathbf{P}\big(\theta_n\in I_{n,(1-\alpha)}^r\big)=\mathbf{P}\big(|H_n^r|\leq h_{n,\alpha}^r\big)+O(n^{-2})=(1-\alpha) + O(n^{-2})$. 

Suppose $I_{n,(1-\alpha)}^p$ is the symmetric bootstrap interval based on $\big(H_n^p,H_n^{p*}\big)$. Then through the same line of arguments we can conclude that
\begin{align}\label{eqn:sym2}
\mathbf{P}\big(\theta_n\in I_{n,(1-\alpha)}^p\big)=& \mathbf{P}\Big(H_n^p+n^{-3/2}V_{n,\alpha}^p \leq
h_{n,\alpha}^p\Big)\nonumber\\
&- \mathbf{P}\Big(H_n^p-n^{-3/2}V_{n,\alpha}^p \leq - h_{n,\alpha}^p \Big) + O(n^{-2}),
\end{align}
where $q_2^p$, $\hat{q}_2^p$, $h_{n,\alpha}^p$, $V_{n,\alpha}^P$ are analogous quantities corresponding to the perturbation bootstrap. Note that $V_{n,\alpha}^p = \sqrt{n}\big\{\hat{q}_2^p(z_{\alpha})-q_2^p(z_{\alpha})\big\}$ where $\hat{q}_2^p$ is an odd polynomial with coefficients depending on the first four moments of $H_n^p$. Hence by looking into the construction of $H_n^p$, it is clear that the coefficients in $\hat{q}_2^p$ involve first four moments of $\dfrac{(G_1^*-\mu_{G^*})}{\mu_{G^*}}$. Ideally each of theses four moments of $\dfrac{(G_1^*-\mu_{G^*})}{\mu_{G^*}}$ should be equal to $1$ %(\rd{??})
to make $V_{n,\alpha}^p$ properly centered, since $G_i^*$'s do not have any contribution in the polynomial $q_2^p$. Due to the assumption (A.6)(i), the first three moments of $\dfrac{(G_1^*-\mu_{G^*})}{\mu_{G^*}}$ are all equal to 1. However, $E(G_1^*-\mu_{G^*})^4 > \big[E(G_1^*-\mu_{G^*})^2\big]^2=\mu_{G^*}^4$ since $G_i^*$'s are non-degenerate. As a consequence, the fourth moment of $\dfrac{(G_1^*-\mu_{G^*})}{\mu_{G^*}}$ is more than $1$, resulting $V_{n,\alpha}^p$ not properly centered.
Suppose $\tilde{V}_{n,\alpha}^p=\sqrt{n}\big\{\tilde{q}_2^p(z_{\alpha})-q_2^p(z_{\alpha})\big\}$ is properly centered. Hence if we can replace $V_{n,\alpha}^p$ by $\tilde{V}_{n,\alpha}^p$ in the equation (\ref{eqn:sym2}), then we can achieve the error $O(n^{-2})$. Keeping that view in mind, define the corrected symmetric perturbation bootstrap confidence interval of $\theta_n$ as $\tilde{I}_{n,(1-\alpha)}^p$ such that the event $\{\theta_n \in \tilde{I}_{n,(1-\alpha)}^p\}$ is same as the event $\{|H_n^p|\leq \tilde{h}_{n,\alpha}^p\}$ where $\tilde{h}_{n,\alpha}^p=h_{n,\alpha}^p + n^{-1}\big\{\hat{q}_2^p(z_{\alpha})-\tilde{q}_2^p(z_{\alpha})\big\}$. Note that
\begin{align}\label{eqn:sym3}
\mathbf{P}\big(\theta_n\in \tilde{I}_{n,(1-\alpha)}^p\big)=\mathbf{P}\big(|H_n^p|\leq \tilde{h}_{n,\alpha}^p\big)=& \mathbf{P}\Big(H_n^p+n^{-3/2}\tilde{V}_{n,\alpha}^p \leq
h_{n,\alpha}^p\Big)\nonumber\\
&- \mathbf{P}\Big(H_n^p-n^{-3/2}\tilde{V}_{n,\alpha}^p \leq - h_{n,\alpha}^p \Big) + O(n^{-2}),
\end{align}
Now by arguments of Hall (1988),
\begin{align*}
&\mathbf{P}\Big(H_n^p+n^{-3/2}\tilde{V}_{n,\alpha}^p\leq x\Big)=\mathbf{P}\Big(H_n^p\leq x\Big) + n^{-3/2}s_1^p(x) + O(n^{-2})\\
& \mathbf{P}\Big(H_n^p-n^{-3/2}\tilde{V}_{n,\alpha}^r\leq x\Big)=\mathbf{P}\Big(H_n^p\leq x\Big) - n^{-3/2}s_1^p(x) + O(n^{-2}),
\end{align*}
where $s_1^p$ is an odd polynomial. These expansions along with (\ref{eqn:sym3}) imply $\mathbf{P}\big(\theta_n\in \tilde{I}_{n,(1-\alpha)}^p\big)=(1-\alpha) + O(n^{-2})$. Finding the Edgeworth expansions of $H_{n}^p$ and $H_n^{p*}$ we can get hold of the correction factor $C_{n}^p(z_{\alpha})=n^{-1}\big\{\hat{q}_2^p(z_{\alpha})-\tilde{q}_2^p(z_{\alpha})\big\}$. Define $C_{n}^p(\cdot)$ formally as
\begin{align}\label{eqn:sym4}
 C_{n}^p(x)=-n^{-1}x\Big[\dfrac{\omega_2}{2}+\dfrac{\omega_4}{24}(x^2-3)\Big]   
\end{align}
where 
\begin{align*}
\omega_2 = \bigg[\check{\sigma}_n^{-4}\Big[n^{-1}\sum_{i=1}^{n}\hat{\epsilon}_i^4\Big] - \check{\sigma}_n^{-2}\tilde{\Sigma}_n^{-1}\Big[n^{-1}\sum_{i=1}^{n}\big\{\hat{D}_n^{(1)}\hat{C}_{11,n}^{-1}\bm{x}_i^{(1)}\big\}^2\hat{\epsilon}_i^4\Big]\bigg]\bigg[\dfrac{E(G_1^*-\mu_{G^*})^4}{\mu_{G^*}^4}-2\bigg],
\end{align*}
\begin{align*}
\omega_4 =& \tilde{\Sigma}_n^{-2}\Big[n^{-1}\sum_{i=1}^{n}\big\{\hat{D}_n^{(1)}\hat{C}_{11,n}^{-1}\bm{x}_i^{(1)}\big\}^4\hat{\epsilon}_i^4\Big]\bigg[\dfrac{E(G_1^*-\mu_{G^*})^4}{\mu_{G^*}^4}-1\bigg]\\
& + 4\check{\sigma}_n^{-2}\tilde{\Sigma}_n^{-1}\Big[n^{-1}\sum_{i=1}^{n}\big\{\hat{D}_n^{(1)}\hat{C}_{11,n}^{-1}\bm{x}_i^{(1)}\big\}^2\hat{\epsilon}_i^4\Big]\bigg[\dfrac{E(G_1^*-\mu_{G^*})^4}{\mu_{G^*}^4}-2\bigg]\\
& - 3\check{\sigma}_n^{-4}\Big[n^{-1}\sum_{i=1}^{n}\hat{\epsilon}_i^4\Big]\bigg[\dfrac{E(G_1^*-\mu_{G^*})^4}{\mu_{G^*}^4}-2\bigg]\\
& + 1.
\end{align*}
When $G_1^* \sim Beta (1/2,3/2)$, then $\dfrac{E(G_1^*-\mu_{G^*})^4}{\mu_{G^*}^4}=3$ and hence 
\begin{align*}
\omega_2 = \bigg[\check{\sigma}_n^{-4}\Big[n^{-1}\sum_{i=1}^{n}\hat{\epsilon}_i^4\Big] - \check{\sigma}_n^{-2}\tilde{\Sigma}_n^{-1}\Big[n^{-1}\sum_{i=1}^{n}\big\{\hat{D}_n^{(1)}\hat{C}_{11,n}^{-1}\bm{x}_i^{(1)}\big\}^2\hat{\epsilon}_i^4\Big]\bigg] \;\;\; \text{and}
\end{align*}
\begin{align*}
\omega_4 =& 2\tilde{\Sigma}_n^{-2}\Big[n^{-1}\sum_{i=1}^{n}\big\{\hat{D}_n^{(1)}\hat{C}_{11,n}^{-1}\bm{x}_i^{(1)}\big\}^4\hat{\epsilon}_i^4\Big]
 + 4\check{\sigma}_n^{-2}\tilde{\Sigma}_n^{-1}\Big[n^{-1}\sum_{i=1}^{n}\big\{\hat{D}_n^{(1)}\hat{C}_{11,n}^{-1}\bm{x}_i^{(1)}\big\}^2\hat{\epsilon}_i^4\Big]\\
& - 3\check{\sigma}_n^{-4}\Big[n^{-1}\sum_{i=1}^{n}\hat{\epsilon}_i^4\Big]
 + 1.
\end{align*}
We are now ready to formally state the result.

\begin{theo}\label{thm:symboot} First assume that $\hat{\bm{\beta}}_n$ falls in class I.
Suppose $\mathbf{P}(\bm{A}_n\cap \bm{B}_n)=1-o(n^{-2})$, $\mathbf{P}_*(\bm{A}_n^*\cap \bm{B}_n^*)=1-o_p(n^{-2})$ for residual bootstrap and $\mathbf{P}_*(\bm{A}_n^{**}\cap \bm{B}_n^{**})=1-o_p(n^{-2})$ for perturbation bootstrap. Also suppose the conditions (A.3), (A.5) (A.7)(iii) hold with $r=8$.

Again if $\hat{\bm{\beta}}_n$ falls in class II or class III, then respectively the conditions in \ref{thm:bootII} with probabilities $o(n^{-2})$ \& $o_p(n^{-2})$ instead of $o(n^{-1/2})$ \& $o_p(n^{-1/2})$ in (\ref{eqn:initialmoderate}) or assumptions in Theorem \ref{thm:bootIII} hold with $r=8$. 

Then for any $\alpha \in (0,1)$, we have
\begin{align*}
%&\sup_{\delta_7<\alpha<1-\delta_7}\big{|}\mathbf{P}\big(\theta_n\in \tilde{I}_{n,(1-\alpha)}^p\big)-(1-\alpha)\big{|}=O(n^{-2})\;\;\;\text{and}\\
\mathbf{P}\big(\theta_n\in I_{n,(1-\alpha)}^r\big)=1-\alpha +O(n^{-2})\;\; \text{and}\;\;
\mathbf{P}\big(\theta_n\in \tilde{I}_{n,(1-\alpha)}^p\big)=1-\alpha +O(n^{-2})
%&\sup_{\delta_8<\alpha<1-\delta_8}\big{|}\mathbf{P}\big(\theta_n\in I_{n,(1-\alpha)}^r\big)-(1-\alpha)\big{|}=O(n^{-2}),
\end{align*}
provided $p_0^8=o\big(n(\log n)^{-3}\big)$. The intervals $I_{n,(1-\alpha)}^r$ and $\tilde{I}_{n,(1-\alpha)}^p$ are defined earlier in this subsection. Here $\theta_n=\bm{D}_n\bm{\beta}_n$ with $\bm{D}_n$ being a $p-$dimensional row vector. %$\delta_7$ and $\delta_8$ are two constants belonging to $(0,1)$.
\end{theo}

\begin{remark}
All the results in this section are true when $p$ grows with $n$. The rate of growth of $p$ depends critically on how thin the tail of the distribution of the regression error is. If $\mathbf{E}|\epsilon_1|^{2l+3}< \infty$, then $p$ can grow like $n^l$ for any $l>0$, that is, $p$ can grow polynomially. On the otherhand $p$ can grow exponentially if mgf of $\epsilon_1$ exists. For example if $\epsilon_1$ is sub-exponential or sub-gaussian then $p$ can grow like $e^{n^{(\delta_1-\gamma\delta_2)}}$. As far as the perturbation bootstrap is concerned, $G_1^*$ must be sub-gaussian if $\epsilon_1$ is sub-gaussian and $G_1^*$ must be sub-exponential if $\epsilon_1$ is sub-exponential. Hence $Beta(1/2,$ $3/2)$ is an appropriate choice for the distribution of $G_i^*$'s when the errors are sub-gaussian and the distribution of $(M_1+M_2)$ is an appropriate choice for the distribution of $G_i^*$'s when the errors are sub-exponential where $M_1$ and $M_2$ are independent and $M_1$ is a Gamma random variable with shape and scale parameters $0.008652$ and $2$ respectively and $M_2$ is a Beta random variable with both the parameters $0.036490$. For details on the  growth rate of $p$ for the purpose of inference in adaptive lasso, see the section 5.2.2.2 in Das et al. (2019). 
\end{remark}

%\section{Simulation results}
%\label{sec:simulation}

%\section{Data analysis}
%\label{sec:dataanalysis}

\section{Proofs}\label{sec:proofs}
\par
\subsection{Notations}
We denote the true parameter vector as $\bm{\beta}_n = (\beta_{1,n},\dots,\beta_{p,n} )'$, where the subscript $n$ emphasizes that the dimension $p:=p_n$ may grow with the sample size $n$.
Set $\mathcal{A}_n=\{j: \beta_{j,n}\neq 0\}$ and  $p_0:=p_{0,n}=|\mathcal{A}_n|$. For simplicity, we shall suppress the subscript $n$ in the notations $p_n$ and $p_{0n}$. Without loss of generality, we shall assume that $\mathcal{A}_n=\{1,\dots,p_0\}$. Let $\bm{C}_n=n^{-1}\sum_{i=1}^{n}\bm{x}_i\bm{x}'_i$ and partition it according to $\mathcal{A}_n = \{1,\dots,p_0\}$ as 
\begin{equation*}
\bm{C}_{n} = \begin{bmatrix}
\bm{C}_{11,n} \;\;\;\bm{C}_{12,n}\\
\bm{C}_{21,n}\;\;\; \bm{C}_{22,n}
\end{bmatrix},
\end{equation*}
where $\bm{C}_{11,n}$ is of dimension $p_0\times p_0$. Define $\tilde{\bm{x}}_i=\bm{C}_n^{-1}\bm{x}_i$ (when $p\leq n$) and $sgn(x) =-1, 0 ,1$ according as $x<0$, $x=0$, $x>0$, respectively. Suppose $\bm{D}_n$ is a known $q\times p$ matrix with $\text{tr}(\bm{D}_n\bm{D}'_n)=O(1)$ and $q$ is not dependent on $n$. Let $\bm{D}_n^{(1)}$ contains the first $p_0$ columns of $\bm{D}_n$.
%Define $\bm{W}_n=n^{-1/2}\sum_{i=1}^{n}\bm{x}_i\epsilon_i$ and $\bm{W}_n=\big(\bm{W}_n^{(1)'},\bm{W}_n^{(2)'}\big)'=(W_{1n},\ldots,W_{pn})'$, where $\bm{W}_n^{(1)}$ consists of the first $p_0$ components of $\bm{W}_n$. Also define a $p_0\times 1$ vector $\tilde{\bm{s}}_n^{(1)}$ which has $j$th component $\sgn(\beta_{j,n})|\tilde{\beta}_{j,n}|$.  Similarly define the bootstrap versions of $\bm{W}_n$ and $\tilde{\bm{s}}_n^{(1)}$ as $\breve{\bm{W}}_n^*$ and $\tilde{\bm{s}}_n^{*(1)}$  by replacing $\bm{\beta}_n$ and $\tilde{\bm{\beta}}_n$ and $\epsilon_i$ respectively by $\bm{\hat{\beta}}_n$, $\tilde{\bm{\beta}}_n^*$ and $\hat{\epsilon}_i$, where $\hat{\epsilon}_i = Y_i - \x_i'\hat{\bm{\beta}}_n$ are Alasso residuals.
Define
\begin{align*}
\bm{S}_{n} = \begin{bmatrix}
\bm{D}_n^{(1)}\bm{C}_{11,n}^{-1}\bm{D}_n^{(1)'}.\sigma^2 \;\;\;\bm{D}_n^{(1)}\bm{C}_{11,n}^{-1}\bar{\bm{x}}^{(1)}_n.\mu_{3}\\
\bar{\bm{x}}^{(1)'}_n\bm{C}_{11,n}^{-1}\bm{D}_n^{(1)'}.\mu_3\;\;\;\;\;\;\;\;\; (\mu_4-\sigma^4)
\end{bmatrix},
\end{align*}
where $\bar{\bm{x}}_n=n^{-1}\sum_{i=1}^{n}\bm{x}_i=(\bar{\bm{x}}^{(1)\prime}_n,\bar{\bm{x}}^{(2)\prime}_n)^\prime$, $\sigma^2=\mathbf{Var}(\epsilon_1)=\mathbf{E}(\epsilon_1^2)$, and where $\mu_3$ and $\mu_4$ are, respectively, the third and fourth central moments of $\epsilon_1$. Let $K$ be a generic positive constant, independent of $n$ and $p$. By $\mathbf{P_*}$ and $\mathbf{E_*}$ we denote, respectively, probability and expectation with respect to the distribution of $G^{*}$ conditional upon the observed data.

\par

\par
We denote by $\|\cdot\|$ and $\|\cdot\|_{\infty}$, respectively, the $L^2$ and $L^{\infty}$ norm. 
For a non-negative integer-valued vector $\bm{\alpha} = (\alpha_1, \alpha_2,\ldots,\alpha_l)'$ and a function $f = (f_1,f_2,\ldots,f_l):\ \mathscr{R}^l\ \rightarrow \ \mathscr{R}^l$, $l\geq 1$, write $|\bm{\alpha}| = \alpha_1 +\ldots+ \alpha_l$, $\bm{\alpha}! = \alpha_1!\ldots \alpha_l!$, $f^{\bm{\alpha}} = (f_1^{\alpha_1})\ldots(f_l^{\alpha_l})$, and $D^{\bm{\alpha}}f_1 = D_1^{\alpha_1}\ldots D_l^{\alpha_l}f_1$, where $D_jf_1$ denotes the partial derivative of $f_1$ with respect to the $j$th component of the argument, $1\leq j \leq l$. 
For $\bm{t} =(t_1,\ldots t_l)'\in \mathscr{R}^l$ and $\bm{\alpha}$ as above, define $t^{\bm{\alpha}} = t_1^{\alpha_1}\ldots t_l^{\alpha_l}$. Let $\bm{\Phi}_V$ denote the multivariate Normal distribution with mean $\mathbf{0}$ and dispersion matrix $\bm{V}$ having $j$th row $\bm{V}_{j.}$ and let $\phi_V$ denote the density of $\bm{\Phi}_V$. We write $\bm{\Phi}_V= \bm{\Phi}$ and $\phi_V = \phi$ when $\bm{V}$ is the identity matrix. Also define the polynomial $\chi_{\bm{\alpha}}(\bm{y}:\bm{V})$ by the identity $(-D)^{\bm{\alpha}}\phi(\bm{y}:\bm{V})=\chi_{\bm{\alpha}}(\bm{y}:\bm{V})\phi(\bm{y}:\bm{V})$. We write $\chi_{\bm{\alpha}}(\bm{y}:\bm{V})$ as $\chi_{\bm{\alpha}}(\bm{y})$ when $\bm{V}$ is the identity matrix. For any set $B\subseteq \mathcal{R}^p$ and any $\bm{b}\in \mathcal{R}^p$, $B+\bm{b}=\{\bm{a}+\bm{b}:\bm{a}\in B\}$.

\par
 Define, $\bm{W}_n=n^{-1/2}\sum_{i=1}^{n}\bm{x}_i\epsilon_i$, $\bm{W}_n^{*}=n^{-1/2}\sum_{i=1}^{n}\epsilon_i^*\bm{x}_i$ when underlying bootstrap is residual bootstrap and $\bm{W}_n^{*}=n^{-1/2}\sum_{i=1}^{n}\hat{\epsilon}_i\bm{x}_i$ $(G_i^*-\mu_{G^*})$ when the perturbation bootstrap is considered. 
Now define the sets $\bm{A}_{2n}=$ $\Big{\{}\big{\{}\|\bm{W}_n\|_{\infty}\leq K\sqrt{\log n}\big{\}}\cap \big{\{} ||\sqrt{n}\big(\tilde{\bm{\beta}}-\bm{\beta}\big)||_{\infty}\leq Cn^{\delta_2}\big{\}}\Big{\}}$ and $\bm{A}_{3n}=$ $\Big{\{}\|\bm{W}_n\|_{\infty}\leq K\sqrt{\log n}\Big\}$. Similarly define the bootstrap sets $\bm{A}_{2n}^*=$ $\Big{\{}\big{\{}\|\bm{W}_n^*\|_{\infty}\leq K\sqrt{\log n}\big{\}}\cap \big{\{} \|\sqrt{n}\big(\tilde{\bm{\beta}}_n^*-\hat{\bm{\beta}}_n\big)\|_{\infty}\leq Cn^{\delta_2}\big{\}}\Big{\}}$ and $\bm{A}_{3n}^*=$ $\Big{\{}\|\bm{W}_n^*\|_{\infty}\leq K\sqrt{\log n}\Big\}$. Arguments similar to the proofs of Lemma 8.1 of Chatterjee and Lahiri (2013) and Das et al. (2018) imply that $\mathbf{P}(\bm{A}_{2n}),\mathbf{P}(\bm{A}_{3n}) \geq 1-o_p(n^{-1/2})$ and $\mathbf{P}_*(\bm{A}_{2n}^*),\mathbf{P}_*(\bm{A}_{3n}^*) \geq 1-o_p(n^{-1/2})$, under the assumptions of Theorem \ref{thm:boot}.

%Note that, $\check{\bm{b}}_n=O_p(n^{-\delta_1})$, by Lemma \ref{lem:betahat} and \ref{lem:Sigma}, described below. Suppose, $r_1=\min\{a\in \mathcal{N}:||\check{\bm{b}}_n||^{a+1}=o_p(n^{-1/2})\}$, $\mathcal{N}$ being the set of natural numbers. Define the conditional Lebesgue density of two-term Edgeworth expansion of $\bm{R}_{n}^*$ as
%\begin{align*}
%\xi_n^*(\bm{x})=&\phi(\bm{x})\Bigg[1+\sum_{k=1}^{r_1}\dfrac{1}{k!}\big{\{}\sum_{|\bm{\alpha}|=k}\check{\bm{b}}_n^{\alpha}H_{\bm{\alpha}}(\bm{x})\big{\}}+\dfrac{1}{\sqrt{n}}\bigg[ \dfrac{1}{6}\sum_{|\bm{\alpha}|=3}\bm{t}^{\bm{\alpha}}\bar{\bm{\xi}}_n^{*(1)}({\bm{\alpha}})H_{\bm{\alpha}}(\bm{x})\\ 
%&-\dfrac{1}{2\hat{\sigma}_n^2}\Big{\{}\sum_{|\bm{\alpha}|=1}\bm{t}^{\bm{\alpha}}\bar{\bm{\xi}}_n^{*(3)}({\bm{\alpha}})H_{\bm{\alpha}}(\bm{x})
%+\sum_{|\bm{\alpha}|=1}\sum_{|\bm{\zeta}|=2}\bm{t}^{\bm{\alpha}+\bm{\zeta}}\bar{\bm{\xi}}_n^{*(3)}({\bm{\alpha}})\bar{\bm{\xi}}_n^{*(1)}({\bm{\zeta}})H_{\bm{\alpha}+\bm{\zeta}}(\bm{x})\Big{\}}\bigg]\Bigg],
%\end{align*}
%where $x\in \mathcal{R}^q$, $\bar{\bm{\xi}}_{n}^{*(j)}(\bm{\alpha})=n^{-1}\sum_{i=1}^{n}\Big(\check{\bm{\xi}}_i^{(0)}\hat{\epsilon}_i^j\Big)^{\bm{\alpha}}$, $j=0,1,\ldots$ and $H_{\bm{\alpha}}(\bm{x})=(-D)^{\bm{\alpha}}\phi(\bm{x})$, where $\phi(\bm{\cdot})$ is the standard normal density on $\mathcal{R}^q$. 

\subsection{Preliminary Lemmas} Lemmas necessary for the proofs of the results, are stated in this section, along with their proofs.

\begin{lem}\label{lem:concentration}
Suppose $Y_1,\dots,Y_n$ are zero mean independent r.v.s and $\mathbf{E}(|Y_i|^t)< \infty$ for $i = 1,\dots,n$ and $\sum_{i = 1}^{n}\mathbf{E}(|Y_i|^t) = \sigma_t$; $S_n = \sum_{i = 1}^{n}Y_i$. Then, for any $t\geq 2$ and $x>0$
\begin{equation*}
P[|S_n|>x]\leq C[\sigma_t x^{-t} + exp(-x^2/\sigma_2)] 
\end{equation*}
\end{lem}
%Proof of Lemma 8.1. 
\textbf{Proof of Lemma \ref{lem:concentration}}.
This inequality was proved in Fuk and Nagaev (1971).

 \begin{lem}\label{lem:betainitialOLS}
 Suppose $p\leq n$ and the initial estimator $\tilde{\bm{\beta}}_n$ is the OLS estimator. Then under conditions (A.1)(iv), (A.5)(i) and (A.6)(i) with r=2, we have 
 \begin{align*}
 &\mathbf{P}\Big(\|\tilde{\bm{\beta}}_n - \bm{\beta}_n\|=O\big{(}n^{-1/2}(\log n)^{1/2}\big{)}\Big) \geq 1 - o\big{(}n^{-1/2}\big{)}\\
&\mathbf{P_*}\Big(\|\tilde{\bm{\beta}}_n^{*} - \hat{\bm{\beta}}_n||=O\big{(}n^{-1/2}(\log n)^{1/2}\big{)}\Big) \geq 1 - o_p\big{(}n^{-1/2}\big{)}
\end{align*} 
that is (\ref{eqn:initialmoderate}) holds.
 \end{lem}

\textbf{Proof of Lemma \ref{lem:betainitialOLS}}.
This lemma follows through the similar argument as in part (iii) of Lemma 8.1 of Chatterjee and Lahiri (2013).

 \begin{lem}\label{lem:betahat}
Suppose $\hat{\bm{\beta}}_n$ belongs to class I or class II. Then under the respective assumptions of Theorem \ref{thm:boot}, we have
 \begin{align*}
 \|\bm{\hat{\beta}}_n-\bm{\beta}_n\|_{\infty}=O_p(n^{-1/2}), \; \;\; \|\bm{\hat{\beta}}_n^*-\hat{\bm{\beta}}_n\|_{\infty}=O_{p_*}(n^{-1/2})\; \text{and}\; \|\bm{\hat{\beta}}_n^{**}-\hat{\bm{\beta}}_n\|_{\infty}=O_{p_*}(n^{-1/2}).
 \end{align*}
 \end{lem}
 
\textbf{Proof of Lemma \ref{lem:betahat}}.
This lemma follows through the same line of arguments as in the proof of Lemma 4 of Das et al. (2018).

\begin{lem}\label{lem:betahathat}
Suppose $\hat{\bm{\beta}}_n$ is the Lasso estimator of $\bm{\beta}$. $\bm{\beta}_n^{(1)}$ and $\hat{\bm{\beta}}_n^{(1)}$ respectively consist of the non-zero components of $\bm{\beta}_n$ and the non-zero components of $\hat{\bm{\beta}}_n$. $\hat{\bm{C}}_{11,n}$ is the submatrix of $\bm{C}_n$ as defined in Section \ref{sec:mainresults}. Then under the conditions of part 3. of Theorem \ref{thm:boot}, we have
 \begin{align*}
 &\mathbf{P}\Big(\|\hat{\bm{\beta}}_n - \bm{\beta}_n\|=O\big{(}p_0^{3/2}n^{-1}\lambda_n\big{)}\Big) \geq 1 - o\big{(}n^{-1/2}\big{)}\\
&\mathbf{P}\Big(\|\hat{\bm{C}}_{11,n}\big(\hat{\bm{\beta}}_n^{(1)} - \bm{\beta}_n^{(1)}\big)\|=O\big{(}p_0^{1/2}n^{-1}\lambda_n\big{)}\Big) \geq 1 - o\big{(}n^{-1/2}\big{)}%\\
%&\mathbf{P_*}\Big(\|\hat{\bm{\beta}}_n^* - \hat{\bm{\beta}}_n^*\|=O\big{(}p_0^{3/2}n^{-1}\lambda_n\big{)}\Big) \geq 1 - o\big{(}n^{-1/2}\big{)}\\
%&\mathbf{P_*}\Big(\|\hat{\bm{C}}_{11,n}\big(\hat{\bm{\beta}}_n^{*(1)} - \hat{\bm{\beta}}_n^{(1)}\big)\|=O\big{(}p_0^{1/2}n^{-1}\lambda_n\big{)}\Big) \geq 1 - o\big{(}n^{-1/2}\big{)}
\end{align*} 
 \end{lem}
 
\textbf{Proof of Lemma \ref{lem:betahathat}}.
In the proof of Theorem \ref{thm:oracleIII}, it is shown that $\big(\hat{\bm{\beta}}_n - \bm{\beta}_n\big)$ has only first $p_0$ components non-zero. $\big(\hat{\bm{\beta}}_n^{(1)} - \bm{\beta}_n^{(1)}\big)$ consists of these non-zero components for sufficiently large $n$ and $$\sqrt{n}\big(\hat{\bm{\beta}}_n^{(1)} - \bm{\beta}_n^{(1)}\big)=\bm{C}_{11,n}^{-1}\Big[\bm{W}_n^{(1)}-\dfrac{\lambda_n}{2\sqrt{n}}\bm{c}_n^{(1)}\Big],$$
where $\bm{W}_n^{(1)}=n^{-1/2}\sum_{i=1}^{n}\bm{x}_i^{(1)}\epsilon_i$ and $\bm{c}_{n}^{(1)}=(c_{1n},\ldots, c_{p_0n})$ with $c_{j,n}=sgn(\beta_{j,n})$. Now the lemma follows easily due to lemma \ref{lem:concentration} and noting that $\log n=O(n^{-1}\lambda_n^2)$ [cf. condition (A.7) (ii)$^\prime$].

%This lemma follows from Markov inequality and using the condition $\max\{n^{-1}\sum_{i=1}^{n}\big|(\bm{C}_{11,n}^{-1})_{j.}\bm{x}_{i}^{(1)}\big|^{2r}:1\leq j \leq p_0\} = O(1)$ [stated in assumption (A.1)(ii)], after observing the form of the Alasso estimator obtained in (8.5) of Theorem 8.2 (a) of Chatterjee and Lahiri (2013) and the solution $\hat{\bm{u}}_n^*$ of the equation \ref{eqn:KKTinproofs} obtained in the proof of part (a) of Lemma \ref{lem:Rstarconverge}.
 
\begin{lem}\label{lem:Sigma} The matrices $\Sigma_n$, $\tilde{\Sigma}_n$ are defined in Section \ref{sec:mainresults}. The matrix $\bm{\Upsilon}_n$ is defined in the proof of Theorem \ref{thm:oracleII}. The matrices $\tilde{\bm{\Upsilon}}_n$, $\hat{\bm{\Upsilon}}_n$ are defined in the proof of Theorem \ref{thm:boot}. 
Under the assumptions of Theorem \ref{thm:boot}, we have
\begin{align*}
||\hat{\bm{\Upsilon}}_n-\bm{\Upsilon}_n||=o_p(n^{-1/2}),
||\tilde{\bm{\Sigma}}_n-\sigma^2\bm{\Sigma}_n||=O_p(n^{-1/2}), ||\tilde{\bm{\Upsilon}}_n-\sigma^2\bm{\Upsilon}_n||=O_p(n^{-1/2}),
\end{align*}
where $\delta_1$ is as defined in assumption \emph{(A.6)}.
\end{lem}
 
\textbf{Proof of Lemma \ref{lem:Sigma}}. 
First we show that 
 $||\hat{\bm{\Upsilon}}_n-\bm{\Upsilon}_n||=o_p(n^{-1/2})$. Note that 
 \begin{align*}
 \hat{\bm{\Upsilon}}_n-\bm{\Upsilon}_n=&n^{-1}\sum_{i=1}^{n}\big(\hat{\bm{\eta}}_i^{(0)}-\bm{\eta}_i^{(0)}\big)\big(\bm{\xi}_i^{(0)}+\bm{\eta}_i^{(0)}\big)'+n^{-1}\sum_{i=1}^{n}\big(\bm{\xi}_i^{(0)}+\hat{\bm{\eta}}_i^{(0)}\big)'\big(\hat{\bm{\eta}}_i^{(0)}-\bm{\eta}_i^{(0)}\big)'
 \end{align*}
 where, due to conditions (A.2)(iv), (A.3)(i), (A.7)(i) and Lemma \ref{lem:betahat}, for large enough $n$ we have
 \begin{align*}
 &n^{-1}\sum_{i=1}^{n}||{\hat{\bm{\eta}}}_i^{(0)}-\bm{\eta}_i^{(0)}||^2\\
 &\leq K^2(\gamma)\cdot \|\bm{D}_n^{(1)}\bm{C}_{11,n}^{-1/2}\|^2\cdot \|\bm{C}_{11,n}^{-1/2}\|^2\cdot\lambda_n^2\Big(\operatorname*{\max}_{1\leq j \leq p}n^{-1}\sum_{i=1}^{n}\big|\tilde{x}_{i,j}\big|^2\Big)\cdot\|\bm{\hat{\beta}}_n^{(1)}-\bm{\beta}_n^{(1)}\|^2\cdot \tilde{P}_3^2\\
 & =o_p\big(n^{-1}\big)
 \end{align*}
 and 
 \begin{align*}
 &n^{-1}\sum_{i=1}^{n}||\bm{\xi}_i^{(0)}||^2 + n^{-1}\sum_{i=1}^{n}||\bm{\eta}_i^{(0)}||^2 +n^{-1}\sum_{i=1}^{n}||\hat{\bm{\eta}}_i^{(0)}||^2\\
 & \leq tr\big(\bm{D}_n^{(1)}\bm{C}_{11,n}^{-1}\bm{D}_n^{(1)}\big)+K\cdot \|\bm{D}_n^{(1)}\bm{C}_{11,n}^{-1}\|^2\cdot \lambda_n^2\Big(\operatorname*{\max}_{1\leq j \leq p}n^{-1}\sum_{i=1}^{n}\big|\tilde{x}_{i,j}\big|^2\Big)\cdot p_0\cdot O_p(\tilde{P}_2+\tilde{P}_3/\sqrt{n})\\
 & =O(1),
 \end{align*}
 since $\|\bm{C}_{11,n}^{-1/2}||^2\leq K \min\{p_0, n^a\}$ and $||\bm{D}_n^{(1)}\bm{C}_{11,n}^{-1/2}\|^2\leq q\|\bm{D}_n^{(1)}\bm{C}_{11,n}^{-1}\bm{D}_n^{(1)'}\|
 $.
 
 Therefore, by the Cauchy-Schwarz inequality, we have $||\hat{\bm{\Upsilon}}_n-\bm{\Upsilon}_n||=o_p(n^{-1/2})$.
 \par
 Now to prove other two, note that for large enough $n$,
  \begin{align*}
 \tilde{\bm{\Sigma}}_n-\sigma^2\bar{\bm{\Sigma}}_n=& n^{-1}\sum_{i=1}^{n}\bm{\xi}_i^{(0)}\bm{\xi}_i^{(0)'}(\hat{\epsilon}_i^2-\sigma^2)
 \end{align*}
 and
 \begin{align*}
 \tilde{\bm{\Upsilon}}_n-\sigma^2\bm{\Upsilon}_n=& n^{-1}\sum_{i=1}^{n}\bm{\xi}_i^{(0)}\bm{\xi}_i^{(0)'}(\hat{\epsilon}_i^2-\sigma^2) + n^{-1}\sum_{i=1}^{n}\big(\hat{\bm{\eta}}_i^{(0)}-\bm{\eta}_i^{(0)}\big)\bm{\xi}_i^{(0)'}\hat{\epsilon}_i^2\\ &+n^{-1}\sum_{i=1}^{n}\bm{\eta}_i^{(0)}\bm{\xi}_i^{(0)'}(\hat{\epsilon}_i^2-\sigma^2)+n^{-1}\sum_{i=1}^{n}\bm{\xi}_i^{(0)}\big(\hat{\bm{\eta}}_i^{(0)}-\bm{\eta}_i^{(0)}\big)'\hat{\epsilon}_i^2\\ &+n^{-1}\sum_{i=1}^{n}\bm{\xi}_i^{(0)}\bm{\eta}_i^{(0)'}(\hat{\epsilon}_i^2-\sigma^2) + n^{-1}\sum_{i=1}^{n}\hat{\bm{\eta}}_i^{(0)}\big(\hat{\bm{\eta}}_i^{(0)}-\bm{\eta}_i^{(0)}\big)'\hat{\epsilon}_i^2\\
 & +n^{-1}\sum_{i=1}^{n}\big(\hat{\bm{\eta}}_i^{(0)}-\bm{\eta}_i^{(0)}\big)\bm{\eta}_i^{(0)'}\hat{\epsilon}_i^2+n^{-1}\sum_{i=1}^{n}\bm{\eta}_i^{(0)}\bm{\eta}_i^{(0)'}(\hat{\epsilon}_i^2-\sigma^2).
 \end{align*}
 Now we need to find the order of the term $||n^{-1}\sum_{i=1}^{n}\bm{\xi}_i^{(0)}\bm{\xi}_i^{(0)'}(\hat{\epsilon}_i^2-\sigma^2)||$ or the order of $\|\tilde{\bm{\Sigma}}_n-\sigma^2\bm{\Sigma}_n\|$, to find the order of $\|\tilde{\bm{\Upsilon}}_n-\sigma^2\bm{\Upsilon}_n\|$, since other terms can be shown to be of smaller order by using H\"older's inequality. Note that by Lemma \ref{lem:concentration}, Lemma \ref{lem:betahat}, Lemma \ref{lem:betahathat} and conditions (A.2)(ii), (A.3)(i), (A.3)(ii) \& (A.5)(i), we have
 \begin{align*}
 \mathbf{P}\Big(\Big{\{}\Big\|\sum_{i=1}^{n}\bm{\xi}_i^{(0)}\bm{\xi}_i^{(0)'}(\epsilon_i^2-\sigma^2)\Big\|>K. n^{1/2}\Big{\}} \Big)\rightarrow 0\; \text{as}\; K\rightarrow \infty
 \end{align*}
 and
  \begin{align*}
 \mathbf{P}\Big(\Big{\{}\Big\|\sum_{i=1}^{n}\bm{\xi}_i^{(0)}\bm{\xi}_i^{(0)'}(\hat{\epsilon}_i^2-\epsilon_i^2)\Big\|>K.n^{1/2}\Big{\}}\Big)\rightarrow 0\; \text{as}\; K\rightarrow \infty
  \end{align*}
Therefore Lemma \ref{lem:Sigma} follows. 

\begin{lem}\label{lem:Normapprox}
Suppose $\bm{X}_1,\dots,\bm{X}_n$ are $n$ random vectors in $\mathcal{R}^k$ satisfying $\mathbf{E}\bm{X}_j=0$, for all $j\in\{1,\dots,n\}$, and $\bm{V}_n=n^{-1}\sum_{i=1}^{n}\mathbf{Var}(\bm{X}_j)$ where $\bm{V}_n$ is a positive definite matrix. Suppose $\gamma_n$ is the smallest eigen value of $\bm{V}_n$. Define, $\rho_3=n^{-1}\sum_{i=1}^{n}\mathbf{E}||\bm{X}_j||^3$. If $\rho_3<\infty$, then we have
\begin{align*}
\sup_{\bm{C}\in \mathcal{C}_q}\Big|\mathbf{P}\big(n^{-1/2}\sum_{i=1}^{n}\bm{X}_i\in \bm{C}\big)-\Phi\big(\bm{C}; \sigma^2\bm{V}_n\big)\Big|\leq k\gamma_n^{-3/2}\rho_3n^{-1/2}
\end{align*}
\end{lem}
\textbf{Proof of Lemma \ref{lem:Normapprox}}.
This result is stated as Corollary 17.2 in Bhattacharya and Rao (1986).

\begin{lem}\label{lem:residualCramer}
Suppose for class I, $\mathbf{P}(\bm{A}_n\cap \bm{B}_n)=1-o(n^{-1/2})$ is true and conditions (A.2)(ii), (A.5) holds with $r=3$. For class II, suppose all the assumptions corresponding to residual bootstrap of Theorem 6 except (A.3) hold. For Lasso, which belongs to class III, suppose all the assumptions of Theorem 7 except (A.3) hold. Then for any $\delta_2>0$ and $K\in(0,\infty)$, there exists $\delta_3 \in (0,1)$ such that 
$$\sup\Big\{\hat{\omega}_n(t_1,t_2):\delta_2^2\leq t_1^2+t_2^2\leq n^K\Big\}=1-\delta_3+o_p(1),$$
where $\hat{\omega}_n(t_1,t_2)=\mathbf{E}_*\exp\Big(it_1\epsilon_1^*+it_2\epsilon_1^{*2}\Big)$.
\end{lem}
\textbf{Proof of Lemma \ref{lem:residualCramer}}.
This follows through the same line arguments as in the proof of Lemma 2 in Babu and Singh (1984).

\subsection{Proof of Results}

\textbf{Proof of Theorem 1:} Note that on the set $\bm{A}_n\cap \bm{B}_n$,
\begin{align*}
\bm{T}_n=\sqrt{n}\bm{D}_n(\hat{\bm{\beta}}_n - \bm{\beta}) = n^{-1/2}\sum_{i=1}^{n}\bm{D}_n^{(1)}\bm{C}_{11,n}^{-1}\bm{x}_i^{(1)}\epsilon_i = n^{-1/2}\sum_{i=1}^{n}\bm{\xi}_i^{(0)}\epsilon_i.
\end{align*}

Now note that under the conditions (A.3)(ii), (A.3)(iii), (A.5)(i) with $r=3/2$, we can employ Lemma \ref{lem:Normapprox} to obtain
\begin{align*}
\sup_{\bm{B}\in \mathcal{C}_q}\Big|\mathbf{P}\big(n^{-1/2}\sum_{i=1}^{n}\bm{\xi}_i^{(0)}\epsilon_i\in \bm{B}\big)-\Phi\big(\bm{B}; \sigma^2\bm{\Sigma}_n\big)\Big|=O\big(n^{-1/2}\big),
\end{align*}
where $\bm{\Sigma}_n=n^{-1}\sum_{i=1}^{n}\bm{\xi}_i^{(0)}\bm{\xi}_i^{(0)\prime}=\bm{D}_n^{(1)}\bm{C}_{11,n}^{-1}\bm{D}_n^{(1)}$. Now since $\mathbf{P}(\bm{A}_n\cap \bm{B}_n)=1-o(n^{-1/2})$, therefore we have

%The first three cumulants of $\bm{t}'\bm{T}_{n}$ are given by

%$\kappa_1\big(\bm{t}'\bm{T}_{n}\big)=\mathbf{E}\big(\bm{t}'\bm{T}_{n}\big)=0$

%$\kappa_2\big(\bm{t}'\bm{T}_{n}\big)=\mathbf{Var}\big(\bm{t}'\bm{T}_{n}\big)=\sigma^2\bm{t}^\prime\Big(n^{-1}\sum_{i=1}^{n}\xi_i^{(0)}\xi_i^{(0)\prime}\Big)\bm{t}=\sigma^2\bm{t}'\Sigma_n\bm{t}$

%$\kappa_3\big(\bm{t}'\bm{T}_{n}\big)=\mathbf{E}\big(\bm{t}'\bm{T}_{n}\big)^3-3\mathbf{E}\big(\bm{t}'\bm{T}_{n}\big)^2. \mathbf{E}\big(\bm{t}'\bm{T}_{n}\big)+2\Big(\mathbf{E}\big(\bm{t}'\bm{T}_{n}\big)\Big)^3= \dfrac{\mu_3}{\sqrt{n}}\sum_{|\bm{\alpha}|=3}\bm{t}^{\bm{\alpha}}\big(n^{-1}\sum_{i=1}^{n}(\bm{\xi}_i^{(0)})^{\bm{\alpha}}\big)$.\\

%Therefore the Lebesgue density of two term EE of $\bm{T}_n$ on the set $\bm{A}_n\cap \bm{B}_n$ is given by

%\begin{align}\label{eqn:sole}
%\bm{\Psi}_n(\bm{x})=\Big[1-\dfrac{\mu_3}{\sqrt{n}}\sum_{|\bm{\alpha}|=3}\big(n^{-1}\sum_{i=1}^{n}(\bm{\xi}_i^{(0)})^{\bm{\alpha}}\big)D^{\bm{\alpha}}\Big]\phi(\bm{x}:\sigma^2\bm{\Sigma}_n)
%\end{align}

%\begin{align*}
%\sup_{\bm{B}\in \mathcal{C}_q}|\mathbf{P}(\bm{T}_n\in \bm{B})-\int_{\bm{B}} \Psi(\bm{x})d\bm{x}|=o(n^{-1/2})
%\end{align*}

\begin{align*}
\Delta_n & = \sup_{\bm{B}\in \mathcal{C}_q}|\mathbf{P}(\bm{T}_n\in \bm{B})-\Phi(\bm{B}; \sigma^2\bm{\Sigma}_n)|\\
&\leq \sup_{\bm{B}\in \mathcal{C}_q}\Big|\mathbf{P}\big(n^{-1/2}\sum_{i=1}^{n}\bm{\xi}_i^{(0)}\epsilon_i\in \bm{B}\big)-\Phi\big(\bm{B}; \sigma^2\bm{\Sigma}_n\big)\Big|+1-\mathbf{P}(\bm{A}_n\cap \bm{B}_n)\\
%& \leq \sup_{\bm{B}\in \mathcal{C}_q}|\mathbf{P}(\bm{T}_n\in \bm{B})-\int_{\bm{B}} \Psi(\bm{x})d\bm{x}| + \sup_{\bm{B}\in \mathcal{C}_q}|\Phi(\bm{B}; \sigma^2\bm{\Sigma}_n)-\int_{\bm{B}} \Psi(\bm{x})d\bm{x}|\\
%& \leq \dfrac{\mu_3}{\sqrt{n}}\sum_{|\bm{\alpha}|=3}n^{-1}\sum_{i=1}^{n}\int_{\bm{B}}\big|\big(\bm{\xi}_i^{(0)}\big)^{\bm{\alpha}}D^{\bm{\alpha}}\phi(\bm{x}:\sigma^2\Sigma_n)\big|d\bm{x} +o(n^{-1/2})\\
& = O(n^{-1/2}).
\end{align*}

\textbf{Proof of Theorem \ref{thm:oracleII} :}
Note that under the set up of Theorem \ref{thm:oracleII}, the penalized estimator $\hat{\bm{\beta}}_n$ is defined as 
 \begin{equation*}
\bm{\hat{\beta}}_n = \operatorname*{arg\,min}_{\bm{t}}\Bigg[\sum_{i=1}^{n}(y_i - \bm{x}'_i \bm{t})^2
+n\lambda_n\sum_{j=1}^{p}\tilde{P}^\prime(|\tilde{\beta}_{j,n}|)(|t_j|)\Bigg]
\end{equation*}

Now, writing $\bm{\hat{u}}_{2n}=\sqrt{n}\big(\bm{\hat{\beta}}_n-\bm{\beta}_n\big)$ and $\bm{W}_n=n^{-1/2}\sum_{i=1}^{n}\bm{x}_i\epsilon_i$, we have
 
 \begin{align}
\bm{\hat{u}}_{2n}&= \operatorname*{arg\,min}_{\bm{v}}\Bigg[\bm{v}^{\prime}\bm{C}_n\bm{v} -2\bm{v}^{\prime}\bm{W}_n+n\lambda_n\sum_{j=1}^{p}\tilde{P}^\prime(|\tilde{\beta}_{j,n}|)\Big(|\beta_{j,n}+\dfrac{v_{j}}{\sqrt{n}}|-|\beta_{j,n}|\Big)\Bigg]\nonumber\\\label{eqn:Z2n}
&= \operatorname*{arg\,min}_{\bm{v}} \bm{Z}_{2n}(\bm{v})\;\;\;\;\; \text{(say)}.
\end{align}

Note that $\bm{Z}_{2n}(\bm{v})$ is convex in $\bm{v}$. Hence, the KKT condition is necessary and sufficient. The KKT condition corresponding to (\ref{eqn:Z2n}) is given by
\begin{align}\label{eqn:KKTinproofs2}
2\bm{C}_n\bm{v}-2\bm{W}_n+\sqrt{n}\lambda_n\bm{\Gamma}_n\bm{l}_n=\bm{0}
\end{align}
for some $l_{j,n}\in [-1,1]$ for all $j\in\{1,\ldots,p\}$, where $\bm{l}_n=(l_{1,n},\ldots, l_{p,n})'$ and $\bm{\Gamma}_n=diag\big(\tilde{P}^\prime(|\tilde{\beta}_{1,n}|),$ $\ldots, \tilde{P}^\prime(|\tilde{\beta}_{p_0,n}|)\big)$. It is easy to show that under the conditions stated in theorem \ref{thm:oracleII}, on the set $\bm{A}_{2n}$, $\Big(\big(\bm{\hat{u}}_{2n}^{(1)}\big)', \bm{0}'\Big)'$, where $\bm{\hat{u}}_{2n}^{(1)}=\bm{C}_{11,n}^{-1}\Big[\bm{W}_n^{(1)}-\dfrac{\sqrt{n}\lambda_n}{2}\tilde{\bm{s}}_n^{(1)}\Big]$ is the unique solution of (\ref{eqn:KKTinproofs2}) and $\bm{W}_n^{(1)}$ consists of the first $p_0$ components of $\bm{W}_n$. Hence $\bm{\hat{u}}_n=\Big(\big(\bm{\hat{u}}_{2n}^{(1)}\big)', \bm{0}'\Big)'$, is the unique solution of the minimization problem (\ref{eqn:Z2n}), where $\tilde{\bm{s}}_n^{(1)}=(\tilde{s}_{1n},\ldots, \tilde{s}_{p_0n})$ and $\tilde{s}_{j,n}=sgn(\beta_{j,n})\tilde{P}^\prime(|\tilde{\beta}_{j,n}|)$ on the set $A_{2n}$. Therefore on the set $A_{2n}$ we have
\begin{align*}
\bm{T}_n&=\sqrt{n}\bm{D}_n(\hat{\bm{\beta}}_n - \bm{\beta}_n)\\
&=\bm{D}_n^{(1)}\bm{C}_{11,n}^{-1}\Big[\bm{W}_n^{(1)}-\dfrac{\sqrt{n}\lambda_n}{2}\tilde{\bm{s}}_n^{(1)}\Big]
\end{align*}
When $p<n$ and $\tilde{\bm{\beta}}_n$ is the OLS, then we have
\begin{align*}
\bm{T}_n &= \bm{D}_n^{(1)}\bm{C}_{11,n}^{-1}\Big[\bm{W}_n^{(1)}-\dfrac{\sqrt{n}\lambda_n}{2}\bm{s}_n^{(1)}-\dfrac{\lambda_n}{2}\bm{L}_n^{(1)}\Big] +\bm{Q}_{1n}\\
&= \bm{T}_{1n}+\bm{Q}_{1n}\;\;\;\; (say)
\end{align*}
where $\bm{s}_n^{(1)}=(s_{1,n},\ldots, s_{p_0,n})$, $s_{j,n}=sgn(\beta_{j,n})\tilde{P}^\prime(|\beta_{j,n}|)$ and $\bm{L}_n^{(1)} = (L_{1,n},\ldots, L_{p_0,n})$ with 
\begin{align*}
L_{j,n}&= \sqrt{n}(\tilde{\beta}_{j,n}-\beta_{j,n})sgn(\beta_{j,n})\tilde{P}^{\prime\prime}(|\beta_{j,n}|)\\
&= n^{-1/2}\sum_{i=1}^{n}\big[\tilde{x}_{i,j}\epsilon_i\big]sgn(\beta_{j,n})\tilde{P}^{\prime\prime}(|\beta_{j,n}|).
\end{align*}
$\bm{Q}_{1n}$ is the remainder term. By condition (A.7)(i) and the continuity of %\rd{$\max_j\tilde{P}^{\prime\prime\prime}(t_j)$ at $\bm{t}=(t_1,\dots,t_{p_0})^\prime=\Big(|\beta_{1,n}|,\dots,|\beta_{p_0,n}|\Big)^\prime$}
$\max_j\tilde{P}^{\prime\prime\prime}(t_j)$ at $\bm{t}=(t_1,\dots,t_{p_0})^\prime=\Big(|\beta_{1,n}|,\dots,|\beta_{p_0,n}|\Big)^\prime$, we have $$P(||\bm{Q}_{1n}||=o(n^{-1/2}))=1-o(n^{-1/2}).$$
Therefore,
\begin{align*}
\bm{T}_{1n}= n^{-1/2}\sum_{i=1}^{n}\big(\bm{\xi}_i^{(0)}+\bm{\eta}_i^{(0)}\big)\epsilon_i+ \bm{b}_n = \bm{T}_{2n}+ \bm{b}_n \;\;\; (say),
\end{align*}
where
$$\bm{\xi}_i^{(0)}=\bm{D}_n^{(1)}\bm{C}_{11,n}^{-1}\bm{x}_i^{(1)}, \bm{\eta}_i^{(0)}=\bm{D}_n^{(1)}\bm{C}_{11,n}^{-1}\bm{\eta}_i, \bm{b}_n= \dfrac{-\sqrt{n}\lambda_n}{2}\bm{D}_n^{(1)}\bm{C}_{11,n}^{-1}\bm{s}_n^{(1)}$$
with $\bm{\eta}_i=(\eta_{i,1},\dots,\eta_{i,p_0})^\prime$ and $\eta_{i,j}=\dfrac{-\lambda_n}{2}\big[\tilde{x}_{i,j}\big]sgn(\beta_{j,n})\tilde{P}^{\prime\prime}(|\beta_{j,n}|)$.\\

Now let us consider $\bm{T}_{2n}=n^{-1/2}\sum_{i=1}^{n}\big(\bm{\xi}_i^{(0)}+\bm{\eta}_i^{(0)}\big)\epsilon_i$. Due to (A.3)(ii) with $r=3/2$, (A.5)(i) and (A.7)(i), we have $\bm{\Upsilon}_n=\mathbf{Var}\big(\bm{T}_n\big)=n^{-1}\sum_{i=1}^{n}\big(\bm{\xi}_i^{(0)}+\bm{\eta}_i^{(0)}\big)\big(\bm{\xi}_i^{(0)}+\bm{\eta}_i^{(0)}\big)^\prime$ is a positive definite matrix and $n^{-1}\sum_{i=1}^{n}\big{\|}\bm{\xi}_i^{(0)}+\bm{\eta}_i^{(0)}\big{\|}^3\mathbf{E|\epsilon_1|^3}<\infty$. Therefore by Lemma \ref{lem:Normapprox}, we have 
\begin{align*}
\sup_{\bm{B}\in \mathcal{C}_q}\Big|\mathbf{P}\big(\bm{T}_{2n}\in \bm{B}\big)-\Phi\big(\bm{B}; \sigma^2\bm{\Upsilon}_n\big)\Big|=O\big(n^{-1/2}\big).
\end{align*}

%Similar to (\ref{eqn:sole}), the Lebesgue density of the two term EE of $\bm{T}_{2n}$ is

%\begin{align*}
%\bm{\Psi}_{2n}(\bm{x})=\Big[1-\dfrac{\mu_3}{\sqrt{n}}\sum_{|\bm{\alpha}|=3}\Big(n^{-1}\sum_{i=1}^{n}\big(\bm{\xi}_i^{(0)}+\bm{\eta}_i^{(0)}\big)^{\bm{\alpha}}\Big)D^{\bm{\alpha}}\Big]\phi(\bm{x}:\sigma^2\bm{\Upsilon}_n)
%\end{align*}
%where $\bm{\Upsilon}_n=n^{-1}\sum_{i=1}^{n}\big(\bm{\xi}_i^{(0)}+\bm{\eta}_i^{(0)}\big)\big(\bm{\xi}_i^{(0)}+\bm{\eta}_i^{(0)}\big)^\prime$. 
Again note that
$$||\bm{b}_n||\leq ||\bm{D}_n^{(1)}\bm{C}_{11,n}^{-1/2}||\cdot||\bm{C}_{11,n}^{-1/2}||\cdot||\bm{s}_n^{(1)}||\cdot \dfrac{\sqrt{n}\lambda_n}{2}\leq k\cdot p_0\cdot \tilde{P}_1\cdot \dfrac{\sqrt{n}\lambda_n}{2}=O(n^{-\delta}).$$
Suppose, $r=\min\{a\in \mathcal{N}:||\bm{b}_n||^{a+1}=O(n^{-1/2})\}$, $\mathcal{N}$ being the set of natural numbers.
%Therefore the Lebesgue density of the two term EE of $\bm{T}_{1n}$ is given by
%$\begin{align}\label{eqn:ole}
%\bm{\Psi}_{1n}(\bm{x})=\bigg[1+ \sum_{|\bm{\alpha}|=1}^{r}\bm{b}_n^{\bm{\alpha}}(-D)^{\bm{\alpha}}+\dfrac{\mu_3}{\sqrt{n}}\sum_{|\bm{\alpha}|=3}\Big(n^{-1}\sum_{i=1}^{n}\big(\bm{\xi}_i^{(0)}+\bm{\eta}_i^{(0)}\big)^{\bm{\alpha}}\Big)(-D)^{\bm{\alpha}}\bigg]\phi(\bm{x}:\sigma^2\bm{\Upsilon}_n)
%\end{align}
 Since $P(\bm{A}_{2n})=1-o(n^{-1/2})$, we have
 \begin{align}\label{eqn:orIIclt}
\sup_{\bm{B}\in \mathcal{C}_q}\Big|\mathbf{P}\big(\bm{T}_{n}\in \bm{B}\big)-\int_{\bm{B}}\tilde{\phi}_n(\bm{y})d\bm{y}\Big|=O\big(n^{-1/2}\big).
\end{align}
where $\tilde{\phi}_n(\bm{y})=\phi(\bm{y}:\sigma^2\bm{\Upsilon}_n)+ \sum_{|\bm{\alpha}|=1}^{r}(-\bm{b}_n)^{\bm{\alpha}}\chi_{\bm{\alpha}}(\bm{y}:\sigma^2\bm{\Upsilon}_n)\phi(\bm{y}:\sigma^2\bm{\Upsilon}_n)$.
%, with $\chi_{\bm{\alpha}}(\bm{y}:\sigma^2\bm{\Upsilon}_n)\phi(\bm{y}:\sigma^2\bm{\Upsilon}_n)$ is defined by the identity $$(-D)^{\bm{\alpha}}\phi(\bm{y}:\sigma^2\bm{\Upsilon}_n)=\chi_{\bm{\alpha}}(\bm{y}:\sigma^2\bm{\Upsilon}_n)\phi(\bm{y}:\sigma^2\bm{\Upsilon}_n).$$

%\begin{equation*}
%\sup\limits_{B \in  \mathcal{C}_q} \big|\mathbf{P}(\bm{T}_n\in B) - \int_B\bm{\Psi}_{1n}(\bm{x})d\bm{x}\big| = o(n^{-1/2}),
%\end{equation*} 
Hence we have
\begin{align*}
 \Delta_n =& \sup_{\bm{B}\in \mathcal{C}_q}|\mathbf{P}(\bm{T}_n\in \bm{B})-\Phi(\bm{B}; \sigma^2\bm{\Sigma}_n)|\\
 =&\sup_{\bm{B}}\bigg|\int_{\bm{B}}\big[\phi(\bm{x}:\sigma^2\bm{\Upsilon}_n)-\phi(\bm{x}:\sigma^2\bm{\Sigma}_n)\big]d\bm{x}\\
 &+\int_{\bm{B}}\sum_{|\bm{\alpha}|=1}^{r}(-\bm{b}_n)^{\bm{\alpha}}\chi_{\bm{\alpha}}(\bm{y}:\sigma^2\bm{\Upsilon}_n)\phi(\bm{y}:\sigma^2\bm{\Upsilon}_n)d{\bm{y}}\bigg|+O\big(n^{-1/2}\big)\\
 \leq & k.||\bm{\Upsilon}_n-\bm{\Sigma}_n|| + O\big(n^{-1/2}+||\bm{b}_n||\big)
\end{align*}
Now defining $\bm{\Lambda}_n=diag\big(sgn(\beta_{1,n})\tilde{P}^{\prime\prime}(|\beta_{1,n}|),\dots, \beta_{p_0,n})\tilde{P}^{\prime\prime}(|\beta_{p_0,n}|)\big)$, through the same line of the proof of Theorem 3.1 in Chatterjee and Lahiri(2013), it can be shown that
$$||\bm{\Upsilon}_n-\bm{\Sigma}_n|| = O\Big(\lambda_nn^a\tilde{P}_2+\lambda_n^2p_0^2\tilde{P}_2^2\Big).$$ Therefore Theorem \ref{thm:oracleII} follows.

\textbf{Proof of Theorem \ref{thm:oracleIII}}:
The Lasso estimator $\hat{\bm{\beta}}_n$ is defined as 
 \begin{equation*}
\bm{\hat{\beta}}_n = \operatorname*{arg\,min}_{\bm{t}}\Bigg[\sum_{i=1}^{n}(y_i - \bm{x}'_i \bm{t})^2
+\lambda_n\sum_{j=1}^{p}(|t_j|)\Bigg]
\end{equation*}

Now, writing $\bm{\hat{u}}_{3n}=\sqrt{n}\big(\bm{\hat{\beta}}_n-\bm{\beta}_n\big)$, we have
 
 \begin{align}
\bm{\hat{u}}_{3n}&= \operatorname*{arg\,min}_{\bm{v}}\Bigg[\bm{v}^{\prime}\bm{C}_n\bm{v} -2\bm{v}^{\prime}\bm{W}_n+\lambda_n\sum_{j=1}^{p}\Big(|\beta_{j,n}+\dfrac{v_{j}}{\sqrt{n}}|-|\beta_{j,n}|\Big)\Bigg]\nonumber\\\label{eqn:Z3n}
&= \operatorname*{arg\,min}_{\bm{v}} \bm{Z}_{3n}(\bm{v})\;\;\;\;\; \text{(say)}.
\end{align}

Note that $\bm{Z}_{3n}(\bm{v})$ is convex in $\bm{v}$. Hence, the KKT condition is necessary and sufficient. The KKT condition corresponding to (\ref{eqn:Z3n}) is given by
\begin{align}\label{eqn:KKTinproofs3}
2\bm{C}_n\bm{v}-2\bm{W}_n+\dfrac{\lambda_n}{\sqrt{n}}\bm{l}_n=\bm{0}
\end{align}
for some $l_{j,n}\in [-1,1]$ for all $j\in\{1,\ldots,p\}$, where $\bm{l}_n=(l_{1,n},\ldots, l_{p,n})'$. It is easy to show that under the conditions stated in theorem \ref{thm:oracleIII}, on the set $\bm{A}_{3n}$, $\Big(\big(\bm{\hat{u}}_{3n}^{(1)}\big)', \bm{0}'\Big)'$, where $\bm{\hat{u}}_{3n}^{(1)}=\bm{C}_{11,n}^{-1}\Big[\bm{W}_n^{(1)}-\dfrac{\lambda_n}{2\sqrt{n}}\bm{c}_{n}^{(1)}\Big]$ is the unique solution of (\ref{eqn:KKTinproofs3}) and hence $\bm{\hat{u}}_{3n}=\Big(\big(\bm{\hat{u}}_{3n}^{(1)}\big)', \bm{0}'\Big)'$, is the unique solution of the minimization problem (\ref{eqn:Z3n}), where $\bm{c}_{n}^{(1)}=(c_{1n},\ldots, c_{p_0n})$ and $c_{j,n}=sgn(\beta_{j,n})$ on the set $\bm{A}_{3n}$. Therefore on the set $\bm{A}_{3n}$ we have
\begin{align*}
\bm{T}_n&=\sqrt{n}\bm{D}_n(\hat{\bm{\beta}}_n - \bm{\beta})\\
&=\bm{D}_n^{(1)}\bm{C}_{11,n}^{-1}\Big[\bm{W}_n^{(1)}-\dfrac{\lambda_n}{2\sqrt{n}}\bm{c}_n^{(1)}\Big]\\
&= \bm{D}_n^{(1)}\bm{C}_{11,n}^{-1}\bm{W}_n^{(1)}-\dfrac{\lambda_n}{2\sqrt{n}}\bm{D}_n^{(1)}\bm{C}_{11,n}^{-1}\bm{c}_n^{(1)}\\
&=n^{-1/2}\sum_{i=1}^{n}\bm{\xi}_i^{(0)}\epsilon_i-\dfrac{\lambda_n}{2\sqrt{n}}\bm{D}_n^{(1)}\bm{C}_{11,n}^{-1}\bm{c}_n^{(1)}\\
&= \bm{T}_{3n}+\bm{b}^{\dagger}_n\;\;\;\; (say).
\end{align*}
Clearly the form of $\bm{T}_{3n}$ is same as $\bm{T}_n$ of Theorem 1 and hence under the conditions (A.3)(ii), (A.3)(iii), (A.5)(i) with $r=3/2$, Lemma \ref{lem:Normapprox} implies

%\begin{align}\label{eqn:sole}
%\bm{\Psi}_{3n}(\bm{x})=\Big[1-\dfrac{\mu_3}{\sqrt{n}}\sum_{|\bm{\alpha}|=3}\big(n^{-1}\sum_{i=1}^{n}(\bm{\xi}_i^{(0)})^{\bm{\alpha}}\big)D^{\bm{\alpha}}\Big]\phi(\bm{x}:\sigma^2\bm{\Sigma}_n),
%\end{align}
\begin{align}\label{eqn:sole}
\sup_{\bm{B}\in \mathcal{C}_q}|\mathbf{P}(\bm{T}_{3n}\in \bm{B})-\Phi(\bm{B}; \sigma^2\bm{\Sigma}_n)|=O(n^{-1/2}),
\end{align}
where $\bm{\Sigma}_n=\bm{D}_n^{(1)}\bm{C}_{11,n}^{-1}\bm{D}_n^{(1)}$. %Therefore we have
%\begin{equation}\label{eqn:EE3}
%\sup\limits_{B \in  \mathcal{C}_q} \big|\mathbf{P}(\bm{T}_{3n}\in B) - \int_B\bm{\Psi}_{3n}(\bm{x})d\bm{x}\big| = o(n^{-1/2}).
%\end{equation} 
Since we have assumed $\max\Big|\Big(\bm{\Sigma}_n^{-1/2}\Big)_{j\cdot}\Big[\bm{D}_n^{(1)}\bm{C}_{11,n}^{-1}\bm{c}_n^{(1)}\Big]\Big|>\kappa$ for some $\kappa >0$, without loss of generality we can consider $\Big|\Big(\bm{\Sigma}_n^{-1/2}\Big)_{1\cdot}\Big[\bm{D}_n^{(1)}\bm{C}_{11,n}^{-1}\bm{c}_n^{(1)}\Big]\Big|>\kappa$, where $\Big(\bm{\Sigma}_n^{-1/2}\Big)_{j\cdot}$ is the $j$th row of $\bm{\Sigma}_n^{-1/2}$ and $\bm{\Sigma}_n=\bm{D}_n^{(1)}\bm{C}_{11,n}^{-1}\bm{D}_n^{(1)}$. Then consider a set $\bm{B}_n$ in $\mathcal{R}^q$ as
\[
\bm{B}_n  = \left\{ \begin{array}{ll} \Big(-\infty,-\dfrac{\kappa\lambda_n}{4\sqrt{n}}\Big)\times \mathcal{R}\times\cdots\times \mathcal{R}, \;\;\; & \text{if $\Big(\bm{\Sigma}_n^{-1/2}\Big)_{1\cdot}\Big[\bm{D}_n^{(1)}\bm{C}_{11,n}^{-1}\bm{c}_n^{(1)}\Big]>\kappa$} \\\Big(\dfrac{\kappa\lambda_n}{4\sqrt{n}}, \infty\Big)\times \mathcal{R}\times\cdots\times \mathcal{R}, & \text{\text{if $\Big(\bm{\Sigma}_n^{-1/2}\Big)_{1\cdot}\Big[\bm{D}_n^{(1)}\bm{C}_{11,n}^{-1}\bm{c}_n^{(1)}\Big]<-\kappa$}} \end{array}  \right. \quad,
\]
%Since, $\bm{\Sigma}_n=\bm{D}_n^{(1)}\bm{C}_{11,n}^{-1}\bm{D}_n^{(1)}$ converges to a positive definite matrix, 
Since $\bm{B}_n$ is convex, $\tilde{\bm{B}}_n = \bm{\Sigma}_n^{1/2}\bm{B}_n=\{\bm{y}\in\mathcal{R}^q:\bm{y}=\bm{\Sigma}_n^{1/2}\bm{x}\; \text{for some}\; \bm{x}\in \bm{B}_n\}$ is also a convex set.  %the first row of $\bm{D}_n^{(1)}$ is $(sgn(\hat{\beta}_{1,n}),\dots,sgn(\hat{\beta}_{p_0,n}))^\prime=\hat{\bm{c}}_n^{(1)\prime}$. Note that on the set $\bm{A}_{3n}$, due to conditions (A.2)(ii), (A.2)(iii) and (A.7)(ii)$'$, we have
%\begin{align*}
%&b_{1,n}^{\dagger}=\dfrac{\lambda_n}{2\sqrt{n}}\hat{\bm{c}}_n^{(1)\prime}\bm{C}_{11,n}^{-1}\bm{c}_n^{(1)} \geq \dfrac{\lambda_n}{2\sqrt{n}}\cdot\bm{c}_n^{(1)\prime}\bm{c}_n^{(1)}\cdot p_0^{-1} = \dfrac{\lambda_n}{2\sqrt{n}} \geq \delta_1 p_0 \sqrt{\log n},\\
%& (\bm{\Sigma}_n)_{11} = \bm{c}_n^{(1)\prime}\bm{C}_{11,n}^{-1}\bm{c}_n^{(1)} \geq \bm{c}_n^{(1)\prime}\bm{c}_n^{(1)}\cdot p_0^{-1} =1\;\;\;\;\; \text{and}\\
%& (\bm{\Sigma}_n)_{11}=\bm{c}_n^{(1)\prime}\bm{C}_{11,n}^{-1}\bm{c}_n^{(1)} \leq  \bm{c}_n^{(1)\prime}\bm{c}_n^{(1)}\cdot||\bm{C}_{11,n}^{-1}|| \leq p_0^2,
%\end{align*}
%where $\bm{b}_n^{\dagger}=(b_{1,n}^{\dagger},\dots, b_{p_0,n}^{\dagger})^\prime$, $(\bm{\Sigma}_n)_{11}$ is the (1,1)th element of $\bm{\Sigma}_n$ and $\delta_1$ is a constant between $0$ and $1$. Hence 
%\begin{align}\label{eqn:41}
%&\dfrac{-(\log n)^{1/4}+b_{1,n}^{\dagger}}{\sqrt{(\bm{\Sigma}_n)_{11}}} \geq \dfrac{-(\log n)^{1/4}+\delta_1 p_0 \sqrt{\log n}}{p_0}\;\;\;\; \text{and}\\\nonumber
%&\dfrac{-(\log n)^{1/4}}{\sqrt{(\bm{\Sigma}_n)_{11}}} \leq -(\log n)^{1/4} 
%\end{align}
Therefore due to (\ref{eqn:sole}) and the fact that $\mathbf{P}(\bm{A}_{3n})=1-o(n^{-1/2})$, we have
\begin{align*}
\Delta_n =& \sup_{\bm{B}\in \mathcal{C}_q}\Big|\mathbf{P}(\bm{T}_n\in \bm{B})-\Phi(\bm{B}; \sigma^2\bm{\Sigma}_n)\Big| \\
\geq & \Big|\mathbf{P}(\bm{T}_n\in \tilde{\bm{B}}_n)-\Phi(\tilde{\bm{B}_n}; \sigma^2\bm{\Sigma}_n)\Big| \\
\geq & \Big|\mathbf{P}\Big(\{\bm{T}_n\in \tilde{\bm{B}}_n\}\cap \bm{A}_{3n}\Big)-\Phi(\tilde{\bm{B}}_n; \sigma^2\bm{\Sigma}_n)\Big|-\mathbf{P}(\bm{A}_{3n}^c)\\
= & \Big|\mathbf{P}\Big(\{\bm{T}_{3n}+\bm{b}^{\dagger}_n\in \tilde{\bm{B}}_n\}\cap \bm{A}_{3n}\Big)-\Phi(\tilde{\bm{B}}_n; \sigma^2\bm{\Sigma}_n)\Big|-\mathbf{P}(\bm{A}_{3n}^c)\\
\geq & \Big|\mathbf{P}\Big(\bm{T}_{3n}+\bm{b}^{\dagger}_n\in \tilde{\bm{B}}_n\Big)-\Phi(\tilde{\bm{B}}_n; \sigma^2\bm{\Sigma}_n)\Big|-2\mathbf{P}(\bm{A}_{3n}^c)\\
%\geq & \Big|\Phi(\tilde{\bm{B}}_n-\bm{b}^{\dagger}_n; \sigma^2\bm{\Sigma}_n)-\Phi(\tilde{\bm{B}}_n; \sigma^2\bm{\Sigma}_n)\Big|\\
%&-\dfrac{\mu_3}{\sqrt{n}}\sum_{|\bm{\alpha}|=3}n^{-1}\sum_{i=1}^{n}\int_{\mathcal{R}^{q}}\big|\big(\bm{\xi}_i^{(0)}\big)^{\bm{\alpha}}D^{\bm{\alpha}}\phi(\bm{x}:\sigma^2\Sigma_n)\big|d\bm{x} - o(n^{-1/2})\\
\geq & \Big|\Phi(\tilde{\bm{B}}_n-\bm{b}^{\dagger}_n; \sigma^2\bm{\Sigma}_n)-\Phi(\tilde{\bm{B}}_n; \sigma^2\bm{\Sigma}_n)\Big| - O(n^{-1/2})\\
=& \Big|\Phi(\bm{B}_n-\bm{\Sigma}_n^{-1/2}\bm{b}^{\dagger}_n; \sigma^2\bm{I}_q)-\Phi(\bm{B}_n; \sigma^2\bm{I}_q)\Big| - O(n^{-1/2})
%=& \bigg[\int_{-(\log n)^{1/4}}^{-(\log n)^{1/4}+b_{1,n}^{\dagger}}(2\pi)^{-1/2}\big[(\bm{\Sigma}_n)_{11}\big]^{-1/2}\exp\Big(-y^2/\Big[2\sqrt{(\bm{\Sigma}_n)_{11}}\Big]\Big)dy\bigg] - O(n^{-1/2})\\
%\geq& \bigg[\int_{g_n}^{h_n}(2\pi)^{-1/2}\exp\big(-y^2/2\big)dy\bigg] - O(n^{-1/2})
\end{align*}
Again note that $\Big|\Big(\bm{\Sigma}_n^{-1/2}\Big)_{1\cdot}\bm{b}_n^{\dagger}\Big|\geq \dfrac{\kappa\lambda_n}{2\sqrt{n}} \geq \delta_1\cdot \kappa \cdot \sqrt{\log n}$, due to condition (A.7)(ii)$'$. Therefore, 
\begin{align*}
&\Big|\Phi(\bm{B}_n-\bm{\Sigma}_n^{-1/2}\bm{b}^{\dagger}_n; \sigma^2\bm{I}_q)-\Phi(\bm{B}_n; \sigma^2\bm{I}_q)\Big|\\
\geq & \bigg[\int_{-4^{-1}\delta_1\kappa\sqrt{\log n}}^{4^{-1}\delta_1\kappa\sqrt{\log n}}(2\pi)^{-1/2}\exp\big(-y^2/2\big)dy\bigg]\\
\rightarrow& 1\;\;\;\;\; \text{as}\; n \rightarrow \infty. 
\end{align*}
Therefore Theorem \ref{thm:oracleIII} follows.

\textbf{Proof of Theorem \ref{thm:boot}}:\\
A. Proof of part 1: Note that on the set $\bm{A}_n^*\cap \bm{B}_n^*$,
\begin{align*}
\bm{T}_n^*=\sqrt{n}\bm{D}_n(\hat{\bm{\beta}}_n^* - \hat{\bm{\beta}}_n) = n^{-1/2}\sum_{i=1}^{n}\bm{D}_n^{(1)}\bm{C}_{11,n}^{-1}\bm{x}_i^{(1)}\epsilon_i^* = n^{-1/2}\sum_{i=1}^{n}\bm{\xi}_i^{(0)}\epsilon_i^*.
\end{align*}
Simularly, on the set $\bm{A}_n^{**}\cap \bm{B}_n^{**}$,
\begin{align*}
\bm{T}_n^{**}=\sqrt{n}\bm{D}_n(\hat{\bm{\beta}}_n^{**} - \hat{\bm{\beta}}_n) &= n^{-1/2}\sum_{i=1}^{n}\bm{D}_n^{(1)}\bm{C}_{11,n}^{-1}\bm{x}_i^{(1)}\hat{\epsilon}_i(G_i^*-\mu_{G^*})\mu_{G^*}^{-1}\\
&= n^{-1/2}\sum_{i=1}^{n}\bm{\xi}_i^{(0)}\hat{\epsilon}_i(G_i^*-\mu_{G^*})\mu_{G^*}^{-1}
\end{align*}

Then by the conditions (A.3)(ii), (A.3)(iii) and (A.5)(i) with $r=3$ and the fact that $\mathbf{P}_*(\bm{A}_n^*\cap \bm{B}_n^*)=1-o_p(n^{-1/2})$ or $\mathbf{P}_*(\bm{A}_n^{**}\cap \bm{B}_n^{**})=1-o_p(n^{-1/2})$, Lemma \ref{lem:Normapprox} yields
\begin{align*}
&\sup_{\bm{B}\in \mathcal{B}_q}\Big|\mathbf{P}_*\big(\bm{T}_n^*\in \bm{B}\big)-\Phi\big(\bm{B}; \hat{\sigma}_n^2\bm{\Sigma}_n\big)\Big|=O\big(n^{-1/2}\big)\;\;\; \text{and}\\
&\sup_{\bm{B}\in \mathcal{B}_q}\Big|\mathbf{P}_*\big(\bm{T}_n^{**}\in \bm{B}\big)-\Phi\big(\bm{B}; \tilde{\bm{\Sigma}}_n\big)\Big|=O\big(n^{-1/2}\big)
\end{align*}
where $\bm{\Sigma}_n=n^{-1}\sum_{i=1}^{n}\bm{\xi}_i^{(0)}\bm{\xi}_i^{(0)\prime}=\bm{D}_n^{(1)}\bm{C}_{11,n}^{-1}\bm{D}_n^{(1)}$, $\hat{\sigma}_n^2=n^{-1}\sum_{i=1}^{n}(\hat{\epsilon}_i-\bar{\epsilon}_n)^2$ and $\tilde{\bm{\Sigma}}_n=n^{-1}\sum_{i=1}^{n}\bm{\xi}_i^{(0)}\bm{\xi}_i^{(0)\prime}\hat{\epsilon}_i^2$. Therefore noting that $|\hat{\sigma}_n^2-\sigma^2|=O_p(n^{-1/2})$ and due to Lemma \ref{lem:Sigma} we have
\begin{align*}
\Delta_n^* & = \sup_{\bm{B}\in \mathcal{C}_q}|\mathbf{P}_*(\bm{T}_n^*\in \bm{B})-\mathbf{P}(\bm{T}_n\in \bm{B})|\\
&\leq \sup_{\bm{B}\in \mathcal{C}_q}\Big|\Phi\big(\bm{B}; \hat{\sigma}_n^2\bm{\Sigma}_n\big)-\Phi\big(\bm{B}; \sigma^2\bm{\Sigma}_n\big)\Big|+O_p(n^{-1/2})\\
& \leq k\cdot|\hat{\sigma}_n^2-\sigma^2|+O_p(n^{-1/2})\\
& = O_p(n^{-1/2})\;\;\;\; \text{and}
\end{align*}
\begin{align*}
\Delta_n^{**} & = \sup_{\bm{B}\in \mathcal{C}_q}|\mathbf{P}_*(\bm{T}_n^{**}\in \bm{B})-\mathbf{P}(\bm{T}_n\in \bm{B})|\\
&\leq \sup_{\bm{B}\in \mathcal{C}_q}\Big|\Phi\big(\bm{B}; \tilde{\bm{\Sigma}}_n\big)-\Phi\big(\bm{B}; \sigma^2\bm{\Sigma}_n\big)\Big|+O_p(n^{-1/2})\\
& \leq k\cdot\|\tilde{\bm{\Sigma}}_n-\sigma^2\bm{\Sigma}_n\|+O_p(n^{-1/2})\\
& = O_p(n^{-1/2})
\end{align*}
B. Proof of part 2: Through the same line of arguments as in the proof of Theorem \ref{thm:oracleII}, it can be shown that on a set $\bm{A}_{2n}^*$,
\begin{align*}
\bm{T}_n^*&=\sqrt{n}\bm{D}_n(\hat{\beta}_n^*-\hat{\beta}_n)\\
&= n^{-1/2}\sum_{i=1}^{n}\big(\hat{\bm{\xi}}_i^{(0)}+\hat{\bm{\eta}}_i^{(0)}\big)\epsilon_i^*+ \hat{\bm{b}}_n + \bm{Q}_{1n}^*\\
&= \bm{T}_{2n}^*+ \hat{\bm{b}}_n + \bm{Q}_{1n}^*\;\;\; \text{and}\\
\bm{T}_n^{**}&=\sqrt{n}\bm{D}_n(\hat{\beta}_n^{**}-\hat{\beta}_n)\\
&= n^{-1/2}\sum_{i=1}^{n}\big(\hat{\bm{\xi}}_i^{(0)}+\hat{\bm{\eta}}_i^{(0)}\big)\hat{\epsilon}_i(G_i^*-\mu_{G^*})\mu_{G^*}^{-1}+ \hat{\bm{b}}_n + \bm{Q}_{1n}^{**}\\
&= \bm{T}_{2n}^{**}+ \hat{\bm{b}}_n + \bm{Q}_{1n}^{**}
\end{align*}
where 
$$\hat{\bm{\xi}}_i^{(0)}=\hat{\bm{D}}_n^{(1)}\hat{\bm{C}}_{11,n}^{-1}\hat{\bm{x}}_i^{(1)}, \hat{\bm{\eta}}_i^{(0)}=\bm{D}_n^{(1)}\bm{C}_{11,n}^{-1}\hat{\bm{\eta}}_i, \hat{\bm{b}}_n= \dfrac{-\sqrt{n}\lambda_n}{2}\hat{\bm{D}}_n^{(1)}\hat{\bm{C}}_{11,n}^{-1}\hat{\bm{s}}_n^{(1)}$$
with $\hat{\bm{\eta}}_i=(\hat{\eta}_{i,1},\dots,\hat{\eta}_{i,p_0})^\prime$, $\hat{\eta}_{i,j}=\dfrac{-\lambda_n}{2}\big[\tilde{x}_{i,j}\big]sgn(\hat{\beta}_{j,n})\tilde{P}^{\prime\prime}(|\hat{\beta}_{j,n}|)$, $\hat{\bm{s}}_n^{(1)}=(\hat{s}_{1,n},\ldots, \hat{s}_{p_0,n})$ and $\hat{s}_{j,n}=sgn(\hat{\beta}_{j,n})\tilde{P}^\prime(|\hat{\beta}_{j,n}|)$. $\hat{\bm{D}}_n^{(1)}$, $\hat{\bm{C}}_{11,n}$ and $\hat{\bm{x}}_i^{(1)}$ are as defined in Section \ref{sec:mainresults}. Note that for large enough $n$, due to Lemma \ref{lem:betahat}, condition (A.7)(i) and the continuity of %\rd{$\max_j\tilde{P}^{\prime\prime\prime}(t_j)$ at $\bm{t}=(t_1,\dots,t_{p_0})^\prime=\Big(|\beta_{1,n}|,\dots,|\beta_{p_0,n}|\Big)^\prime$}
$\max_j\tilde{P}^{\prime\prime\prime}(t_j)$ at $\bm{t}=(t_1,\dots,t_{p_0})^\prime=\Big(|\beta_{1,n}|,\dots,|\beta_{p_0,n}|\Big)^\prime$, we have
$$||\hat{\bm{b}}_n||\leq ||\bm{D}_n^{(1)}\bm{C}_{11,n}^{-1/2}||\cdot||\bm{C}_{11,n}^{-1/2}||\cdot \dfrac{\sqrt{n}\lambda_n}{2}\cdot||\hat{\bm{s}}_n^{(1)}||\leq k\cdot p_0\cdot \dfrac{\sqrt{n}\lambda_n}{2}\cdot O_p(\tilde{P}_1+\tilde{P}_2/\sqrt{n}+\tilde{P}_3/n)=O_p(n^{-\delta}).$$
Suppose, $\hat{r}=\min\{a\in \mathcal{N}:||\hat{\bm{b}}_n||^{a+1}=O_p(n^{-1/2})\}$, $\mathcal{N}$ being the set of natural numbers. Therefore using Lemma \ref{lem:Normapprox}, it is easy to see
 \begin{align}\label{eqn:orIIbootclt}
&\sup_{\bm{B}\in \mathcal{C}_q}\Big|\mathbf{P}_*\big(\bm{T}_{n}^*\in \bm{B}\big)-\int_{\bm{B}}\tilde{\phi}_n^*(\bm{y})d\bm{y}\Big|=O_p\big(n^{-1/2}\big)\;\;\; \text{and}\\\nonumber
&\sup_{\bm{B}\in \mathcal{C}_q}\Big|\mathbf{P}_*\big(\bm{T}_{n}^{**}\in \bm{B}\big)-\int_{\bm{B}}\tilde{\phi}_n^{**}(\bm{y})d\bm{y}\Big|=O_p\big(n^{-1/2}\big)
\end{align}
where $\tilde{\phi}_n^*(\bm{y})=\phi(\bm{y}:\sigma_n^2\hat{\bm{\Upsilon}}_n)+ \sum_{|\bm{\alpha}|=1}^{\hat{r}}(-\hat{\bm{b}}_n)^{\bm{\alpha}}\chi_{\bm{\alpha}}(\bm{y}:\sigma_n^2\hat{\bm{\Upsilon}}_n)\phi(\bm{y}:\sigma_n^2\hat{\bm{\Upsilon}}_n)$ and $\tilde{\phi}_n^{**}(\bm{y})=\phi(\bm{y}:\tilde{\bm{\Upsilon}}_n)+ \sum_{|\bm{\alpha}|=1}^{\hat{r}}(-\hat{\bm{b}}_n)^{\bm{\alpha}}\chi_{\bm{\alpha}}(\bm{y}:\tilde{\bm{\Upsilon}}_n)\phi(\bm{y}:\tilde{\bm{\Upsilon}}_n)$ where $\hat{\bm{\Upsilon}}_n=n^{-1}\sum_{i=1}^{n}(\hat{\bm{\xi}}_i^{(0)}+\hat{\bm{\eta}}_i^{(0)})(\hat{\bm{\xi}}_i^{(0)}+\hat{\bm{\eta}}_i^{(0)})^\prime$ and $\tilde{\bm{\Upsilon}}_n=n^{-1}\sum_{i=1}^{n}(\hat{\bm{\xi}}_i^{(0)}+\hat{\bm{\eta}}_i^{(0)})(\hat{\bm{\xi}}_i^{(0)}+\hat{\bm{\eta}}_i^{(0)})^\prime\hat{\epsilon}_i^2$.
%The quantities $\chi_{\bm{\alpha}}(\bm{y}:\sigma_n^2\bm{\Upsilon}_n)\phi(\bm{y}:\sigma_n^2\hat{\bm{\Upsilon}}_n)$ and $\chi_{\bm{\alpha}}(\bm{y}:\tilde{\bm{\Upsilon}}_n)\phi(\bm{y}:\tilde{\bm{\Upsilon}}_n)$ are defined by the identities $$(-D)^{\bm{\alpha}}\phi(\bm{y}:\sigma_n^2\hat{\bm{\Upsilon}}_n)=\chi_{\bm{\alpha}}(\bm{y}:\sigma_n^2\hat{\bm{\Upsilon}}_n)\phi(\bm{y}:\sigma_n^2\hat{\bm{\Upsilon}}_n)\;\;\; \text{and}$$ $$(-D)^{\bm{\alpha}}\phi(\bm{y}:\tilde{\bm{\Upsilon}}_n)=\chi_{\bm{\alpha}}(\bm{y}:\tilde{\bm{\Upsilon}}_n)\phi(\bm{y}:\tilde{\bm{\Upsilon}}_n).$$
Now note that $|\sigma_n^2-\sigma^2|=O_p(n^{-1/2})$ and by Lemma \ref{lem:betahat} for large enough $n$ we have
\begin{align*}
\|\hat{\bm{b}}_n-\bm{b}_n\|&\leq ||\bm{D}_n^{(1)}\bm{C}_{11,n}^{-1/2}||\cdot||\bm{C}_{11,n}^{-1/2}||\cdot||\hat{\bm{s}}_n^{(1)}-\bm{s}_n^{(1)}||\cdot \dfrac{\sqrt{n}\lambda_n}{2}\\
&\leq \dfrac{\sqrt{n}\lambda_n}{2}\cdot p_0 \cdot O_p(\tilde{P}_2/\sqrt{n}+\tilde{P}_3/n)\\
&=o_p(n^{-1/2}).
\end{align*}
Again by Lemma \ref{lem:Sigma}, $\|\hat{\bm{\Upsilon}}_n-\bm{\Upsilon}_n\|=o_p(n^{-1/2})$ and $\|\tilde{\bm{\Upsilon}}_n-\sigma^2\bm{\Upsilon}_n\|=O_p(n^{-1/2})$. Therefore by comparing (\ref{eqn:orIIclt}) and (\ref{eqn:orIIbootclt}), part 2 of Theorem  \ref{thm:boot} follows.\\

C. Proof of part 3: Similar to original Lasso case, it can be shown that there exists a bootstrap set $\bm{A}_{3n}^*$ with $\mathbf{P}_*(\bm{A}_{3n}^*)=1-o_p(n^{-1/2})$ such that on the set $\bm{A}_{3n}^*$,
\begin{align*}
\bm{T}_n^*&=\sqrt{n}\bm{D}_n(\hat{\bm{\beta}}_n^* - \hat{\bm{\beta}}_n)\\
%&=\hat{\bm{D}}_n^{(1)}\hat{\bm{C}}_{11,n}^{-1}\Big[\bm{W}_n^{(1)}-\dfrac{\lambda_n}{2\sqrt{n}}\bm{c}_n^{(1)}\Big]\\
&= n^{-1/2}\sum_{i=1}^{n}\hat{\bm{\xi}}_i^{(0)}\epsilon_i^*-\dfrac{\lambda_n}{2\sqrt{n}}\hat{\bm{D}}_n^{(1)}\hat{\bm{C}}_{11,n}^{-1}\hat{\bm{c}}_n^{(1)}\\
&= \bm{T}_{3n}^*+\hat{\bm{b}}^{\dagger}_n\;\;\;\; \text{and}\\
\bm{T}_n^{**}&=\sqrt{n}\hat{\bm{D}}_n(\hat{\bm{\beta}}_n^{**} - \hat{\bm{\beta}}_n)\\
%&=\hat{\bm{D}}_n^{(1)}\hat{\bm{C}}_{11,n}^{-1}\Big[\bm{W}_n^{(1)}-\dfrac{\lambda_n}{2\sqrt{n}}\bm{c}_n^{(1)}\Big]\\
&= n^{-1/2}\sum_{i=1}^{n}\hat{\bm{\xi}}_i^{(0)}\hat{\epsilon}_i(G_i^*-\mu_{G^*})\mu_{G^*}^{-1}-\dfrac{\lambda_n}{2\sqrt{n}}\hat{\bm{D}}_n^{(1)}\hat{\bm{C}}_{11,n}^{-1}\hat{\bm{c}}_n^{(1)}\\
&= \bm{T}_{3n}^{**}+\hat{\bm{b}}^{\dagger}_n
\end{align*}
where $\hat{\bm{\xi}}_i^{(0)}$, $i\in\{1,\dots,n\}$, are defined earlier and $\hat{\bm{b}}_n^{\dagger}=-\dfrac{\lambda_n}{2\sqrt{n}}\hat{\bm{D}}_n^{(1)}\hat{\bm{C}}_{11,n}^{-1}\hat{\bm{c}}_n^{(1)}$ with $\hat{\bm{c}}_{n}^{(1)}=(\hat{c}_{1n},\ldots, \hat{c}_{p_0n})$ and $\hat{c}_{j,n}=sgn(\hat{\beta}_{j,n})$. Note that on the set $\bm{A}_{3n}$ with $\mathbf{P}(\bm{A}_{3n})=1-o(n^{-1/2})$, $\hat{c}_{j,n}=c_{j,n}$ for all $j\in\{1,\dots,p_0\}$ and hence $\hat{\bm{b}}_n^{\dagger}=\bm{b}_n^{\dagger}$. Therefore for any $k>0$ and sufficiently large $n$,
\begin{align*}
&\mathbf{P}\bigg[\Delta_n^* = \sup_{\bm{B}\in \mathcal{C}_q}|\mathbf{P}_*(\bm{T}_n^{*}\in \bm{B})-\mathbf{P}(\bm{T}_n\in \bm{B})|>k\cdot n^{-1/2}\bigg]\\
\leq & \mathbf{P}\bigg[\sup_{\bm{B}\in \mathcal{C}_q}\Big|\mathbf{P}_*\Big(\{\bm{T}_n^{*}\in \bm{B}\}\cap \bm{A}_{3n}^*\Big)-\mathbf{P}\Big(\{\bm{T}_n\in \bm{B}\}\cap \bm{A}_{3n}\Big)\Big|+\mathbf{P}_{*}(\bm{A}_{3n}^{*c})+\mathbf{P}(\bm{A}_{3n}^c)>k\cdot n^{-1/2}\bigg]\\
\leq & \mathbf{P}\bigg[\sup_{\bm{B}\in \mathcal{C}_q}\Big|\mathbf{P}_*\Big(\{\bm{T}_n^{*}\in \bm{B}\}\Big)-\mathbf{P}\Big(\{\bm{T}_n\in \bm{B}\}\Big)\Big|+2\mathbf{P}_{*}(\bm{A}_{3n}^{*c})+2\mathbf{P}(\bm{A}_{3n}^c)>k\cdot n^{-1/2}\bigg]\\
\leq & \mathbf{P}\bigg[\sup_{\bm{B}\in \mathcal{C}_q}\Big|\mathbf{P}_*\Big(\bm{T}_{3n}^{*}+\hat{\bm{b}}_n^{\dagger}\in \bm{B}\Big)-\mathbf{P}\Big(\bm{T}_{3n}+ \bm{b}_n^{\dagger}\in\bm{B}\Big)\Big|+2\mathbf{P}_{*}(\bm{A}_{3n}^{*c})>(k/2)\cdot n^{-1/2}\bigg]\\
\leq & \mathbf{P}\bigg[\Big\{\sup_{\bm{B}\in \mathcal{C}_q}\Big|\mathbf{P}_*\Big(\bm{T}_{3n}^{*}+\hat{\bm{b}}_n^{\dagger}\in \bm{B}\Big)-\mathbf{P}\Big(\bm{T}_{3n}+ \bm{b}_n^{\dagger}\in\bm{B}\Big)\Big|>(k/4)\cdot n^{-1/2}\Big\}\cap \bm{A}_{3n}\bigg]+\mathbf{P}(\bm{A}_{3n}^c)\\
\leq & \mathbf{P}\bigg[\Big\{\sup_{\bm{B}\in \mathcal{C}_q}\Big|\mathbf{P}_*\Big(\bm{T}_{3n}^{*}+\bm{b}_n^{\dagger}\in \bm{B}\Big)-\mathbf{P}\Big(\bm{T}_{3n}+ \bm{b}_n^{\dagger}\in\bm{B}\Big)\Big|>(k/4)\cdot n^{-1/2}\Big\}\bigg]+o(n^{-1/2})\\
= & \mathbf{P}\bigg[\Big\{\sup_{\bm{C}\in \mathcal{C}_q}\Big|\mathbf{P}_*\Big(\bm{T}_{3n}^{*}\in \bm{C}\Big)-\mathbf{P}\Big(\bm{T}_{3n}\in\bm{C}\Big)\Big|>(k/4)\cdot n^{-1/2}\Big\}\bigg]+o(n^{-1/2})
\end{align*}
Fourth inequality follows from the fact that $\mathbf{P}_{*}\Big(\bm{A}_{3n}^{*c}\Big)=o(n^{-1/2})$ on the set $\bm{A}_{3n}$. Therefore, to prove the result corresponding to residual bootstrap in part 3. of Theorem \ref{thm:boot}, it is enough to show
\begin{align*}
\sup_{\bm{C}\in \mathcal{C}_q}\Big|\mathbf{P}_*\Big(\bm{T}_{3n}^{*}\in \bm{C}\Big)-\mathbf{P}\Big(\bm{T}_{3n}\in\bm{C}\Big)\Big| = O_p(n^{-1/2}).
\end{align*}
Similar arguments will also ensure that the following is enough to establish the result corresponding to perturbation bootstrap:
\begin{align*}
\sup_{\bm{C}\in \mathcal{C}_q}\Big|\mathbf{P}_*\Big(\bm{T}_{3n}^{**}\in \bm{C}\Big)-\mathbf{P}\Big(\bm{T}_{3n}\in\bm{C}\Big)\Big| = O_p(n^{-1/2}).
\end{align*}
Now similar to (\ref{eqn:sole}), by Lemma \ref{lem:Normapprox} we have
\begin{align}\label{eqn:soleboot}
&\sup_{\bm{B}\in \mathcal{B}_q}\Big|\mathbf{P}_*\big(\bm{T}_n^*\in \bm{B}\big)-\Phi\big(\bm{B}; \hat{\sigma}_n^2\bm{\Sigma}_n\big)\Big|=O\big(n^{-1/2}\big)\;\;\; \text{and}\\\nonumber
&\sup_{\bm{B}\in \mathcal{B}_q}\Big|\mathbf{P}_*\big(\bm{T}_n^{**}\in \bm{B}\big)-\Phi\big(\bm{B}; \tilde{\bm{\Sigma}}_n\big)\Big|=O\big(n^{-1/2}\big)
\end{align}
where $\bm{\Sigma}_n=n^{-1}\sum_{i=1}^{n}\bm{\xi}_i^{(0)}\bm{\xi}_i^{(0)\prime}=\bm{D}_n^{(1)}\bm{C}_{11,n}^{-1}\bm{D}_n^{(1)}$, $\hat{\sigma}_n^2=n^{-1}\sum_{i=1}^{n}(\hat{\epsilon}_i-\bar{\epsilon}_n)^2$ and $\tilde{\bm{\Sigma}}_n=n^{-1}\sum_{i=1}^{n}\bm{\xi}_i^{(0)}\bm{\xi}_i^{(0)\prime}\hat{\epsilon}_i^2$. Now comparing (\ref{eqn:sole}) and (\ref{eqn:soleboot}) and noting that $|\hat{\sigma}_n^2-\sigma^2|=O_p(n^{-1/2})$ and $\|\tilde{\bm{\Sigma}}_n-\sigma^2\bm{\Sigma}_n\|=O_p(n^{-1/2})$, part 3. follows.

\textbf{Proof of Theorem \ref{thm:bootI}:} Note that on the set $\bm{A}_n\cap \bm{B}_n$,
\begin{align}
\bm{R}_n=&\bm{T}_n/\hat{\sigma}_n\nonumber\\
=&\Big[\sigma^{-1}-2^{-1}\sigma^{-3}(\hat{\sigma}_n^2-\sigma^2)\Big]\sqrt{n}\bm{D}_n(\hat{\bm{\beta}}_n - \bm{\beta}_n) +\bm{Q}_{2n}\nonumber\\
=& \Big[\sigma^{-1}-2^{-1}\sigma^{-3}(\hat{\sigma}_n^2-\sigma^2)\Big]n^{-1/2}\sum_{i=1}^{n}\bm{D}_n^{(1)}\bm{C}_{11,n}^{-1}\bm{x}_i^{(1)}\epsilon_i +\bm{Q}_{2n}\nonumber\\
=& n^{-1/2}\sum_{i=1}^{n}\bm{\xi}_i^{(0)}\epsilon_i/\sigma -2^{-1}\sigma^{-3}\Big[n^{-1}\sum_{i=1}^{n}(\epsilon_i^2-\sigma^2)\Big]n^{-1/2}\sum_{i=1}^{n}\bm{\xi}_i^{(0)}\epsilon_i +\bm{Q}_{3n}\nonumber\\\label{eqn:stuor1}
=& \tilde{\bm{R}}_n +\bm{Q}_{3n}\;\;\;\; \text{(say)},
\end{align}
where due to conditions (A.2)(ii), (A.3)(ii), (A.5)(i) and (A.7)(iii), $\mathbf{P}\Big(\|\bm{Q}_{2n}\|+\|\bm{Q}_{3n}\|=o(n^{-1/2})\Big)=1-o(n^{-1/2})$. Therefore two term EE of $\bm{R}_n$ and $\tilde{\bm{R}}_n$ agree up to $o(n^{-1/2})$. Now we are going to use the transformation technique of Bhattacharya and Ghosh (1978) to find two term EE of $\tilde{\bm{R}}_n$. However before employing that technique we need to find two term EE of $\bm{R}_n^{\dagger}=\Big[n^{-1/2}\sum_{i=1}^{n}\big(\bm{\xi}_i^{(0)\prime}\epsilon_i,(\epsilon_i^2-\sigma^2)\big)^\prime\Big]$. First three cumulants of $\bm{t}^\prime\bm{R}_{n}^{\dagger}$ are given by
\begin{align*}
\kappa_1\big(\bm{t}^\prime\bm{R}_{n}^{\dagger}\big)&=\mathbf{E}\big(\bm{t}^\prime\bm{R}_{n}^{\dagger}\big)=\bm{0}\\
\kappa_2\big(\bm{t}^\prime\bm{R}_{n}^{\dagger}\big)&=\mathbf{Var}\big(\bm{t}^\prime\bm{R}_{n}^{\dagger}\big)=\bm{t}^\prime\begin{bmatrix}
\Big(n^{-1}\sum_{i=1}^{n}\bm{\xi}_i^{(0)}\bm{\xi}_i^{(0)\prime}\Big)\sigma^2 \;\;\;\Big(n^{-1}\sum_{i=1}^{n}\bm{\xi}_i^{(0)}\Big)\mu_{3}\\
\Big(n^{-1}\sum_{i=1}^{n}\bm{\xi}_i^{(0)\prime}\Big)\mu_3\;\;\;\;\;\;\;\;\;\;\;\;\;\; (\mu_4-\sigma^4)
\end{bmatrix}\bm{t}=\bm{t}'\bm{S}_n\bm{t}\\
\kappa_3\big(\bm{t}^\prime\bm{R}_{n}^{\dagger}\big)&=\mathbf{E}\big(\bm{t}^\prime\bm{R}_{n}^{\dagger}\big)^3-3\mathbf{E}\big(\bm{t}^\prime\bm{R}_{n}^{\dagger}\big)^2. \mathbf{E}\big(\bm{t}^\prime\bm{R}_{n}^{\dagger}\big)+2\Big(\mathbf{E}\big(\bm{t}^\prime\bm{R}_{n}^{\dagger}\big)\Big)^3\\
&= n^{-1/2}\mu_3\bigg[\sum_{|\bm{\alpha}|\leq 3}\bm{t}^{\bm{\alpha}}\bar{\bm{\xi}}_n^{(0)}(\bm{\alpha})\mathbf{E}\Big[\epsilon_1^{|\bm{\alpha}|}(\epsilon_1^2-\sigma^2)^{3-|\bm{\alpha}|}\Big]\bigg]
\end{align*}
where $\bar{\bm{\xi}}_n^{(0)}(\bm{\alpha})=n^{-1}\sum_{i=1}^{n}(\bm{\xi}_i^{(0)})^{\bm{\alpha}}$. Hence by Theorem 20.6 of Bhattacharya and Rao (1986) we have
\begin{align}\label{eqn:orEE0}
\sup_{\bm{B}\in \mathcal{C}_{q+1}}\Big|\mathbf{P}(\bm{R}_n^{\dagger}\in \bm{B})-\int_{\bm{B}} \bm{\psi}_{n}^{(0)}(\bm{x})d\bm{x}\Big|=o(n^{-1/2})
\end{align}
where
\begin{align*}
\bm{\psi}_{n}^{(0)}(\bm{x})=&\phi(\bm{x}:\bm{\Sigma}_n)\Bigg[1+\dfrac{\mu_3}{\sqrt{n}}\sum_{|\bm{\alpha}|\leq 3}\bar{\bm{\xi}}_n^{(0)}(\bm{\alpha})\mathbf{E}\Big[\epsilon_1^{|\bm{\alpha}|}(\epsilon_1^2-\sigma^2)^{3-|\bm{\alpha}|}\Big]\chi_{\bm{\alpha}}(\bm{x}:\bm{S}_n)\Bigg],
\end{align*}
provided there exists $\delta_4 \in (0,1)$, independent of $n$, such that for all $\upsilon \leq \delta_4$,
\begin{align}\label{eqn:orIEEres1}
n^{-1}\sum_{i=1}^{n}\mathbf{E}\Big\|\big(\bm{\xi}_i^{(0)\prime}\epsilon_i,(\epsilon_i^2-\sigma^2)\big)^\prime\Big\|^{3}\mathbf{1}\Big(\Big\|\big(\bm{\xi}_i^{(0)\prime}\epsilon_i,(\epsilon_i^2-\sigma^2)\big)^\prime\Big\|> \upsilon \sqrt{n}\Big)=o(1)
\end{align}
and
\begin{align}\label{eqn:orIEEres2}
\max_{|\bm{\alpha}|\leq q+1}\int_{\|\bm{t}\|\geq \upsilon\sqrt{n}}\Big|D^{\bm{\alpha}}\mathbf{E}\exp(i\bm{t}^\prime\bm{R}_{1n}^{\dagger})\Big|=o\Big(n^{-1/2}\Big)
\end{align}
where $\bm{R}_{1n}^{\dagger}=n^{-1/2}\sum_{i=1}^{n}\big(\bm{Z}_{i}-\mathbf{E}\bm{Z}_i\big)$ with $$\bm{Z}_{i}=\big(\bm{\xi}_i^{(0)\prime}\epsilon_i,(\epsilon_i^2-\sigma^2)\big)^\prime\mathbf{1}\Big(\Big\|\big(\bm{\xi}_i^{(0)\prime}\epsilon_i,(\epsilon_i^2-\sigma^2)\big)^\prime\Big\|\leq \upsilon \sqrt{n}\Big).$$

First consider (\ref{eqn:orIEEres1}). Note that due to the condition (A.3)(ii) with $r=4$,
$$\max\Big\{\|\bm{\xi}_i^{(0)}\|: i\in\{1,\dots,n\}\Big\} = O\big(n^{1/8}\big).$$
Therefore, due to (A.5)(i) with $r=4$, we have for any $\upsilon >0$,
\begin{align*}
&n^{-1}\sum_{i=1}^{n}\mathbf{E}\Big\|\big(\bm{\xi}_i^{(0)\prime}\epsilon_i,(\epsilon_i^2-\sigma^2)\big)^\prime\Big\|^{3}\mathbf{1}\Big(\Big\|\big(\bm{\xi}_i^{(0)\prime}\epsilon_i,(\epsilon_i^2-\sigma^2)\big)^\prime\Big\|> \upsilon \sqrt{n}\Big)\\
\leq& n^{-1}\sum_{i=1}^{n}\mathbf{E}\Big(\big\|\bm{\xi}_i^{(0)}\big\|^2\epsilon_1^2+(\epsilon_1^2-\sigma^2)^2\Big)^{3/2}\mathbf{1}\Big(\big\|\bm{\xi}_i^{(0)}\big\|^2\epsilon_1^2+(\epsilon_1^2-\sigma^2)^2> \upsilon^2 n\Big)\\
\leq& n^{-1}\sum_{i=1}^{n}\Big(1+\big\|\bm{\xi}_i^{(0)}\big\|^2\Big)^2\mathbf{E}\bigg[\Big(\epsilon_1^2+(\epsilon_i^2-\sigma^2)^2\Big)^{3/2}\mathbf{1}\Big(\epsilon_1^2+(\epsilon_1^2-\sigma^2)^2> k\upsilon^2 n^{3/4}\Big)\bigg]\\
=& o(1).
\end{align*}
Now consider (\ref{eqn:orIEEres2}). Note that for any $|\bm{\alpha}|\leq q+1$,  $|D^{\bm{\alpha}}\mathbf{E}\exp(i\bm{t}^\prime\bm{R}_{1n}^{\dagger})|$ is bounded above by a sum of $n^{|\alpha|}$-terms, each of which is bounded above by
\begin{align}\label{eqn:orIEEres3}
C(\alpha) \cdot n^{-|\bm{\alpha}|/2} \max\{\mathbf{E}\|\bm{Z}_{k}-\mathbf{E}\bm{Z}_k\|^{|\bm{\alpha}|}: k \in \bm{I}_n\} \cdot \prod_{k \in \bm{I}^{\mathsf{c}}_n}|\mathbf{E}\exp(i\bm{t}^\prime\bm{Z}_{k}/\sqrt{n})|
\end{align}
where $\bm{I}_n\subset \{1,\dots, n\}$ is of size $|\bm{\alpha}|$ and $\bm{I^{\mathsf{c}}_n}=\{1,\dots,n\}\backslash \bm{I}_n$ and $C(\bm{\alpha})$ is a constant which depends only on $\bm{\alpha}$. Now for any $\omega>0$ and $\bm{t}\in \mathcal{R}^{q}$, define the set $$\bm{B}_n(\bm{t},\omega) = \Big\{k:1\leq k \leq n  \; \text{and}\; \big(\bm{t}^\prime\bm{\xi}_k^{(0)}\big)^2 > \omega^2\Big\}.$$ Hence for any $\bm{t}\in \mathcal{R}^{q+1}$ writing $\bm{t}=\Big(\bm{t}_q^\prime,t_{q+1}\Big)$, we have
\begin{align*}
&\sup\Bigg{\{}\prod_{k \in \bm{I}^{\mathsf{c}}_n}|\mathbf{E}\exp(i\bm{t}^\prime\bm{Z}_{k}/\sqrt{n})|:\|\bm{t}\|\geq \upsilon \sqrt{n} \Bigg{\}} \\
=&\sup\Bigg{\{}\prod_{k \in \bm{I}^{\mathsf{c}}_n}|\mathbf{E}\exp(i\bm{t}^\prime\bm{Z}_{k})|:\|\bm{t}\|^2\geq \upsilon^2  \Bigg{\}} \\
\leq & \max\Bigg\{\sup\bigg{\{}\prod_{k \in \bm{I}^{\mathsf{c}}_n\cap \bm{B}_n\Big(\dfrac{\bm{t}_q}{\|\bm{t}_q\|}, \upsilon/\sqrt{2}\Big)}\Big[|\mathbf{E}\exp\Big(i\bm{t}_q^\prime\bm{\xi}_{k}^{(0)}\epsilon_k\Big)|+\mathbf{P}\Big(\epsilon_1^2+(\epsilon_1^2-\sigma^2)^2> k\upsilon^2 n^{3/4}\Big)\Big]\\
&:\|\bm{t}_q\|\geq \upsilon/\sqrt{2}  \bigg{\}}, \sup\bigg{\{}\prod_{k \in \bm{I}^{\mathsf{c}}_n}\Big[|\mathbf{E}\exp\Big(it_{q+1}(\epsilon_k^2-\sigma^2)\Big)|+\mathbf{P}\Big(\epsilon_1^2+(\epsilon_1^2-\sigma^2)^2> k\upsilon^2 n^{3/4}\Big)\Big]\\
&:t_{q+1}\geq \upsilon/\sqrt{2}  \bigg{\}}\Bigg\}
\end{align*}
Now since $\big|\bm{I}^{\mathsf{c}}_n\big|\geq \Big|\bm{I}^{\mathsf{c}}_n\cap \bm{B}_n\Big(\dfrac{\bm{t}_q}{\|\bm{t}_q\|}, \upsilon/\sqrt{2}\Big)\Big|\geq \big|\bm{B}_n\Big(\dfrac{\bm{t}_q}{\|\bm{t}_q\|}, \upsilon/\sqrt{2}\Big)\big|-|\alpha|$, due to Cramer's condition (A.5)(ii), for some $\theta \in (0,1)$, we have
\begin{align}\label{eqn:orIEEres4}
&\sup\Bigg{\{}\prod_{k \in \bm{I}^{\mathsf{c}}_n}|\mathbf{E}\exp(i\bm{t}^\prime\bm{Z}_{k}/\sqrt{n})|:\|\bm{t}\|\geq \upsilon \sqrt{n} \Bigg{\}} \leq \theta^{\Big|\bm{B}_n\Big(\dfrac{\bm{t}_q}{\|\bm{t}_q\|}, \upsilon/\sqrt{2}\Big)\Big|-\big|\alpha\big|}
\end{align}
Next note that for any $\bm{u}\in \mathcal{R}^q$ with $\|u\|=1$, due to conditions (A.3)(ii) \& (A.3)(iii), for sufficiently large $n$ we have 
\begin{align*}
\dfrac{n\delta}{2}\leq & \sum_{i=1}^{n}\Big|\bm{u}^\prime\bm{\xi}_i^{(0)}\Big|^2\\
\leq & \max\Big\{\big\|\bm{\xi}_i^{(0)}\big\|^2:1\leq i \leq n\Big\}\cdot |\bm{B}_n(\bm{u},\omega)|+\Big(n-|\bm{B}_n(\bm{u},\omega)|\Big)\cdot \omega^2\\
\leq & n^{1/4}\cdot |\bm{B}_n(\bm{u},\omega)|+\Big(n-|\bm{B}_n(\bm{u},\omega)|\Big)\cdot \omega^2\\
\leq & k\cdot n^{1/4}\cdot |\bm{B}_n(\bm{u},\omega)| + n\omega^2
\end{align*}
which implies $|\bm{B}_n(\bm{u},\omega)|\geq k_1\cdot n^{3/4}$ whenever $\omega < \sqrt{\delta/2}$. Therefore taking $\delta_4=\sqrt{\delta/3}$, (\ref{eqn:orIEEres2}) follows from (\ref{eqn:orIEEres3}) and (\ref{eqn:orIEEres4}). Now we are ready to apply the transformation technique of Bhattacharya and Ghosh (1978) to find two term EE of $\tilde{\bm{R}}_n$ and hence of $\bm{R}_n$. Now the first three cumulants of $\bm{t}'\tilde{\bm{R}}_{n}$ are given by

\begin{align*}
\kappa_1\big(\bm{t}'\tilde{\bm{R}}_{n}\big)&=\mathbf{E}\big(\bm{t}'\tilde{\bm{R}}_{n}\big)=-n^{-1/2}\dfrac{\mu_3}{2\sigma^3}\sum_{|\bm{\alpha}|=1}\bm{t}^{\bm{\alpha}}\bar{\bm{\xi}}_n^{(0)}(\bm{\alpha})+o(n^{-1/2})\\
\kappa_2\big(\bm{t}'\tilde{\bm{R}}_{n}\big)&=\mathbf{Var}\big(\bm{t}'\tilde{\bm{R}}_{n}\big)=\bm{t}^\prime\Big(n^{-1}\sum_{i=1}^{n}\bm{\xi}_i^{(0)}\bm{\xi}_i^{(0)\prime}\Big)\bm{t}=\bm{t}'\bm{\Sigma}_n\bm{t}\\
\kappa_3\big(\bm{t}'\tilde{\bm{R}}_{n}\big)&=\mathbf{E}\big(\bm{t}'\tilde{\bm{R}}_{n}\big)^3-3\mathbf{E}\big(\bm{t}'\tilde{\bm{R}}_{n}\big)^2. \mathbf{E}\big(\bm{t}'\tilde{\bm{R}}_{n}\big)+2\Big(\mathbf{E}\big(\bm{t}'\tilde{\bm{R}}_{n}\big)\Big)^3\\
&= n^{-1/2}\dfrac{\mu_3}{\sigma^3}\Big[\sum_{|\bm{\alpha}|=3}\bm{t}^{\bm{\alpha}}\bar{\bm{\xi}}_n^{(0)}(\bm{\alpha})-3\sum_{|\bm{\alpha}|=1}\sum_{|\gamma|=2}\bm{t}^{\bm{\alpha}+\bm{\gamma}}\bar{\bm{\xi}}_n^{(0)}(\bm{\alpha})\bar{\bm{\xi}}_n^{(0)}(\bm{\gamma})\Big]
\end{align*}
where $\bar{\bm{\xi}}_n^{(0)}(\bm{\alpha})=n^{-1}\sum_{i=1}^{n}(\bm{\xi}_i^{(0)})^{\bm{\alpha}}$. Therefore the Lebesgue density of two term EE of $\tilde{\bm{R}}_n$ is given by

\begin{align*}
\bm{\psi}_{1n}(\bm{x})=&\phi(\bm{x}:\bm{\Sigma}_n)\Bigg[1+\dfrac{1}{\sqrt{n}}\bigg[ -\dfrac{\mu_3}{2\sigma^3}\sum_{|\bm{\alpha}|=1}\bar{\bm{\xi}}_n({\bm{\alpha}})\chi_{\bm{\alpha}}(\bm{x}:\bm{\Sigma}_n)\\
&+\dfrac{\mu_3}{6\sigma^3}\Big{\{}\sum_{|\bm{\alpha}|=3}\bar{\bm{\xi}}_n({\bm{\alpha}})\chi_{\bm{\alpha}}(\bm{x}:\bm{\Sigma}_n)-3\sum_{|\bm{\alpha}|=3}\sum_{|\bm{\zeta}|=1}\bar{\bm{\xi}}_n({\bm{\alpha}})\bar{\bm{\xi}}_n({\bm{\zeta}})\chi_{\bm{\alpha}+\bm{\zeta}}(\bm{x}:\bm{\Sigma}_n)\Big{\}}\bigg]\Bigg],
\end{align*}
%where  $$(-D)^{\bm{\alpha}}\phi(\bm{y}:\bm{\Sigma}_n)=\chi_{\bm{\alpha}}(\bm{y}:\bm{\Sigma}_n)\phi(\bm{y}:\bm{\Sigma}_n).$$
Hence due to (\ref{eqn:stuor1}) and since $\mathbf{P}(\bm{A}_n\cap \bm{B}_n)=1-o(n^{-1/2})$, we have
\begin{align}\label{eqn:orEE1}
\sup_{\bm{B}\in \mathcal{C}_q}|\mathbf{P}(\bm{R}_n\in \bm{B})-\int_{\bm{B}} \bm{\psi}_{1n}(\bm{x})d\bm{x}|=o(n^{-1/2})
\end{align}
Through the same line of arguments, it can be shown that
\begin{align}\label{eqn:orEE2}
\sup_{\bm{B}\in \mathcal{C}_q}|\mathbf{P}(\check{\bm{R}}_n\in \bm{B})-\int_{\bm{B}} \check{\bm{\psi}}_{1n}(\bm{x})d\bm{x}|=o(n^{-1/2})
\end{align}
where 
\begin{align*}
\check{\bm{\psi}}_{1n}(\bm{x})=&\phi(\bm{x})\Bigg[1+\dfrac{1}{\sqrt{n}}\bigg[ -\dfrac{\mu_3}{2\sigma^3}\sum_{|\bm{\alpha}|=1}\bar{\bm{\xi}}_n^{\dagger(0)}({\bm{\alpha}})\chi_{\bm{\alpha}}(\bm{x})\\
&+\dfrac{\mu_3}{6\sigma^3}\Big{\{}\sum_{|\bm{\alpha}|=3}\bar{\bm{\xi}}^{\dagger(0)}_n({\bm{\alpha}})\chi_{\bm{\alpha}}(\bm{x})-3\sum_{|\bm{\alpha}|=3}\sum_{|\bm{\zeta}|=1}\bar{\bm{\xi}}^{\dagger(0)}_n({\bm{\alpha}})\bar{\bm{\xi}}^{\dagger(0)}_n({\bm{\zeta}})\chi_{\bm{\alpha}+\bm{\zeta}}(\bm{x})\Big{\}}\bigg]\Bigg],
\end{align*}
with $\bar{\bm{\xi}}_n^{\dagger(0)}(\bm{\alpha})=n^{-1}\sum_{i=1}^{n}\Big(\bm{\Sigma}_n^{-1/2}\bm{\xi}_i^{(0)}\Big)^{\bm{\alpha}}$. % and $(-D)^{\bm{\alpha}}\phi(\bm{y})=\chi_{\bm{\alpha}}(\bm{y})\phi(\bm{y}).$

Now let us look into $\bm{R}_n^*$, the residual bootstrapped versions of $\bm{R}_n$ and $\check{\bm{R}}_n^*$, the perturbation bootstrapped versions of $\check{\bm{R}}_n$. Note that similar to (\ref{eqn:stuor1}) we have
\begin{align}
\bm{R}_n^*
=& n^{-1/2}\sum_{i=1}^{n}\bm{\xi}_i^{(0)}\epsilon_i^*/\hat{\sigma}_n -2^{-1}\hat{\sigma}_n^{-3}\Big[n^{-1}\sum_{i=1}^{n}(\epsilon_i^{*2}-\hat{\sigma}_n^2)\Big]n^{-1/2}\sum_{i=1}^{n}\bm{\xi}_i^{(0)}\epsilon_i^* +\bm{Q}_{2n}^*\nonumber\\\label{eqn:stuboot1}
=& \tilde{\bm{R}}_{1n}^* +\bm{Q}_{2n}^*\\
\check{\bm{R}}_n^*
=& n^{-1/2}\sum_{i=1}^{n}\tilde{\bm{\Sigma}}_n^{-1/2}\bm{\xi}_i^{(0)}\hat{\epsilon}_i(G_i^*-\mu_{G^*})\mu_{G^*}^{-1}\nonumber\\
&-2^{-1}\check{\sigma}_n^{-2}\Big[\mu_{G^*}^{-2}n^{-1}\sum_{i=1}^{n}\hat{\epsilon}_i^{2}\big[(G_i^*-\mu_{G^*})^2-\sigma_{G^*}^2\big]\Big]n^{-1/2}\sum_{i=1}^{n}\tilde{\bm{\Sigma}}_n^{-1/2}\bm{\xi}_i^{(0)}\hat{\epsilon}_i(G_i^*-\mu_{G^*})\mu_{G^*}^{-1} +\bm{Q}_{3n}^*\nonumber\\\label{eqn:stuboot2}
=& \check{\bm{R}}_{1n}^* +\bm{Q}_{3n}^*
\end{align}
where due to conditions (A.2)(ii), (A.3)(ii), (A.5)(i), (A.6)(i), (A.7)(iii) and Lemma \ref{lem:betahat}, we have $\mathbf{P}_*\Big(\|\bm{Q}_{2n}^*\|=o(n^{-1/2})\Big)=1-o_p(n^{-1/2})$ and $\mathbf{P}_*\Big(\|\bm{Q}_{3n}^*\|=o(n^{-1/2})\Big)=1-o_p(n^{-1/2})$. Therefore it is enough to find the two term EE of $\tilde{\bm{R}}_{1n}^*$ and $\check{\bm{R}}_{1n}^*$. Due to the Cramer's condition (A.6)(ii) on $\big((G_i^*-\mu_{G^*}), (G_i^*-\mu_{G^*})^2\big)$, one can find the EE of $\check{\bm{R}}_{1n}^*$ through the same line of arguments as in the original case, i.e. first finding two term EE of $\Big(n^{-1/2}\sum_{i=1}^{n}\tilde{\bm{\Sigma}}_n^{-1/2}\bm{\xi}_i^{(0)\prime}\hat{\epsilon}_i(G_i^*-\mu_{G^*}), n^{-1/2}\sum_{i=1}^{n}\hat{\epsilon}_i^{2}\big[(G_i^*-\mu_{G^*})^2-\sigma_{G^*}^2\big]\Big)^\prime$ by applying Theorem 20.6 of Bhattacharya and Rao (1986) and then using the transformation technique of Bhattacharya and Ghosh (1978). However in case of residual bootstrap, one can not use Theorem 20.6 of Bhattacharya and Rao (1986) directly to obtain two term EE of $\Big(n^{-1/2}\sum_{i=1}^{n}\bm{\xi}_i^{(0)\prime}\epsilon_i^*, n^{-1/2}\sum_{i=1}^{n}(\epsilon_i^{*2}-\hat{\sigma}_n^2)\Big)$, since only the conditional Cramer's condition on $(\epsilon_1^*,\epsilon_1^{*2})$, viz Lemma \ref{lem:residualCramer}, holds under (A.5)(ii). Here one needs to use a smoothing kernel which vanishes outside a compact set to take advantage of the conditional Cramer's condition, instead of using the smoothing kernel used in Theorem 20.6 of Bhattacharya and Rao (1986). The arguments of Theorem 20.8 of Bhattacharya and Rao (1986) or Theorem 2 of Babu and Singh (1984) can be followed to obtain the two term EE of $\Big(n^{-1/2}\sum_{i=1}^{n}\bm{\xi}_i^{(0)\prime}\epsilon_i^*, n^{-1/2}\sum_{i=1}^{n}(\epsilon_i^{*2}-\hat{\sigma}_n^2)\Big)$ and then, similar to the original case, one can use the transformation technique of Bhattacharya and Ghosh (1978) to come up with the two term EE of $\bm{R}_n^*$.

Writing $\hat{\mu}_3=n^{-1}\sum_{i=1}^{n}(\hat{\epsilon}_i-\bar{\epsilon}_n)^3$, the first three conditional cumulants of $\bm{t}'\tilde{\bm{R}}_{1n}^*$ given $\{\hat{\epsilon}_1,\dots, \hat{\epsilon}_n\}$ are given by
\begin{align*}
\kappa_1\big(\bm{t}'\tilde{\bm{R}}_{1n}^*\big)&=\mathbf{E}\big(\bm{t}'\tilde{\bm{R}}_{1n}^*\big)=-n^{-1/2}\dfrac{\hat{\mu}_3}{2\hat{\sigma}_n^3}\sum_{|\bm{\alpha}|=1}\bm{t}^{\bm{\alpha}}\bar{\bm{\xi}}_n^{(0)}(\bm{\alpha})+o(n^{-1/2})\\
\kappa_2\big(\bm{t}'\tilde{\bm{R}}_{1n}^*\big)&=\mathbf{Var}\big(\bm{t}'\tilde{\bm{R}}_{1n}^*\big)=\bm{t}^\prime\Big(n^{-1}\sum_{i=1}^{n}\bm{\xi}_i^{(0)}\bm{\xi}_i^{(0)\prime}\Big)\bm{t}=\bm{t}'\bm{\Sigma}_n\bm{t}\\
\kappa_3\big(\bm{t}'\tilde{\bm{R}}_{1n}^*\big)&=\mathbf{E}\big(\bm{t}'\tilde{\bm{R}}_{1n}^*\big)^3-3\mathbf{E}\big(\bm{t}'\tilde{\bm{R}}_{1n}^*\big)^2. \mathbf{E}\big(\bm{t}'\tilde{\bm{R}}_{1n}^*\big)+2\Big(\mathbf{E}\big(\bm{t}'\tilde{\bm{R}}_{1n}^*\big)\Big)^3\\
&= n^{-1/2}\dfrac{\hat{\mu}_3}{\hat{\sigma}_n^3}\Big[\sum_{|\bm{\alpha}|=3}\bm{t}^{\bm{\alpha}}\bar{\bm{\xi}}_n^{(0)}(\bm{\alpha})-3\sum_{|\bm{\alpha}|=1}\sum_{|\gamma|=2}\bm{t}^{\bm{\alpha}+\bm{\gamma}}\bar{\bm{\xi}}_n^{(0)}(\bm{\alpha})\bar{\bm{\xi}}_n^{(0)}(\bm{\gamma})\Big]
\end{align*}
Therefore, the Lebesgue density of two term conditional EE of $\tilde{\bm{R}}_{1n}^*$ is given by
\begin{align*}
\bm{\psi}_{1n}^*(\bm{x})=&\phi(\bm{x}:\bm{\Sigma}_n)\Bigg[1+\dfrac{1}{\sqrt{n}}\bigg[ -\dfrac{\hat{\mu}_3}{2\hat{\sigma}_n^3}\sum_{|\bm{\alpha}|=1}\bar{\bm{\xi}}_n({\bm{\alpha}})\chi_{\bm{\alpha}}(\bm{x}:\bm{\Sigma}_n)\\
&+\dfrac{\hat{\mu}_3}{6\hat{\sigma}_n^3}\Big{\{}\sum_{|\bm{\alpha}|=3}\bar{\bm{\xi}}_n({\bm{\alpha}})\chi_{\bm{\alpha}}(\bm{x}:\bm{\Sigma}_n)-3\sum_{|\bm{\alpha}|=3}\sum_{|\bm{\zeta}|=1}\bar{\bm{\xi}}_n({\bm{\alpha}})\bar{\bm{\xi}}_n({\bm{\zeta}})\chi_{\bm{\alpha}+\bm{\zeta}}(\bm{x}:\bm{\Sigma}_n)\Big{\}}\bigg]\Bigg].
\end{align*}
%where  $$(-D)^{\bm{\alpha}}\phi(\bm{y}:\bm{\Sigma}_n)=\chi_{\bm{\alpha}}(\bm{y}:\bm{\Sigma}_n)\phi(\bm{y}:\bm{\Sigma}_n).$$

Again the first three conditional cumulants of $\bm{t}'\check{\bm{R}}_{1n}^*$ given $\{\hat{\epsilon}_1,\dots, \hat{\epsilon}_n\}$ are given by

$\kappa_1\big(\bm{t}'\check{\bm{R}}_{1n}^*\big)=-\dfrac{1}{\sqrt{n}}. \dfrac{1}{2\check{\sigma}_n^2}\sum_{|\bm{\alpha}|=1}\bm{t}^{\bm{\alpha}}\bar{\bm{\xi}}_n^{*(3)}({\bm{\alpha}}) +o_p(n^{1/2})$

$\kappa_2\big(\bm{t}'\check{\bm{R}}_{1n}^*\big)=\mathbf{Var_*}\big(\bm{t}'\bm{R}_{1n}^*\big)=\bm{t}'\bm{t}+o_p(n^{-1/2})$

$\kappa_3\big(\bm{t}'\check{\bm{R}}_{1n}^*\big)=\mathbf{E_*}\big(\bm{t}'\bm{R}_{1n}^*\big)^3-3\mathbf{E_*}\big(\bm{t}'\bm{R}_{1n}^*\big)^2. \mathbf{E_*}\big(\bm{t}'\bm{R}_{1n}^*\big)+2\Big(\mathbf{E_*}\big(\bm{t}'\bm{R}_{1n}^*\big)\Big)^3$\\
\hspace*{10mm} $=\dfrac{1}{\sqrt{n}}\bigg[\sum_{|\bm{\alpha}|=3}\bm{t}^{\bm{\alpha}}\bar{\bm{\xi}}_n^{*(1)}({\bm{\alpha}})-\dfrac{3}{\hat{\sigma}_n^2}\sum_{|\bm{\alpha}|=1}\sum_{|\bm{\zeta}|=2}\bm{t}^{\bm{\alpha}+\bm{\zeta}}\bar{\bm{\xi}}_n^{*(3)}({\bm{\alpha}})\bar{\bm{\xi}}_n^{*(1)}({\bm{\zeta}})\Bigg]+o_p(n^{-1/2})$,

where $\bar{\bm{\xi}}_{n}^{*(j)}(\bm{\alpha})=n^{-1}\sum_{i=1}^{n}\Big(\tilde{\bm{\Sigma}}_n^{-1/2}\bm{\xi}_i^{(0)}\hat{\epsilon}_i^j\Big)^{\bm{\alpha}}$, $j=0,1,2,3$.

Therefore the Lebesgue density of two term conditional EE of $\tilde{\bm{R}}_{2n}^*$ is given by
\begin{align*}
\check{\bm{\psi}}_{1n}^*(\bm{x})=&\phi(\bm{x})\Bigg[1+\dfrac{1}{\sqrt{n}}\bigg[ \dfrac{1}{6}\sum_{|\bm{\alpha}|=3}\bar{\bm{\xi}}_n^{*(1)}({\bm{\alpha}})\chi_{\bm{\alpha}}(\bm{x})\\ 
&-\dfrac{1}{2\hat{\sigma}_n^2}\Big{\{}\sum_{|\bm{\alpha}|=1}\bar{\bm{\xi}}_n^{*(3)}({\bm{\alpha}})\chi_{\bm{\alpha}}(\bm{x})
+\sum_{|\bm{\alpha}|=1}\sum_{|\bm{\zeta}|=2}\bar{\bm{\xi}}_n^{*(3)}({\bm{\alpha}})\bar{\bm{\xi}}_n^{*(1)}({\bm{\zeta}})\chi_{\bm{\alpha}+\bm{\zeta}}(\bm{x})\Big{\}}\bigg]\Bigg].
\end{align*}
%where  $$(-D)^{\bm{\alpha}}\phi(\bm{y})=\chi_{\bm{\alpha}}(\bm{y})\phi(\bm{y}).$$

Hence due to (\ref{eqn:stuboot1}), (\ref{eqn:stuboot2})  and since $\mathbf{P}_*(\bm{A}_n^*\cap \bm{B}_n^*)=1-o_p(n^{-1/2})$, we have
\begin{align}\label{eqn:bootEE1}
&\sup_{\bm{B}\in \mathcal{C}_q}|\mathbf{P}_*(\bm{R}_n^*\in \bm{B})-\int_{\bm{B}} \bm{\psi}^*_{1n}(\bm{x})d\bm{x}|=o_p(n^{-1/2})\;\;\;\; and
\end{align}
\begin{align}\label{eqn:bootEE2}
\sup_{\bm{B}\in \mathcal{C}_q}|\mathbf{P}_*(\check{\bm{R}}_n^*\in \bm{B})-\int_{\bm{B}} \check{\bm{\psi}}^*_{1n}(\bm{x})d\bm{x}|=o_p(n^{-1/2})
\end{align}

Now note that $|\hat{\mu}_3-\mu_3|=o_p(1)$, $|\hat{\sigma}_n^2-\sigma^2|=o_p(1)$ $\|\tilde{\Sigma}_n-\sigma^2\Sigma_n\|=o_p(1)$. Theorem X.3.8 of Bhatia(1996) and $\|\tilde{\Sigma}_n-\sigma^2\Sigma_n\|=o_p(1)$ imply that $\|\tilde{\Sigma}_n^{1/2}-\sigma\Sigma_n^{1/2}\|=o_p(1)$, which again implies $\|\tilde{\Sigma}_n^{-1/2}-\sigma^{-1}\Sigma_n^{-1/2}\|=o_p(1)$ by noting that
$$A^{-1/2} = B^{-1/2} + B^{-1/2}(B^{1/2} - A^{1/2})B^{-1/2}
+B^{-1/2}(B^{1/2} - A^{1/2})B^{-1/2}(B^{1/2} -A^{1/2})A^{-1/2}.$$ Hence comparing (\ref{eqn:orEE1}) \& (\ref{eqn:bootEE1}) and comparing (\ref{eqn:orEE2}) \& (\ref{eqn:bootEE2}), Theorem  \ref{thm:bootI} follows due to conditions (A.3)(i) and (A.3)(iii) with $r=4$.

\textbf{Proof of Theorem \ref{thm:bootII}}:
In the proof of Theorem \ref{thm:oracleII}, we have seen that on the set $A_{2n}$,
\begin{align*}
\bm{T}_n &= \bm{D}_n^{(1)}\bm{C}_{11,n}^{-1}\Big[\bm{W}_n^{(1)}-\dfrac{\sqrt{n}\lambda_n}{2}\bm{s}_n^{(1)}-\dfrac{\lambda_n}{2}\bm{L}_n^{(1)}\Big] +\bm{Q}_{1n}\\
&= \bm{D}_n^{(1)}\bm{C}_{11,n}^{-1}\Big[\bm{W}_n^{(1)}-\dfrac{\sqrt{n}\lambda_n}{2}\bm{s}_n^{(1)}\Big]-\dfrac{\lambda_n}{2}\bm{D}_n^{(1)}\bm{C}_{11,n}^{-1}\bm{L}_n^{(1)} +\bm{Q}_{1n}\\
&=\Big[\bm{D}_n^{(1)}\bm{C}_{11,n}^{-1}\bm{W}_n^{(1)}+\bm{b}_n\Big]+\Big[-\dfrac{\lambda_n}{2}\bm{D}_n^{(1)}\bm{C}_{11,n}^{-1}\bm{L}_n^{(1)} +\bm{Q}_{1n}\Big]\\
&= \tilde{\bm{T}}_{1n}+ \tilde{\bm{Q}}_{1n},
\end{align*}
where  $\tilde{\bm{Q}}_{1n}=-\dfrac{\lambda_n}{2}\bm{D}_n^{(1)}\bm{C}_{11,n}^{-1}\bm{L}_n^{(1)} +\bm{Q}_{1n}$, $\bm{s}_n^{(1)}=(s_{1,n},\ldots, s_{p_0,n})$, $s_{j,n}=sgn(\beta_{j,n})\tilde{P}^\prime(|\beta_{j,n}|)$ and $\bm{L}_n^{(1)} = (L_{1,n},\ldots, L_{p_0,n})$ with 
\begin{align*}
L_{j,n}&= \sqrt{n}(\tilde{\beta}_{j,n}-\beta_{j,n})sgn(\beta_{j,n})\tilde{P}^{\prime\prime}(|\beta_{j,n}|).
%&= n^{-1/2}\sum_{i=1}^{n}\big[(\bm{C}_{n}^{-1})_{j.}\bm{x}_{i}\epsilon_i\big]sgn(\beta_{j,n})\tilde{P}^{\prime\prime}(|\beta_{j,n}|).
\end{align*}
Now we know that $P(||\bm{Q}_{1n}||=o(n^{-1/2}))=1-o(n^{-1/2})$. Again due to conditions (A.2)(iii), (A.3)(i), (A.7)(i) and the condition (\ref{eqn:initialmoderate}) on the initial estimator $\tilde{\bm{\beta}}_n$, on the set $\bm{A}_{2n}$ we have  $$\|\tilde{\bm{Q}}_{1n}-\bm{Q}_{1n}\|=o(n^{-1/2}).$$
Now since $\mathbf{P}\big(|n^{-1}\sum_{i=1}^{n}(\epsilon_i^2-\sigma^2)|>k\cdot \sqrt{\log n/n}\big)=o_p(n^{-1/2})$ and $\|\bm{b}_n\|=O(n^{-\delta})$, we have 
\begin{align}
\bm{R}_n=&\bm{T}_n/\hat{\sigma}_n\nonumber\\
=&\Big[\sigma^{-1}-2^{-1}\sigma^{-3}(\hat{\sigma}_n^2-\sigma^2)\Big]\tilde{\bm{T}}_{1n} + \tilde{\bm{Q}}_{2n}\nonumber\\
=&\Big[\sigma^{-1}-2^{-1}\sigma^{-3}\Big(n^{-1}\sum_{i=1}^{n}(\epsilon_i^2-\sigma^2)\Big)\Big]\bm{D}_n^{(1)}\bm{C}_{11,n}^{-1}\bm{W}_n^{(1)} + \bm{b}_n/\sigma + \tilde{\bm{Q}}_{3n}\nonumber\\
=&\Big[\tilde{\bm{R}}_n + \bm{b}_n/\sigma\Big] + \tilde{\bm{Q}}_{3n}\nonumber\\\label{eqn:orIIstuor}
=& \tilde{\bm{R}}_{1n} +\tilde{\bm{Q}}_{3n}\;\;\;\; \text{(say)},
\end{align}
where $\mathbf{P}\big(\|\tilde{\bm{Q}}_{2n}\|+\|\tilde{\bm{Q}}_{3n}\|=o(n^{-1/2})\big)=1-o(n^{-1/2})$. Now in the proof of Theorem 5, we have seen that
\begin{align*}
\sup_{\bm{B}\in \mathcal{C}_q}|\mathbf{P}(\tilde{\bm{R}}_{n}\in \bm{B})-\int_{\bm{B}} \bm{\psi}_{1n}(\bm{x})d\bm{x}|=o(n^{-1/2})
\end{align*}
where 
\begin{align*}
\bm{\psi}_{1n}(\bm{x})=&\phi(\bm{x}:\bm{\Sigma}_n)\Bigg[1+\dfrac{1}{\sqrt{n}}\bigg[ -\dfrac{\mu_3}{2\sigma^3}\sum_{|\bm{\alpha}|=1}\bar{\bm{\xi}}_n({\bm{\alpha}})\chi_{\bm{\alpha}}(\bm{x}:\bm{\Sigma}_n)\\
&+\dfrac{\mu_3}{6\sigma^3}\Big{\{}\sum_{|\bm{\alpha}|=3}\bar{\bm{\xi}}_n({\bm{\alpha}})\chi_{\bm{\alpha}}(\bm{x}:\bm{\Sigma}_n)-3\sum_{|\bm{\alpha}|=3}\sum_{|\bm{\zeta}|=1}\bar{\bm{\xi}}_n({\bm{\alpha}})\bar{\bm{\xi}}_n({\bm{\zeta}})\chi_{\bm{\alpha}+\bm{\zeta}}(\bm{x}:\bm{\Sigma}_n)\Big{\}}\bigg]\Bigg],
\end{align*}
Hence if  $r=\min\{a\in \mathcal{N}:||\bm{b}_n||^{a+1}=O(n^{-1/2})\}$, $\mathcal{N}$ being the set of natural numbers, then we have
 \begin{align}\label{eqn:orIIorEE1}
\sup_{\bm{B}\in \mathcal{C}_q}\Big|\mathbf{P}\big(\bm{R}_{n}\in \bm{B}\big)-\int_{\bm{B}}\bm{\psi}_{2n}(\bm{y})d\bm{y}\Big|=o\big(n^{-1/2}\big).
\end{align}
where 
\begin{align*}
\bm{\psi}_{1n}(\bm{x})=&\phi(\bm{x}:\bm{\Sigma}_n)\Bigg[1+\sum_{|\bm{\alpha}|=1}^{r+1}\big(-\bm{b}_n/\sigma\big)^{\bm{\alpha}}\chi_{\bm{\alpha}}(\bm{x}:\bm{\Sigma}_n)+\dfrac{1}{\sqrt{n}}\bigg[ -\dfrac{\mu_3}{2\sigma^3}\sum_{|\bm{\alpha}|=1}\bar{\bm{\xi}}_n({\bm{\alpha}})
\chi_{\bm{\alpha}}(\bm{x}:\bm{\Sigma}_n)\\
&+\dfrac{\mu_3}{6\sigma^3}\Big{\{}\sum_{|\bm{\alpha}|=3}\bar{\bm{\xi}}_n({\bm{\alpha}})\chi_{\bm{\alpha}}(\bm{x}:\bm{\Sigma}_n)-3\sum_{|\bm{\alpha}|=3}\sum_{|\bm{\zeta}|=1}\bar{\bm{\xi}}_n({\bm{\alpha}})\bar{\bm{\xi}}_n({\bm{\zeta}})\chi_{\bm{\alpha}+\bm{\zeta}}(\bm{x}:\bm{\Sigma}_n)\Big{\}}\bigg]\Bigg],
\end{align*}
with $\bar{\bm{\xi}}_n^{(0)}(\bm{\alpha})=n^{-1}\sum_{i=1}^{n}(\bm{\xi}_i^{(0)})^{\bm{\alpha}}$.
%and $\chi_{\bm{\alpha}}(\bm{y}:\sigma^2\bm{\Sigma}_n)\phi(\bm{y}:\sigma^2\bm{\Sigma}_n)$ being defined by the identity $$(-D)^{\bm{\alpha}}\phi(\bm{y}:\bm{\Sigma}_n)=\chi_{\bm{\alpha}}(\bm{y}:\bm{\Sigma}_n)\phi(\bm{y}:\bm{\Sigma}_n).$$

Through the same line of arguments, it can be shown that
\begin{align}\label{eqn:orIIorEE2}
\sup_{\bm{B}\in \mathcal{C}_q}|\mathbf{P}(\check{\bm{R}}_n\in \bm{B})-\int_{\bm{B}} \check{\bm{\psi}}_{2n}(\bm{x})d\bm{x}|=o(n^{-1/2})
\end{align}
where if $\check{\bm{b}}_n=\bm{\Sigma}_n^{-1/2}\bm{b}_n$ and $\check{r}=\min\{a\in \mathcal{N}:\big\|\check{\bm{b}}_n\big\|^{a+1}=O(n^{-1/2})\}$, $\mathcal{N}$ being the set of natural numbers, then
\begin{align*}
\check{\bm{\psi}}_{1n}(\bm{x})=&\phi(\bm{x})\Bigg[1+\sum_{|\bm{\alpha}|=1}^{\check{r}+1}\big(-\check{\bm{b}}_n/\sigma\big)^{\bm{\alpha}}\chi_{\bm{\alpha}}(\bm{x})+\dfrac{1}{\sqrt{n}}\bigg[ -\dfrac{\mu_3}{2\sigma^3}\sum_{|\bm{\alpha}|=1}\bar{\bm{\xi}}_n^{\dagger(0)}({\bm{\alpha}})\chi_{\bm{\alpha}}(\bm{x})\\
&+\dfrac{\mu_3}{6\sigma^3}\Big{\{}\sum_{|\bm{\alpha}|=3}\bar{\bm{\xi}}^{\dagger(0)}_n({\bm{\alpha}})\chi_{\bm{\alpha}}(\bm{x})-3\sum_{|\bm{\alpha}|=3}\sum_{|\bm{\zeta}|=1}\bar{\bm{\xi}}^{\dagger(0)}_n({\bm{\alpha}})\bar{\bm{\xi}}^{\dagger(0)}_n({\bm{\zeta}})\chi_{\bm{\alpha}+\bm{\zeta}}(\bm{x})\Big{\}}\bigg]\Bigg],
\end{align*}
with $\bar{\bm{\xi}}_n^{\dagger(0)}(\bm{\alpha})=n^{-1}\sum_{i=1}^{n}\Big(\bm{\Sigma}_n^{-1/2}\bm{\xi}_i^{(0)}\Big)^{\bm{\alpha}}$.% and $(-D)^{\bm{\alpha}}\phi(\bm{y})=\chi_{\bm{\alpha}}(\bm{y})\phi(\bm{y}).$

Now look into the bootstrap versions. Note that similar to (\ref{eqn:orIIstuor}) we have
\begin{align}
\bm{R}_n^*
=& \Big[\hat{\sigma}_n^{-1} -2^{-1}\hat{\sigma}_n^{-3}\big[n^{-1}\sum_{i=1}^{n}(\epsilon_i^{*2}-\hat{\sigma}_n^2)\big]\Big]n^{-1/2}\sum_{i=1}^{n}\bm{\xi}_i^{(0)}\epsilon_i^* +\hat{\bm{b}}_n/\hat{\sigma}_n + \bm{Q}_{4n}^*\nonumber\\\label{eqn:orIIstuboot1}
=& \tilde{\bm{R}}_{2n}^* +\bm{Q}_{4n}^*\;\;\;\; (say)\\
\check{\bm{R}}_n^*
=& \bm{Q}_{5n}^* + \tilde{\bm{\Sigma}}_n^{-1/2}\hat{\bm{b}}_n+n^{-1/2}\sum_{i=1}^{n}\tilde{\bm{\Sigma}}_n^{-1/2}\bm{\xi}_i^{(0)}\hat{\epsilon}_i(G_i^*-\mu_{G^*})\mu_{G^*}^{-1}\nonumber\\
&-2^{-1}\check{\sigma}_n^{-2}\Big[\mu_{G^*}^{-2}n^{-1}\sum_{i=1}^{n}\hat{\epsilon}_i^{2}\big[(G_i^*-\mu_{G^*})^2-\sigma_{G^*}^2\big]\Big]n^{-1/2}\sum_{i=1}^{n}\tilde{\bm{\Sigma}}_n^{-1/2}\bm{\xi}_i^{(0)}\hat{\epsilon}_i(G_i^*-\mu_{G^*})\mu_{G^*}^{-1} \nonumber\\\label{eqn:orIIstuboot2}
=& \bm{Q}_{5n}^* + \check{\bm{R}}_{2n}^*\;\;\;\; (say)
\end{align}
where $\hat{\bm{b}}_n= \dfrac{-\lambda_n}{2\sqrt{n}}\hat{\bm{D}}_n^{(1)}\hat{\bm{C}}_{11,n}^{-1}\hat{\bm{s}}_n^{(1)}$
with $\hat{\bm{s}}_n^{(1)}=(\hat{s}_{1,n},\ldots, \hat{s}_{p_0,n})$ and $\hat{s}_{j,n}=sgn(\hat{\beta}_{j,n})\tilde{P}^\prime(|\hat{\beta}_{j,n}|)$. $\hat{\bm{D}}_n^{(1)}$, $\hat{\bm{C}}_{11,n}$ and $\hat{\bm{x}}_i^{(1)}$ are as defined in Section \ref{sec:mainresults}. Due to conditions (A.2)(ii), (A.3)(ii), (A.5)(i), (A.6)(i), (A.7)(iii) and Lemma \ref{lem:betahat}, we have $\mathbf{P}_*\Big(\|\bm{Q}_{4n}^*\|=o(n^{-1/2})\Big)=1-o_p(n^{-1/2})$ and $\mathbf{P}_*\Big(\|\bm{Q}_{5n}^*\|=o(n^{-1/2})\Big)=1-o_p(n^{-1/2})$.

Now we know that $|\hat{\sigma}_n^2-\sigma^2|=O_p(n^{-1/2})$, $\|\hat{\bm{b}}_n\|=O_p(n^{-\delta})$ and by Lemma \ref{lem:Sigma}, $\|\tilde{\bm{\Sigma}}_n-\sigma^2\bm{\Sigma}_n\|=O_p(n^{-1/2})$. Hence we can define $\tilde{r}=\min\{a\in \mathcal{N}:\|\hat{\bm{b}}_n/\hat{\sigma}_n\|^{a+1}=o_p(n^{-1/2})\}$ and $\check{r}=\min\{a\in \mathcal{N}:\big\|\tilde{\bm{\Sigma}}_n^{-1/2}\hat{\bm{b}}_n\big\|^{a+1}=o_p(n^{-1/2})\}$, $\mathcal{N}$ being the set of natural numbers. Again note that condition (A.6)(ii) imposes conditional Cramer's condition on $\Big((G_1^*-\mu_{G^*}),(G_1^*-\mu{G^*})^2\Big)$ and Lemma \ref{lem:residualCramer} imposes Cramer's condition on $(\epsilon_1^*,\epsilon_1^{*2})$. Therefore through the same line of arguments as in Theorem 5, we have
 \begin{align}
&\sup_{\bm{B}\in \mathcal{C}_q}\Big|\mathbf{P}_*\big(\bm{R}_{n}^*\in \bm{B}\big)-\int_{\bm{B}}\tilde{\bm{\psi}}_{2n}^*(\bm{y})d\bm{y}\Big|=o_p\big(n^{-1/2}\big)\;\;\; \text{and}\nonumber\\\label{eqn:orIIbootEE}
&\sup_{\bm{B}\in \mathcal{C}_q}\Big|\mathbf{P}_*\big(\check{\bm{R}}_{n}^*\in \bm{B}\big)-\int_{\bm{B}}\check{\bm{\psi}}_{2n}^{**}(\bm{y})d\bm{y}\Big|=o_p\big(n^{-1/2}\big)
\end{align}
where 
\begin{align*}
\bm{\psi}_{2n}^*(\bm{x})=&\phi(\bm{x}:\bm{\Sigma}_n)\Bigg[1+\sum_{|\bm{\alpha}|=1}^{\tilde{r}}\big(-\hat{\bm{b}}_n/\hat{\sigma}_n\big)^{\bm{\alpha}}\chi_{\bm{\alpha}}(\bm{x}:\bm{\Sigma}_n)+\dfrac{1}{\sqrt{n}}\bigg[ -\dfrac{\hat{\mu}_3}{2\hat{\sigma}_n^3}\sum_{|\bm{\alpha}|=1}\bar{\bm{\xi}}_n({\bm{\alpha}})\chi_{\bm{\alpha}}(\bm{x}:\bm{\Sigma}_n)\\
&+\dfrac{\hat{\mu}_3}{6\hat{\sigma}_n^3}\Big{\{}\sum_{|\bm{\alpha}|=3}\bar{\bm{\xi}}_n({\bm{\alpha}})\chi_{\bm{\alpha}}(\bm{x}:\bm{\Sigma}_n)-3\sum_{|\bm{\alpha}|=3}\sum_{|\bm{\zeta}|=1}\bar{\bm{\xi}}_n({\bm{\alpha}})\bar{\bm{\xi}}_n({\bm{\zeta}})\chi_{\bm{\alpha}+\bm{\zeta}}(\bm{x}:\bm{\Sigma}_n)\Big{\}}\bigg]\Bigg]
\end{align*}
and
\begin{align*}
\check{\bm{\psi}}_{2n}^*(\bm{x})=&\phi(\bm{x})\Bigg[1+\sum_{|\bm{\alpha}|=1}^{\check{r}}\big(-\tilde{\bm{\Sigma}}_n^{-1/2}\hat{\bm{b}}_n\big)^{\bm{\alpha}}\chi_{\bm{\alpha}}(\bm{x})+\dfrac{1}{\sqrt{n}}\bigg[ \dfrac{1}{6}\sum_{|\bm{\alpha}|=3}\bar{\bm{\xi}}_n^{*(1)}({\bm{\alpha}})\chi_{\bm{\alpha}}(\bm{x})\\ 
&-\dfrac{1}{2\hat{\sigma}_n^2}\Big{\{}\sum_{|\bm{\alpha}|=1}\bar{\bm{\xi}}_n^{*(3)}({\bm{\alpha}})\chi_{\bm{\alpha}}(\bm{x})
+\sum_{|\bm{\alpha}|=1}\sum_{|\bm{\zeta}|=2}\bar{\bm{\xi}}_n^{*(3)}({\bm{\alpha}})\bar{\bm{\xi}}_n^{*(1)}({\bm{\zeta}})\chi_{\bm{\alpha}+\bm{\zeta}}(\bm{x})\Big{\}}\bigg]\Bigg].
\end{align*}
%with  $$(-D)^{\bm{\alpha}}\phi(\bm{y}:\bm{\Sigma}_n)=\chi_{\bm{\alpha}}(\bm{y}:\bm{\Sigma}_n)\phi(\bm{y}:\bm{\Sigma}_n)\;\;\; and$$   $$(-D)^{\bm{\alpha}}\phi(\bm{y})=\chi_{\bm{\alpha}}(\bm{y})\phi(\bm{y}).$$

Now note that $|\sigma_n^2-\sigma^2|=O_p(n^{-1/2})$, $\|\tilde{\bm{\Sigma}}_n-\sigma^2\bm{\Sigma}_n\|=O_p(n^{-1/2})$, $\|\bm{b}_n\|=O(n^{-\delta})$, $\|\hat{\bm{b}}_n\|=O_p(n^{-1/2})$ and $\|\hat{\bm{b}}_n-\bm{b}_n\|=o_p(n^{-1/2})$. Since one can have
\begin{align*}
&\|\hat{\bm{b}}_n/\hat{\sigma}_n-\bm{b}_n/\sigma\|\leq k\cdot \|\hat{\bm{b}}_n\|\cdot |\sigma_n^2-\sigma^2| + \sigma^{-1}\|\hat{\bm{b}}_n-\bm{b}_n\|\\
and\;\;\;\; &\|\tilde{\bm{\Sigma}}_n^{-1/2}\hat{\bm{b}}_n-\sigma^{-1}\check{\bm{b}}_n\|\leq k\cdot \|\hat{\bm{b}}_n\|\cdot \|\tilde{\bm{\Sigma}}_n^{-1/2}-\sigma^{-1}\bm{\Sigma}_n^{-1/2}\| + \sigma^{-1}\|\bm{\Sigma}_n^{-1/2}\|\cdot\|\hat{\bm{b}}_n-\bm{b}_n\|,
\end{align*}
condition (A.3)(i) and Theorem X.3.8 of Bhatia(1996) imply 
$$\|\hat{\bm{b}}_n/\hat{\sigma}_n-\bm{b}_n/\sigma\|+\|\tilde{\bm{\Sigma}}_n^{-1/2}\hat{\bm{b}}_n-\sigma^{-1}\check{\bm{b}}_n\|=o_p(n^{-1/2})$$
Therefore Theorem \ref{thm:bootII} follows by comparing (\ref{eqn:orIIorEE1}) \& (\ref{eqn:orIIorEE2}) with (\ref{eqn:orIIbootEE}).

\textbf{Proof of Theorem \ref{thm:bootIII}}:
Note that on the set $\bm{A}_{3n}$, for sufficiently large $n$ we have,
\begin{align}
\bm{T}_n-\bm{b}_n^{\dagger}
&= n^{-1/2}\sum_{i=1}^{n}\bm{D}_n^{(1)}\bm{C}_{11,n}^{-1}\bm{x}_i^{(1)}\epsilon_i -\dfrac{\lambda_n}{2\sqrt{n}}\bm{D}_n^{(1)}\bm{C}_{11,n}^{-1}\bm{s}_n^{\dagger(1)}+\bm{D}_n^{(1)}\bm{C}_{11,n}^{-1}\hat{\bm{s}}_n^{\dagger(1)}\dfrac{\lambda_n}{2\sqrt{n}}\nonumber\\
&=n^{-1/2}\sum_{i=1}^{n}\bm{\xi}_i^{(0)}\epsilon_i+\dfrac{\lambda_n}{2\sqrt{n}}\bm{D}_n^{(1)}\bm{C}_{11,n}^{-1}\big(\hat{\bm{s}}_n^{(1)}-\bm{s}_n^{(1)}\big)\nonumber\\\label{eqn:orIIIor}
&=n^{-1/2}\sum_{i=1}^{n}\bm{\xi}_i^{(0)}\epsilon_i+Q_{4n},\;\;\;\; \text(say)
\end{align}
where the $j$th element of $\hat{s}_n^{\dagger(1)}$ is $sgn\big(\hat{\beta}_{j,n}\big)$. Now since $||\hat{\bm{\beta}}_n-\bm{\beta}_n||_{\infty}=O(n^{-1/2})$ on the set $\bm{A}_{3n}$, one can conclude that on the set $\bm{A}_{3n}$, $\mathbf{P_*}\big(\tilde{\bm{s}}_n^{\dagger(1)} = \bm{s}_n^{\dagger(1)}\big)=1$ for sufficiently large $n$. Hence we can conclude that $\mathbf{P}\big(||Q_{4n}||\neq 0\big)=o_p(n^{-1/2})$. Now expanding $\hat{\sigma}_n$ or $\check{\sigma}_n$ around $\sigma$ and by (\ref{eqn:orIIIor}), one has
\begin{align}
\breve{\bm{R}}_n&=\Big[\sigma^{-1}-2^{-1}\sigma^{-3}\Big(n^{-1}\sum_{i=1}^{n}(\epsilon_i^2-\sigma^2)\Big)\Big]\bm{D}_n^{(1)}\bm{C}_{11,n}^{-1}\bm{W}_n^{(1)} +\bm{Q}_{5n}\nonumber\\\label{eqn:orIIIstuor1}
&=\breve{\bm{R}}_{1n}+Q_{5n},\;\;\; \text{(say)}\\
\tilde{\bm{R}}_n&=\Big[\sigma^{-1}-2^{-1}\sigma^{-3}\Big(n^{-1}\sum_{i=1}^{n}(\epsilon_i^2-\sigma^2)\Big)\Big]\bm{\Sigma}_n^{-1/2}\bm{D}_n^{(1)}\bm{C}_{11,n}^{-1}\bm{W}_n^{(1)} +Q_{6n}\nonumber\\\label{eqn:orIIIstuor2}
&=\tilde{\bm{R}}_{2n}+\bm{Q}_{6n},\;\;\; \text{(say)}
\end{align}
where on the set $\bm{A}_{3n}$,
\begin{align*}
\mathbf{P}\big(\|\bm{Q}_{5n}\|=o(n^{-1/2}) \big)=1-o(n^{-1/2})\;\; \&\;\; \mathbf{P}\big(\|\bm{Q}_{6n}\|=o(n^{-1/2}) \big)=1-o(n^{-1/2}).
\end{align*}
Now by looking into (\ref{eqn:stuor1}) in the proof of Theorem \ref{thm:bootI}, we have
\begin{align}\label{eqn:orIIIorEE}
\sup_{\bm{B}\in \mathcal{C}_q}|\mathbf{P}(\breve{\bm{R}}_n\in \bm{B})-\int_{\bm{B}} \bm{\psi}_{1n}(\bm{x})d\bm{x}|=o(n^{-1/2})\nonumber\\
\sup_{\bm{B}\in \mathcal{C}_q}|\mathbf{P}(\tilde{\bm{R}}_n\in \bm{B})-\int_{\bm{B}} \check{\bm{\psi}}_{1n}(\bm{x})d\bm{x}|=o(n^{-1/2})
\end{align}
where
\begin{align*}
\bm{\psi}_{1n}(\bm{x})=&\phi(\bm{x}:\bm{\Sigma}_n)\Bigg[1+\dfrac{1}{\sqrt{n}}\bigg[ -\dfrac{\mu_3}{2\sigma^3}\sum_{|\bm{\alpha}|=1}\bar{\bm{\xi}}_n({\bm{\alpha}})\chi_{\bm{\alpha}}(\bm{x}:\bm{\Sigma}_n)\\
&+\dfrac{\mu_3}{6\sigma^3}\Big{\{}\sum_{|\bm{\alpha}|=3}\bar{\bm{\xi}}_n({\bm{\alpha}})\chi_{\bm{\alpha}}(\bm{x}:\bm{\Sigma}_n)-3\sum_{|\bm{\alpha}|=3}\sum_{|\bm{\zeta}|=1}\bar{\bm{\xi}}_n({\bm{\alpha}})\bar{\bm{\xi}}_n({\bm{\zeta}})\chi_{\bm{\alpha}+\bm{\zeta}}(\bm{x}:\bm{\Sigma}_n)\Big{\}}\bigg]\Bigg],
\end{align*}
and 
\begin{align*}
\check{\bm{\psi}}_{1n}(\bm{x})=&\phi(\bm{x})\Bigg[1+\dfrac{1}{\sqrt{n}}\bigg[ -\dfrac{\mu_3}{2\sigma^3}\sum_{|\bm{\alpha}|=1}\bar{\bm{\xi}}_n^{\dagger(0)}({\bm{\alpha}})\chi_{\bm{\alpha}}(\bm{x})\\
&+\dfrac{\mu_3}{6\sigma^3}\Big{\{}\sum_{|\bm{\alpha}|=3}\bar{\bm{\xi}}^{\dagger(0)}_n({\bm{\alpha}})\chi_{\bm{\alpha}}(\bm{x})-3\sum_{|\bm{\alpha}|=3}\sum_{|\bm{\zeta}|=1}\bar{\bm{\xi}}^{\dagger(0)}_n({\bm{\alpha}})\bar{\bm{\xi}}^{\dagger(0)}_n({\bm{\zeta}})\chi_{\bm{\alpha}+\bm{\zeta}}(\bm{x})\Big{\}}\bigg]\Bigg],
\end{align*}
with $\bar{\bm{\xi}}_n^{(0)}(\bm{\alpha})=n^{-1}\sum_{i=1}^{n}\Big(\bm{\xi}_i^{(0)}\Big)^{\bm{\alpha}}$, $\bar{\bm{\xi}}_n^{\dagger(0)}(\bm{\alpha})=n^{-1}\sum_{i=1}^{n}\Big(\bm{\Sigma}_n^{-1/2}\bm{\xi}_i^{(0)}\Big)^{\bm{\alpha}}$. % $(-D)^{\bm{\alpha}}\phi(\bm{y}:\bm{\Sigma}_n)=\chi_{\bm{\alpha}}(\bm{y}:\bm{\Sigma}_n)\phi(\bm{y}:\bm{\Sigma}_n)$ and $(-D)^{\bm{\alpha}}\phi(\bm{y})=\chi_{\bm{\alpha}}(\bm{y})\phi(\bm{y}).$

Now let us look into bias corrected bootstrapped versions $\breve{\bm{R}}_n^*$ and $\tilde{\bm{R}}_n^*$. Note that similar to (\ref{eqn:orIIIstuor1}) and (\ref{eqn:orIIIstuor2}) we have
\begin{align}
\breve{\bm{R}}_n^*
=& \Big[\hat{\sigma}_n^{-1} -2^{-1}\hat{\sigma}_n^{-3}\big[n^{-1}\sum_{i=1}^{n}(\epsilon_i^{*2}-\hat{\sigma}_n^2)\big]\Big]n^{-1/2}\sum_{i=1}^{n}\bm{\xi}_i^{(0)}\epsilon_i^* +\bm{Q}_{6n}^*\nonumber\\\label{eqn:orIIIstuboot1}
=& \breve{\bm{R}}_{1n}^* +\bm{Q}_{6n}^*\;\;\;\; (say)\\
\tilde{\bm{R}}_n^*
=& n^{-1/2}\sum_{i=1}^{n}\tilde{\bm{\Sigma}}_n^{-1/2}\bm{\xi}_i^{(0)}\hat{\epsilon}_i(G_i^*-\mu_{G^*})\mu_{G^*}^{-1}\nonumber\\
&-2^{-1}\check{\sigma}_n^{-2}\Big[\mu_{G^*}^{-2}n^{-1}\sum_{i=1}^{n}\hat{\epsilon}_i^{2}\big[(G_i^*-\mu_{G^*})^2-\sigma_{G^*}^2\big]\Big]n^{-1/2}\sum_{i=1}^{n}\tilde{\bm{\Sigma}}_n^{-1/2}\bm{\xi}_i^{(0)}\hat{\epsilon}_i(G_i^*-\mu_{G^*})\mu_{G^*}^{-1} +\bm{Q}_{7n}^*\nonumber\\\label{eqn:stuboot2}
=& \tilde{\bm{R}}_{3n}^* +\bm{Q}_{7n}^*\;\;\;\; (say),
\end{align}
where due to conditions (A.2)(ii), (A.3)(ii), (A.5)(i), (A.6)(i), (A.7)(iii) and Lemma \ref{lem:betahat}, we have $\mathbf{P}_*\Big(\|\bm{Q}_{6n}^*\|=o(n^{-1/2})\Big)=1-o_p(n^{-1/2})$ and $\mathbf{P}_*\Big(\|\bm{Q}_{7n}^*\|=o(n^{-1/2})\Big)=1-o_p(n^{-1/2})$.

%Now note that condition (A.6)(ii) imposes Cramer's condition on $\Big((G_1^*-\mu_{G^*}),(G_1^*-\mu{G^*})^2\Big)$ and Lemma \ref{lem:residualCramer} imposes Cramer's condition on $(\epsilon_1^*,\epsilon_1^{*2})$. 
Therefore through the same line of arguments of Theorem \ref{thm:bootI} we have
 \begin{align}
&\sup_{\bm{B}\in \mathcal{C}_q}\Big|\mathbf{P}_*\big(\breve{\bm{R}}_{n}^*\in \bm{B}\big)-\int_{\bm{B}}\bm{\psi}_{1n}^*(\bm{y})d\bm{y}\Big|=o_p\big(n^{-1/2}\big)\;\;\; \text{and}\nonumber\\\label{eqn:orIIIbootEE}
&\sup_{\bm{B}\in \mathcal{C}_q}\Big|\mathbf{P}_*\big(\tilde{\bm{R}}_{n}^*\in \bm{B}\big)-\int_{\bm{B}}\check{\bm{\psi}}_{1n}^{*}(\bm{y})d\bm{y}\Big|=o_p\big(n^{-1/2}\big)
\end{align}
where 
\begin{align*}
\bm{\psi}_{1n}^*(\bm{x})=&\phi(\bm{x}:\bm{\Sigma}_n)\Bigg[1+\dfrac{1}{\sqrt{n}}\bigg[ -\dfrac{\hat{\mu}_3}{2\hat{\sigma}_n^3}\sum_{|\bm{\alpha}|=1}\bar{\bm{\xi}}_n({\bm{\alpha}})\chi_{\bm{\alpha}}(\bm{x}:\bm{\Sigma}_n)\\
&+\dfrac{\hat{\mu}_3}{6\hat{\sigma}_n^3}\Big{\{}\sum_{|\bm{\alpha}|=3}\bar{\bm{\xi}}_n({\bm{\alpha}})\chi_{\bm{\alpha}}(\bm{x}:\bm{\Sigma}_n)-3\sum_{|\bm{\alpha}|=3}\sum_{|\bm{\zeta}|=1}\bar{\bm{\xi}}_n({\bm{\alpha}})\bar{\bm{\xi}}_n({\bm{\zeta}})\chi_{\bm{\alpha}+\bm{\zeta}}(\bm{x}:\bm{\Sigma}_n)\Big{\}}\bigg]\Bigg]
\end{align*}
and
\begin{align*}
\check{\bm{\psi}}_{1n}^*(\bm{x})=&\phi(\bm{x})\Bigg[1+\dfrac{1}{\sqrt{n}}\bigg[ \dfrac{1}{6}\sum_{|\bm{\alpha}|=3}\bar{\bm{\xi}}_n^{*(1)}({\bm{\alpha}})\chi_{\bm{\alpha}}(\bm{x})\\ 
&-\dfrac{1}{2\hat{\sigma}_n^2}\Big{\{}\sum_{|\bm{\alpha}|=1}\bar{\bm{\xi}}_n^{*(3)}({\bm{\alpha}})\chi_{\bm{\alpha}}(\bm{x})
+\sum_{|\bm{\alpha}|=1}\sum_{|\bm{\zeta}|=2}\bar{\bm{\xi}}_n^{*(3)}({\bm{\alpha}})\bar{\bm{\xi}}_n^{*(1)}({\bm{\zeta}})\chi_{\bm{\alpha}+\bm{\zeta}}(\bm{x})\Big{\}}\bigg]\Bigg].
\end{align*}
%with $$(-D)^{\bm{\alpha}}\phi(\bm{y}:\bm{\Sigma}_n)=\chi_{\bm{\alpha}}(\bm{y}:\bm{\Sigma}_n)\phi(\bm{y}:\bm{\Sigma}_n)\;\;\; and$$  $$(-D)^{\bm{\alpha}}\phi(\bm{y})=\chi_{\bm{\alpha}}(\bm{y})\phi(\bm{y}).$$
Therefore Theorem \ref{thm:bootIII} follows by comparing (\ref{eqn:orIIIorEE}) with (\ref{eqn:orIIIbootEE}).

\textbf{Proof of Theorem \ref{thm:symboot}}: Let us denote $\bm{A}_{0n}=\bm{A}_n\cap \bm{B}_n$, $\bm{A}_{2n}$ or $\bm{A}_{3n}$ according as $\hat{\bm{\beta}}_n$ fall in class I, II or III. Similarly define the sets $\bm{A}_{0n}^*$ or $\bm{A}_{0,n}^{**}$ in terms of $*$ or $**$ versions of the sets $\bm{A}_n$, $\bm{B}_n$, $\bm{A}_{2n}$, $\bm{A}_{3n}$. 

First let us look into the residual bootstrap, that is the pair $\big(H_n^r, H_n^{r*}\big)$. Under the conditions of the theorem, it is easy to check that for any choice of $\big(H_n^r, H_n^{r*}\big)$, we have as $n\geq N_1$,
\begin{align*}
H_n^r=\hat{\sigma}_n^{-1} n^{-1/2}\sum_{i=1}^{n}\bm{D}_n^{(1)}C_{11,n}^{-1}\bm{x}_i^{(1)}\epsilon_i\;\; \text{and}\;\; H_n^{r*}=\sigma_n^{*-1} n^{-1/2}\sum_{i=1}^{n}\bm{D}_n^{(1)}C_{11,n}^{-1}\bm{x}_i^{(1)}\epsilon_i^*,
\end{align*}
respectively on the sets $\bm{A}_{0n}$ and $\bm{A}_{0n}^*$, where $N_1$ is some natural number. Using the Lemma \ref{lem:concentration} and the assumptions of the theorem, it is easy to show that $\mathbf{P}\big(\bm{A}_{0,n}\big)=1-o(n^{-2})$ and $\mathbf{P}_*\big(\bm{A}_{0,n}^*\big)=1-o_p(n^{-2})$. Since the desired error rate is $O(n^{-2})$,  without loss of generality we can assume that $\bm{A}_{0n}$ and $\bm{A}_{0n}^*$ both have probability equal to 1. Again note that $H_n^r$ is a studentized version of least square estimator of $\theta_n$ and $H_n^{r*}$ is a suitable residual bootstrap version of it. Therefore the set up for residual bootstrap fits with the regression set up considered in Hall (1988). %By definition of $z_{\alpha}$, $\delta_8<\alpha<1-\delta_8$ implies $0<z_{\alpha}<M_3$ for some positive number $M_3$ which does not depend on $n$. 
Hence following the arguments through (3.4) to (3.7) of Hall (1988) (with error $o(n^{-2})$ instead of $O(n^{-5/2})$), we have
\begin{align*}
&\Big|\mathbf{P}\big(\theta_n\in I_{n,(1-\alpha)}^r\big)- (1-\alpha)\Big|\\
= & \Big|\mathbf{P}\big(|H_{n}^r|\leq \hat{h}_{n,\alpha}^r\big)-\mathbf{P}\big(|H_{n}^r|\leq h_{n,\alpha}^r\big)\Big|\\
\leq &  \sup_{x}\Big|\Big[\mathbf{P}\Big(H_n^r+n^{-3/2}V_{n,\alpha}^r\leq x\Big)
-\mathbf{P}\Big(H_n^r-n^{-3/2}V_{n,\alpha}^r\leq -x\Big)\Big] - \mathbf{P}\Big(|H_n^r|\leq x\Big)\Big|+O(n^{-2})\\
=& O(n^{-2}),
\end{align*}
for any $\alpha \in (0,1)$.

Now let us explore the pair $\big(H_n^p,H_n^{p*}\big)$ in terms of coverage accuracy of the symmetric bootstrap confidence interval. Similar to the residual bootstrap, we have
\begin{align*}
&\Big|\mathbf{P}\big(\theta_n\in I_{n,(1-\alpha)}^p\big)- (1-\alpha)\Big|\\
= & \Big|\mathbf{P}\big(|H_{n}^p|\leq \hat{h}_{n,\alpha}^p\big)-\mathbf{P}\big(|H_{n}^p|\leq h_{n,\alpha}^p\big)\Big|\\
\leq &  \sup_{x}\Big|\Big[\mathbf{P}\Big(H_n^p+n^{-3/2}V_{n,\alpha}^p\leq x\Big)
-\mathbf{P}\Big(H_n^p-n^{-3/2}V_{n,\alpha}^p\leq -x\Big)\Big] - \mathbf{P}\Big(|H_n^p|\leq x\Big)\Big|+O(n^{-2}),
\end{align*}
where $V_{n,\alpha}^p = \sqrt{n}\big\{\hat{q}_2^p(z_{\alpha})-q_2^p(z_{\alpha})\big\}$. $q_2^p(x)\phi(x)$ is the coefficient of $n^{-1}$ in the Edgeworth expansion of $\mathbf{P}\big(H_n^p \leq x\big)$. If $\kappa_j,n^p$ is the $j$th cumulant of $H_n^p$, then from the theory of smooth function models we know $$\kappa_{j,n}^p=O\big(n^{-(j-2)/2}\big),$$
see (1.21) in Bhattacharya and Ghosh (1978). Therefore $n^{-1}q_2^p(x)\phi(x)$ involves cumulants of $H_{n}^p$ upto order 4. Similarly, $n^{-1}\hat{q}_2^p(x)\phi(x)$ involves cumulants of $H_{n}^{p*}$ upto order 4 where $\hat{q}_2^p(x)\phi(x)$ is the coefficient of $n^{-1}$ in the Edgeworth expansion of $H_n^{p*}$. Since $\dfrac{\mathbf{E}(G_1^*-\mu_{G^*})^4}{\mu_{G^*}^4}>1$, $\hat{q}_2^p(\cdot)$ is not in general asymptotically unbiased for $q_2^p(\cdot)$ and hence $V_{n,\alpha}^p$ is not properly centered. See also the discussion on this issue just before the statement of Theorem \ref{thm:symboot} in Section \ref{sec:mainresults}. Suppose $\tilde{q}_2^p(\cdot)$ is asymptotically unbiased for $q_2^p(\cdot)$. Then  $\tilde{V}_{n,\alpha}^p=\sqrt{n}\big\{\tilde{q}_2^p(z_{\alpha})-q_2^p(z_{\alpha})\big\}$ is properly centered and aim should be to replace $V_{n,\alpha}^p$ by $\tilde{V}_{n,\alpha}^p$ in the calculations of coverage error. That is why $\tilde{I}_{n,(1-\alpha)}^p$ is the correct symmetric perturbation bootstrap interval, not $I_{n,(1-\alpha)}^p$. Upon correction, we can show that the coverage error is $O(n^{-2})$ through the same arguments as in case of residual bootstrap. Hence it is enough to find the correction term $C_{n}^p(z_{\alpha})=n^{-1}\big\{\hat{q}_2^p(z_{\alpha})-\tilde{q}_2^p(z_{\alpha})\big\}$. 

As mentioned earlier, only first four cumulants are important. And as the anomaly is only due to the fourth moment of $\dfrac{(G_1^*-\mu_{G^*})}{\mu_{G^*}}$, we are going to find terms of order $n^{-1}$, which involve fourth moment of $\epsilon_1$, in the first four cumulants of $H_{n}^p$. We will do the same for $H_n^{p*}$. In respect to this, we can simply ignore the bias terms appearing in the form of $\bm{T}_n$ and $\bm{T}_n^{**}$ for class II and III, since their contribution is $o(n^{-1})$. For example in case of Lasso, the bias term in $\bm{T}_n$ and $\bm{T}_n^{**}$ contribute in terms of $\dfrac{\lambda_n}{n^2},\dfrac{\lambda_n^2}{n^3},\dots$ which are all $o(n^{-1})$, due to the assumption (A.7)(i)$^{\prime}$. Therefore the correction term $C_{n}^p(z_{\alpha})$ has the same form for each of the three classes. Let us concentrate on class I and assume $\big(H_n^p,H_n^{p*}\big)=\big(\check{\bm{R}}_n,\check{\bm{R}}_n^*\big)$. %For simplicity assume $\mathbf{P}\big(\bm{A}_{0,n}\big)=1$ and $\mathbf{P}_*\big(\bm{A}_{0,n}^{**}\big)=1$.
On the set $\bm{A}_{0,n}$, we have
\begin{align*}
\bm{R}_n = \check{\sigma}_n^{-1}\hat{\Sigma}_n^{-1/2}\bm{T}_n=\check{\sigma}_n^{-1}\hat{\Sigma}_n^{-1/2}n^{-1/2}\sum_{i=1}^{n}\bm{D}_n^{(1)}\bm{C}_{11,n}^{-1}\bm{x}_i^{(1)}\epsilon_i
\end{align*}
where $\bm{T}_n=\sqrt{n}\bm{D}_n\big(\hat{\bm{\beta}}_n-\bm{\beta}_n\big)$, $(\check{\sigma}_n^2-\sigma^2)=n^{-1}\sum_{i=1}^{n}(\hat{\epsilon}_i^2-\epsilon_i)^2+n^{-1}\sum_{i=1}^{n}(\epsilon_i^2-\sigma^2)$ $$\text{and}\;\;\;\check{\sigma}_n^{-1}=\sigma^{-1}-(2\sigma^3)^{-1}(\check{\sigma}_n^2-\sigma^2)+3(8\sigma^5)^{-1}(\check{\sigma}_n^2-\sigma^2)^2 + o(n^{-1}).$$
Therefore we have on the set $\bm{A}_{0,n}$,
\begin{align*}
\check{\bm{R}}_n = & \sigma^{-1}\hat{\Sigma}_n^{-1/2}n^{-1/2}\sum_{i=1}^{n}\bm{D}_n^{(1)}\bm{C}_{11,n}^{-1}\bm{x}_i^{(1)}\epsilon_i\\
& -(2\sigma^3)^{-1}\Big[n^{-1}\sum_{i=1}^{n}(\epsilon_i^2-\sigma^2)\Big]\Big[\sigma_n^{-1/2}n^{-1/2}n^{-1/2}\sum_{i=1}^{n}\bm{D}_n^{(1)}\bm{C}_{11,n}^{-1}\bm{x}_i^{(1)}\epsilon_i\Big]\\
& -(2\sigma^3)^{-1} \Big[n^{-1}\sum_{i=1}^{n}\bm{C}_{11,n}^{-1}\bm{X}_i^{(1)}\epsilon_i\Big]^{\prime}\bm{C}_{11,n}\Big[n^{-1}\sum_{i=1}^{n}\bm{C}_{11,n}^{-1}\bm{X}_i^{(1)}\epsilon_i\Big]\Big[\sigma_n^{-1/2}n^{-1/2}n^{-1/2}\sum_{i=1}^{n}\bm{D}_n^{(1)}\bm{C}_{11,n}^{-1}\bm{x}_i^{(1)}\epsilon_i\Big]\\
&-(2\sigma^3)^{-1}\Big[n^{-1}\sum_{i=1}^{n}\bm{C}_{11,n}^{-1}\bm{X}_i^{(1)}\epsilon_i\Big]^{\prime}\Big[2n^{-1}\sum_{i=1}^{n}\bm{x}_i^{(1)}\epsilon_i\Big]\Big[\sigma_n^{-1/2}n^{-1/2}n^{-1/2}\sum_{i=1}^{n}\bm{D}_n^{(1)}\bm{C}_{11,n}^{-1}\bm{x}_i^{(1)}\epsilon_i\Big]\\
&+3(8\sigma^5)^{-1}\Big[n^{-1}\sum_{i=1}^{n}(\epsilon_i^2-\sigma^2)\Big]^2\Big[\sigma_n^{-1/2}n^{-1/2}n^{-1/2}\sum_{i=1}^{n}\bm{D}_n^{(1)}\bm{C}_{11,n}^{-1}\bm{x}_i^{(1)}\epsilon_i\Big]\\
&+o(n^{-1})\\
= & U_n + o(n^{-1})\;\;\; \text{(say)}.
\end{align*}
Hence enough to look into first four cumulants of $U_n$ instead of $\bm{R}_n$. Let us denote by $\tilde{\mu}_{j,n}$ (or $\tilde{\kappa}_{j,n}$) the collection of terms of order $n^{-1}$ which involve fourth moment of $\epsilon_1$, in $j$th raw moment (or cumulant) of $U_n$, $j=1,2,3,4$. Now writing $U_n=a+b+c+d+e$, it can be shown that 
\begin{align*}
&\tilde{\mu}_{1,n} = 0\\
&\tilde{\mu}_{2,n} = \mathbf{E}(b^2 + 2ab + 2ae) = 0\\
&\tilde{\mu}_{3,n} = 0\\
& \tilde{\mu}_{4,n} = \mathbf{E}(a^4 + 4a^3b+4a^3e+6a^2b^2)
\end{align*}
and hence
\begin{align*}
\tilde{\kappa}_{1,n} = \tilde{\mu}_{1,n} &= 0\\
\tilde{\kappa}_{2,n}  = \tilde{\mu}_{2,n} &= 0\\
\tilde{\kappa}_{3,n} = \tilde{\mu}_{3,n} &= 0\\
\tilde{\kappa}_{4,n} = \tilde{\mu}_{4,n} -3 \tilde{\mu}_{2,n}^2 &= \mathbf{E}(a^4 + 4a^3b+4a^3e+6a^2b^2) - 6\mathbf{E}(b^2+2ab+2ae)\\
 & = n^{-1}\Big[\sigma^{-4}\Sigma_n^{-2}n^{-1}\sum_{i=1}^{n}\big(\bm{D}_n^{(1)}\bm{C}_{11,n}^{-1}\bm{x}_i^{(1)}\big)^4\mathbf{E}\epsilon_1^4\Big] + n^{-1}\Big[\dfrac{\mathbf{E}\epsilon_1^4}{\sigma^4}-1\Big]
\end{align*}
Now similar to $\check{\bm{R}}_n$, it can be shown $\check{\bm{R}}_n^*=U_n^*+o_p(n^{-1})$. Define $\tilde{\kappa}_{j,n}^*$, $j=1,2,3,4$, such that $\tilde{\kappa}_{j,n}^*$ to $U_n^*$ is what $\tilde{\kappa}_{j,n}$ is to $U_n$. Then in the same fashion it can be shown that  
\begin{align*}
\tilde{\kappa}_{1,n}^*  = 0,\;\;\; \tilde{\kappa}_{3,n}^*  = 0,
\end{align*}
\begin{align*}
\tilde{\kappa}_{2,n}^*  = n^{-1}\bigg[\Big[\check{\sigma}_n^{-4}n^{-1}\sum_{i=1}^{n}\hat{\epsilon}_i^4\Big] - \Big[\check{\sigma}_n^{-2}\tilde{\Sigma}_n^{-1}n^{-1}\sum_{i=1}^{n}\Big(\bm{D}_n^{(1)}\bm{C}_{11,n}^{-1}\bm{x}_i^{(1)}\Big)^2\hat{\epsilon}_i^4\Big]\bigg]\Big[\dfrac{\mathbf{E}(G_1^*-\mu_{G^*})^4}{\mu_{G^*}^4}-1\Big],
\end{align*}
\begin{align*}
\tilde{\kappa}_{4,n}^* 
 =& n^{-1}\Big[\tilde{\Sigma}_n^{-2}n^{-1}\sum_{i=1}^{n}\Big(\bm{D}_n^{(1)}\bm{C}_{11,n}^{-1}\bm{x}_i^{(1)}\Big)^4\hat{\epsilon}_i^4\Big]\dfrac{\mathbf{E}(G_1^*-\mu_{G^*})^4}{\mu_{G^*}^4}\\
 &+4n^{-1}\Big[\check{\sigma}_n^{-2}\tilde{\Sigma}_n^{-1}n^{-1}\sum_{i=1}^{n}\Big(\bm{D}_n^{(1)}\bm{C}_{11,n}^{-1}\bm{x}_i^{(1)}\Big)^2\hat{\epsilon}_i^4\Big]\Big[\dfrac{\mathbf{E}(G_1^*-\mu_{G^*})^4}{\mu_{G^*}^4}-1\Big]\\
 &-3n^{-1}\Big[\check{\sigma}_n^{-4}n^{-1}\sum_{i=1}^{n}\hat{\epsilon}_i^4\Big]\Big[\dfrac{\mathbf{E}(G_1^*-\mu_{G^*})^4}{\mu_{G^*}^4}-1\Big]
\end{align*}
Therefore the contribution of $\tilde{\kappa}_{j,n}$'s, $j=1,2,3,4$, in the characteristic function $\mathbf{E}\big(e^{it\check{\bm{R}}_n}\big)$ of $\check{\bm{R}}_n$ is $$\Big[\frac{1}{2!}\tilde{\kappa}_{2,n}(it)^2+\frac{1}{4!}\tilde{\kappa}_{4,n}(it)^4\Big]e^{-t^2/2}.$$
The Fourier inversion of this contribution is $$l_2(x)=\Big[\frac{1}{2!}\tilde{\kappa}_{2,n}(-1)^2d^2/dx^2+\frac{1}{4!}\tilde{\kappa}_{4,n}(-1)^4d^4/dx^4\Big]\phi(x).$$ For the time being assume that $n^{-1}q_2^p(x)\phi(x)$ appears in the Edgeworth expansion of $\check{\bm{R}}_n$ only due to $l_2(x)$. We can assume this abuse of the notation of $q_2^p$ since we need to correct only the effect of the fourth moment. Then we have 
\begin{align*}
q_2^p(x)\phi(x) = n\Big[\frac{1}{2!}\tilde{\kappa}_{2,n}(-d/dx)+\frac{1}{4!}\tilde{\kappa}_{4,n}(-1)^3d^3/dx^3\Big]\phi(x),
\end{align*}
that is $q_2^p(x)=-nx\Big[\frac{1}{2!}\tilde{\kappa}_{2,n}+\frac{1}{4!}\tilde{\kappa}_{4,n}(x^2-3)\Big]=-\frac{1}{4!}n\tilde{\kappa}_{4,n}x(x^2-3)$. Similarly abusing the notation $\hat{q}_2^p$, we have
$$\hat{q}_2^p(x)=-nx\Big[\frac{1}{2!}\tilde{\kappa}_{2,n}^*+\frac{1}{4!}\tilde{\kappa}_{4,n}^*(x^2-3)\Big].$$
We can define $\tilde{q}_2^p(\cdot)$ as
$$\tilde{q}_2^p(x)=-nx\Big[\frac{1}{2!}\hat{\kappa}_{2,n}^*+\frac{1}{4!}\hat{\kappa}_{4,n}^*(x^2-3)\Big],$$
where 
\begin{align*}
\hat{\kappa}_{2,n}^*  = n^{-1}\bigg[\Big[\check{\sigma}_n^{-4}n^{-1}\sum_{i=1}^{n}\hat{\epsilon}_i^4\Big] - \Big[\check{\sigma}_n^{-2}\tilde{\Sigma}_n^{-1}n^{-1}\sum_{i=1}^{n}\Big(\bm{D}_n^{(1)}\bm{C}_{11,n}^{-1}\bm{x}_i^{(1)}\Big)^2\hat{\epsilon}_i^4\Big]\bigg]\;\; \text{and}
\end{align*}
\begin{align*}
\hat{\kappa}_{4,n}^* 
 =& n^{-1}\Big[\tilde{\Sigma}_n^{-2}n^{-1}\sum_{i=1}^{n}\Big(\bm{D}_n^{(1)}\bm{C}_{11,n}^{-1}\bm{x}_i^{(1)}\Big)^4\hat{\epsilon}_i^4\Big]\\
 &+4n^{-1}\Big[\check{\sigma}_n^{-2}\tilde{\Sigma}_n^{-1}n^{-1}\sum_{i=1}^{n}\Big(\bm{D}_n^{(1)}\bm{C}_{11,n}^{-1}\bm{x}_i^{(1)}\Big)^2\hat{\epsilon}_i^4\Big]\\
 &-3n^{-1}\Big[\check{\sigma}_n^{-4}n^{-1}\sum_{i=1}^{n}\hat{\epsilon}_i^4\Big] - n^{-1}
\end{align*}
Note that $\check{\sigma}_n^2\rightarrow \sigma^2$, by Lemma \ref{lem:Sigma}, $\tilde{\bm{\Sigma}}_n\rightarrow \sigma^2 \bm{\Sigma}_n$ as $n\rightarrow \infty$. Again $|\hat{\epsilon}_i^2-\epsilon_i^2|\leq \|\bm{x}_i^{(1)}\|\|(\hat{\bm{\beta}}_n-\bm{\beta}_n)\|$ where by Lemma \ref{lem:betahat} or Lemma \ref{lem:betahathat} $\|(\hat{\bm{\beta}}_n-\bm{\beta}_n)\|=o_p(1)$. Hence $\tilde{q}_2^p(\cdot)$ is asymptotically unbiased for $q_2^p(\cdot)$. Therefore the correction needed for constructing symmetric confidence interval is
\begin{align*}
C_{n}^p(z_{\alpha})&=n^{-1}\big\{\hat{q}_2^p(z_{\alpha})-\tilde{q}_2^p(z_{\alpha})\big\}\\
&=-x\Big[\frac{1}{2!}\big(\tilde{\kappa}_{2,n}^*-\hat{\kappa}_{2,n}^*\big)+\frac{1}{4!}\big(\tilde{\kappa}_{4,n}^*-\hat{\kappa}_{4,n}^*\big)(x^2-3)\Big]\\
&=-n^{-1}x\Big[\dfrac{\omega_2}{2}+\dfrac{\omega_4}{24}(x^2-3)\Big],
\end{align*}
where
\begin{align*}
\omega_2 = \bigg[\check{\sigma}_n^{-4}\Big[n^{-1}\sum_{i=1}^{n}\hat{\epsilon}_i^4\Big] - \check{\sigma}_n^{-2}\tilde{\Sigma}_n^{-1}\Big[n^{-1}\sum_{i=1}^{n}\big\{\hat{D}_n^{(1)}\hat{C}_{11,n}^{-1}\bm{x}_i^{(1)}\big\}^2\hat{\epsilon}_i^4\Big]\bigg]\bigg[\dfrac{E(G_1^*-\mu_{G^*})^4}{\mu_{G^*}^4}-2\bigg],
\end{align*}
\begin{align*}
\omega_4 =& \tilde{\Sigma}_n^{-2}\Big[n^{-1}\sum_{i=1}^{n}\big\{\hat{D}_n^{(1)}\hat{C}_{11,n}^{-1}\bm{x}_i^{(1)}\big\}^4\hat{\epsilon}_i^4\Big]\bigg[\dfrac{E(G_1^*-\mu_{G^*})^4}{\mu_{G^*}^4}-1\bigg]\\
& + 4\check{\sigma}_n^{-2}\tilde{\Sigma}_n^{-1}\Big[n^{-1}\sum_{i=1}^{n}\big\{\hat{D}_n^{(1)}\hat{C}_{11,n}^{-1}\bm{x}_i^{(1)}\big\}^2\hat{\epsilon}_i^4\Big]\bigg[\dfrac{E(G_1^*-\mu_{G^*})^4}{\mu_{G^*}^4}-2\bigg]\\
& - 3\check{\sigma}_n^{-4}\Big[n^{-1}\sum_{i=1}^{n}\hat{\epsilon}_i^4\Big]\bigg[\dfrac{E(G_1^*-\mu_{G^*})^4}{\mu_{G^*}^4}-2\bigg]\\
& + 1.
\end{align*}

Hence we are done.

%\section{Conclusion}

%


\end{document}